\documentclass{elsarticle}

\usepackage{jcomp}
\usepackage{framed,multirow,amsmath}
\usepackage{amssymb,cases}
\usepackage{latexsym}
\usepackage{amsthm}
\usepackage{algorithm}
\usepackage{algorithmic}
\usepackage{url}
\usepackage{xcolor}
\definecolor{newcolor}{rgb}{.2,.3,.7}

\usepackage{hyperref}
\hypersetup{
	colorlinks,
	citecolor=newcolor,
	filecolor=newcolor,
	linkcolor=newcolor,
	urlcolor=newcolor
}

\theoremstyle{plain}
\theoremstyle{definition}

\newtheorem{assumption}{Assumption}
\newtheorem{remark}{Remark}
\newtheorem{example}{Example}[section]

\ifpdf
\DeclareGraphicsExtensions{.eps,.pdf, .png,.jpg}
\else
\DeclareGraphicsExtensions{.eps}
\fi

\usepackage{bm}

\makeatletter
\@addtoreset{equation}{section}
\@addtoreset{figure}{section}
\@addtoreset{table}{section}
\makeatother

\usepackage{tikz,pgfplots,pgf}
\usetikzlibrary{matrix,shapes,arrows,positioning}
\usepackage{etoolbox} 
\usetikzlibrary{shapes}
\tikzset{box/.style ={
	rectangle, 
	rounded corners =5pt, 
	minimum width =50pt, 
	minimum height =20pt, 
	inner sep=5pt, 
	draw=blue} 
}

\tikzset{zbox/.style ={
		rectangle, 
		minimum width =50pt, 
		minimum height =20pt, 
		inner sep=5pt, 
		draw=black} 
}

\tikzset{ball/.style ={
		circle, 
		minimum width =20pt, 
		minimum height =20pt, 
		inner sep=0.1pt, 
		draw=blue} 
}

\tikzset{global scale/.style={
		scale=#1,
		every node/.append style={scale=#1}
	}
}

\tikzset{elegant/.style={smooth,thick,samples=50,cyan}}
\tikzset{eaxis/.style={->,>=stealth}}

\graphicspath{{pics/}}
\begin{document}

\verso{Given-name Surname \textit{etal}}

\begin{frontmatter}

\tnotetext[tnote1]{Date: \today}
\title{\textbf{\large Prox-PINNs: A Deep Learning Algorithmic Framework for Elliptic Variational Inequalities}}

\author{\snm{{\small Yu {Gao}}}\footnote{Department of Mathematics, The Hong Kong University of Science and Technology, Clear Water Bay, Hong Kong, China. Email: {\color{newcolor}yugaomath@gmail.com}}}
\author{\snm{{\small Yongcun {Song}}}\footnote{Division of Mathematical Sciences, School of Physical and Mathematical Sciences,
Nanyang Technological University, 21 Nanyang Link, 637371, Singapore. Email: {\color{newcolor}yongcun.song@ntu.edu.sg}}}
\author{\snm{{\small Zhiyu {Tan}}}\footnote{School of Mathematical Sciences and Fujian Provincial Key Laboratory on Mathematical Modeling and High Performance Scientific Computing, Xiamen University, Xiamen 361005, Fujian, China. Email: {\color{newcolor}zhiyutan@xmu.edu.cn}}}
\author{\snm{{\small Hangrui {Yue}}}\corref{cor1}\footnote{School of Mathematical Sciences, Nankai University, Tianjin 300071, China. Email: {\color{newcolor}yuehangrui@gmail.com}}}
\cortext[cor1]{Corresponding author}
\author{\snm{{\small Shangzhi {Zeng}}}\footnote{National Center for Applied Mathematics Shenzhen \& Department of Mathematics, Southern University of Science and Technology, Shenzhen 518055, Guangdong, China.  Email: {\color{newcolor}zengsz@sustech.edu.cn }}}

\begin{abstract}
{\small
\textbf{Abstract:} Elliptic variational inequalities (EVIs) present significant challenges in numerical computation due to their inherent non-smoothness, nonlinearity, and inequality formulations. Traditional mesh-based methods often struggle with complex geometries and high computational costs, while existing deep learning approaches lack generality for diverse EVIs. To alleviate these issues, this paper introduces Prox-PINNs, a novel deep learning algorithmic framework that integrates proximal operators with physics-informed neural networks (PINNs) to solve a broad class of EVIs. The Prox-PINNs reformulate EVIs as nonlinear equations using proximal operators and then approximate the solutions via neural networks that enforce boundary conditions as hard constraints. Then the neural networks are trained by minimizing physics-informed residuals. The Prox-PINNs framework advances the state-of-the-art by unifying the treatment of diverse EVIs within a  mesh-free and scalable computational architecture. The framework is demonstrated on several prototypical applications, including obstacle problems, elasto-plastic torsion, Bingham visco-plastic flows, and simplified friction problems. Numerical experiments validate the method’s accuracy, efficiency, robustness, and flexibility across benchmark examples. 



}
\end{abstract}


\end{frontmatter}



\section{Introduction}

Elliptic variational inequalities (EVIs) constitute a fundamental class of nonlinear problems arising in diverse applications, including contact mechanics, non-Newtonian fluid flows, elasto-plastic deformation, and image processing, where traditional equality-based formulations fail to capture realistic constraints, see e.g., \cite{duvant1976inequalities,Glowinski1984Numerical,glowinski2008lectures,glowinski2015variational,glowinskinumerical1981,kinderlehrer1980introduction,lions1976some,tremolieres2011numerical}. In particular, concrete real applications modeled by EVIs include obstacle problems, elasto-plastic torsion problems, simplified friction problems,  image restoration problems,  image denoising problems, simplified Signorini problems,  Bingham visco-plastic flows, and optimal control of partial differential equations (PDEs), see \cite{Glowinski1984Numerical,glowinski2008lectures,song2019inexact,tremolieres2011numerical}  and references therein.  EVIs typically result in non-smooth or discontinuous solutions, necessitating sophisticated mathematical tools for theoretical analysis. Notable theoretical advancements for EVIs, such as existence, uniqueness, and regularity of solutions, can be referred to \cite{Glowinski1984Numerical,glowinski2008lectures,glowinski2015variational} and references therein. Despite the theoretical advances, solving EVIs numerically is challenging due to the presence of non-smoothness and nonlinearity and the need to resolve possible free boundaries and large-scale systems. Addressing these challenges requires a combination of advanced optimization techniques, tailored discretization strategies, and robust iterative algorithms, making EVIs significantly more demanding than standard elliptic PDEs.
Therefore, algorithmic design for EVIs requires systematic approaches that carefully integrate their inherent structures and characteristics.

Mathematically, EVIs can be formulated as
\begin{equation}\label{EVI}
\text{Find}~u\in V, \quad\text{such that}\quad a(u,v-u)+j(v)-j(u)\geq l(v-u),\quad \forall v\in V.
\end{equation}
In (\ref{EVI}), $V$ is a real Hilbert space defined over a domain $\Omega\subset\mathbb{R}^d (d\geq 1)$, endowed with the inner product $(\cdot,\cdot)$ and the norm $\|\cdot\|$.
 The bilinear functional $a: V\times V\rightarrow \mathbb{R}$ is  continuous and $V$-elliptic, that is, there exist constants $c_1>0$ and $c_2>0$ such that $|a(w,v)|\leq c_1\|w\|\cdot\|v\|$ and $|a(v,v)|\geq c_2\|v\|^2$, $\forall w,v \in V$.  
The functional $l\in V'$ with $V'$ the dual space of $V$, and  the functional $j: V\rightarrow \mathbb{R}\cup\{+\infty\}$ is non-smooth, convex, proper, and lower semi-continuous. Note that $a(\cdot,\cdot)$ is not necessarily symmetric.  As commented in \cite{glowinski2015variational}, if $a(\cdot,\cdot)$ is symmetric, the EVI (\ref{EVI}) is equivalent to an optimization problem, which makes  (\ref{EVI}) easier to solve. 

Typically, problems in the form of (\ref{EVI}) are called EVIs of the second kind (EVI.2). If we consider a closed convex nonempty subset $K\subset V$ and let $j$ be the indicator functional of $K$, then (\ref{EVI}) reduces to
\begin{equation}\label{EVI.1}
\text{Find}~u\in K,\quad\text{such that}\quad a(u,v-u)\geq l(v-u),\quad \forall v\in K,
\end{equation}
which is called EVIs of the first kind (EVI.1). Actually, as commented in \cite{glowinski2015variational}, the distinction between (\ref{EVI}) and (\ref{EVI.1}) is rather artificial, since (\ref{EVI.1}) can be viewed as a special case of (\ref{EVI}). Therefore, we focus on (\ref{EVI}) hereafter and all the results can be applied to (\ref{EVI.1}) directly.

Over the years, the design and analysis of numerical methods for EVIs have been intensively studied in the literature. In particular, some
numerical approaches have been designed for solving some specific cases of \eqref{EVI} with a primary focus on developing iterative schemes that can overcome the difficulty of the
nonsmoothness of $j$. 
For instance, over-relaxation methods were studied for obstacle problems and simplified Signorini problems in \cite{Glowinski1984Numerical}, augmented Lagrangian methods and alternating direction method of multipliers (ADMM) were applied to solve obstacle problems in \cite{Glowinski1984Numerical} and to solve Bingham viscous-plastic fluid flows in \cite{Dean2007On}. Newton-type methods were considered in \cite{de2010numerical,los2011duality} for EVI.2 with applications to Bingham visco-plastic flows, simplified friction problems, and total variation regularization in image processing. Semismooth Newton and augmented Lagrangian methods were studied in \cite{stadler2004semismooth} for a simplified friction problem. Several Moreau-Yosida regularization-based path-following methods for a class of gradient-constrained EVIs were proposed in \cite{hintermuller2015several}.  Moreover,  the $L^1$-penalty method \cite{Tran2014An}, the primal-dual method \cite{zosso2017efficient}, and the operator-splitting method \cite{liu2023fast} were designed for obstacle problems. 
A preconditioned conjugate gradient-based inexact Uzawa method was discussed in \cite{Cheng2003Inexact} for EVI.2.
Note that all these methods are implemented with mesh-based discretization schemes (e.g., finite difference methods (FDM) or finite element methods (FEM)). As a result, these methods are struggling to solve problems in complex domains and high-dimensional spaces. Moreover, large-scale and ill-conditioned algebraic systems are usually required to be solved at each iteration, leading to a high computational burden. 

To alleviate the above-mentioned issues, some deep learning methods have been recently designed for EVIs in the literature, see \cite{alphonse2024neural,bahja2023physics,cao2024hybrid,cheng2023deep,combettes2020deep,darehmiraki2022deep,huang2022deep,huang2020int,schwab2022deep,zhao2022two}.
Compared with traditional iterative methods, which discretize the EVIs using mesh-based schemes, these deep learning methods are usually mesh-free, easy to implement, and effective in solving problems in complex domains and high-dimensional spaces. Moreover, deep learning methods avoid solving algebraic systems completely by taking advantage of automatic differentiation, and could break the curse of dimensionality; computational costs can thus be reduced. In particular, once the neural networks are trained on a fixed set of randomly sampled points, deep learning methods can solve the problem at a new resolution by simply performing a forward pass of the pre-trained networks. In contrast, traditional numerical methods have to recompute the solution from scratch for each resolution, resulting in substantially higher computational costs.  Despite these advantages, it is worth noting that the above-mentioned deep learning methods only apply to some specific cases of  \eqref{EVI} and lack the flexibility to build a general framework for seamlessly tackling various EVIs. For instance, the methods in \cite{bahja2023physics,cheng2023deep,darehmiraki2022deep,zhao2022two} are designed for obstacle problems, and only EVI.1 was considered in \cite{alphonse2024neural}. In \cite{cao2024hybrid,huang2020int}, several deep learning methods were proposed for EVIs in the form of \eqref{EVI} but with symmetric $a(\cdot,\cdot)$ , which limits their applicability domain. To the best of our knowledge, there
seems to be no deep learning approach in the literature that can address the generic EVI model \eqref{EVI} without imposing restrictive constraints. 

In this paper, we develop a novel deep learning algorithmic framework that is capable of  solving a general class of EVIs modeled by \eqref{EVI}.  To this end, we first rewrite \eqref{EVI} as a nonlinear equation by leveraging the proximal operator of $j$.  We then approximate the solution $u$ by a neural network, where the boundary condition of $u$ (e.g., the homogeneous Dirichlet boundary condition when $V=H_0^1(\Omega)$) is imposed as a hard constraint and can be treated separately in the training process.  Inspired by physics-informed neural networks (PINNs) \cite{raissi2019physics}, the residual of the nonlinear equation is used as the loss function to train the neural network. Therefore, the framework is termed Prox-PINNs to signify the integration of the proximal operator and the physics-informed nature originating from PINNs.   Note that while PINNs have been widely applied across scientific domains (see \cite{cuomo2022scientific,faroughi2024physics,lai2023hard,song2024admm,toscano2025pinns} and the references therein), they lack the inherent capability to address inequalities like \eqref{EVI}.  The Prox-PINNs thus substantially extend the applicability of PINNs to EVIs while retaining their advantages: being mesh-free, easy to implement, and adaptable to diverse scenarios.

The Prox-PINNs is a high-level framework that imposes no specific constraints on $a(\cdot,\cdot)$ or $j(\cdot)$ and hence can be applied to various EVIs in the form of \eqref{EVI}, which is distinguished from the existing deep learning methods in the literature.
We demonstrate the numerical implementation of the Prox-PINNs via case studies involving distinct choices of $j$. The framework is then applied to several prototypical EVIs with different $a(\cdot,\cdot)$ and $j(\cdot)$, including obstacle problems, elasto-plastic torsion problems, Bingham visco-plastic flows, and simplified friction problems. For each EVI, we validate the effectiveness, efficiency, accuracy, robustness, and flexibility of the Prox-PINNs through extensive numerical experiments on benchmark examples. To highlight
the advances of the Prox-PINNs, we include some numerical
comparisons with FEM-based high-fidelity traditional numerical methods and other deep learning methods.

The rest of this paper is organized as follows. In Section \ref{se: algorithmic_framework} we present the Prox-PINNs framework.  Then, we elaborate on the numerical implementation of the algorithmic framework through case studies in Section \ref{se:case_studies}, and specific Prox-PINNs methods are derived for different types of EVIs.  In Section \ref{se: numerical},  the effectiveness and efficiency
of the resulting Prox-PINNs methods are demonstrated by extensive numerical studies for several typical EVIs. Finally, some conclusions and perspectives are given in Section \ref{se:conclusion}.

\section{The Prox-PINNs Framework for (\ref{EVI})}\label{se: algorithmic_framework}

In this section, we present the proposed Prox-PINNs framework for \eqref{EVI}. 
Note that machine learning algorithms such as PINNs primarily operate on the strong form of partial differential equations. We shall show that, for EVIs, an analogous strong form can be formulated, enabling their solution via PINNs.
More precisely, we first employ the proximal operator of  $j$ to reformulate \eqref{EVI} as a nonlinear equation in terms of $u$. The solution $u$ is then approximated by constructing a neural network surrogate, which is trained within a physics-informed framework through the minimization of a loss function that encodes the governing equation.

\subsection{Proximal Formulation of (\ref{EVI})}
Since $a(\cdot,\cdot)$ is a bilinear form on $V\times V$, by Riesz representation theorem, there exists $A\in \mathcal{L}(V, V')$ such that
$\langle Au, v\rangle_{V',V} = a(u,v), \forall\ u, v\in V.$
Therefore, the problem (\ref{EVI}) can be rewritten as
\begin{equation*}
	u\in V, \quad\text{such that}\quad  \langle Au-l, v-u\rangle_{V',V} + j(v)-j(u)\geq 0,\quad \forall v\in V,
\end{equation*}
which implies that
\begin{equation}\label{eq:inclusion_o}
u\in V,\quad -(Au-l)\in\partial j(u)\quad \mbox{in}\ \ V',
\end{equation}
where $\partial j(u):=\{\xi\in V' \mid j(v)-j(u) \geq  \langle \xi, v-u \rangle_{V',V}, \forall\ v\in V\}$ is the subdifferential of $j$ at $u$. 

Let $H$ be another Hilbert space with $V$ continuously embedded into $H$ and it satisfies the Gelfand triple $V\subset H=H' \subset V'$. 
\begin{assumption}\label{assum:1}
The following assumptions hold:
\begin{enumerate}[(1)]
\item $j(\cdot)$ can be extended to $H$ as a convex, proper and lower semi-continuous functional. 
\item  $Au, l\in H$.
\end{enumerate}
\end{assumption}
In the following arguments, we suppose that Assumption \ref{assum:1} holds and hence there exists $f\in H$ such that $l(v)=(f,v)_H, \forall v\in H$. Therefore,
we have
\begin{equation}\label{eq:inclusion}
u\in V\subset H,\quad -(Au-f)\in\partial j(u)\quad \mbox{in}\ \ H,
\end{equation}
where $\partial j(u):=\{\xi\in H \mid j(v)-j(u) \geq  (\xi, v-u )_H, \forall\ v\in H\}$ is the subdifferential of $j$ at $u$. 
Note that $H$ can be chosen as $V$.

Let $\eta\in \mathbb{R}$ be a positive constant and we rewrite \eqref{eq:inclusion} as
\begin{equation*}
u\in V\subset H,\quad u-\eta(Au-f)\in u+\eta\partial j(u).
\end{equation*}
We thus have 
\begin{equation}\label{eq:oc_evi}
u\in V\subset H,\quad (I-\eta A)u+\eta f\in(I+\eta\partial j)(u).
\end{equation}

Since $j$ is convex, proper, and lower semi-continuous, $\partial j$ is maximal monotone and hence, the operator $(I+\eta \partial j)^{-1}$ is single-valued (see e.g., \cite{{bauschkeconvex}}). Let $w=(I-\eta A)u+\eta f$, it follows from \eqref{eq:oc_evi} that
$$
0\in \partial j(u)+\frac{1}{\eta}(u-w),
$$
or equivalently
$$
z=\mathop{\arg\min}\limits_{v\in H} j(v)+\frac{1}{2\eta}\|v-w\|_{H}^2:=\text{Prox}_{\eta j}(w),
$$
where $\text{Prox}_{\eta j}(\cdot)$ is the proximal operator of $j$. 
The above result indicates that (\ref{eq:oc_evi}) can be written as
\begin{equation}\label{eq:oc}
u\in V\subset H,\quad	\text{Prox}_{\eta j}((I-\eta A)u+\eta f)=u.
\end{equation}
We thus obtain a nonlinear equation that can be treated as a  strong form associated with the underlying EVI.

\subsection{A Concrete Illustrative Example}\label{se:obstacle}
To provide a concrete illustration of the preceding discussion, we specify \eqref{EVI} as an obstacle problem. For this purpose, we let $V=H_0^1(\Omega)$ with $\Omega\subset\mathbb{R}^d (d\geq 1)$ a bounded domain, $j(\cdot) = I_K(\cdot)$ with $K: = \{v\in H^1_0(\Omega)\mid v\geq \psi,\ \ \mbox{a.e.\ in}\ \Omega\}$ and specify the bilinear functional $a: H_0^1(\Omega)\times H_0^1(\Omega)\rightarrow \mathbb{R}$ as
$$
a(u,v)=\int_\Omega \nabla u\cdot\nabla v ~dx~\text{and}~l(v) = \int_\Omega fv~dx
$$
with $f\in L^2(\Omega)$. 
	
We define the operator $A: H_0^1(\Omega)\rightarrow H^{-1}(\Omega)$ by 
$$Av= - \Delta v, \quad\forall v\in H_0^1(\Omega).$$
The problem (\ref{EVI}) can be rewritten as
\begin{equation*}
	u\in H_0^1(\Omega), \quad\text{such that}\quad \langle - \Delta u-f, v-u\rangle_{H^{-1}(\Omega), H_0^1(\Omega)} + j(v) - j(u) \geq 0,\quad \forall v\in H_0^1(\Omega),
\end{equation*}
which implies that
\begin{equation}\nonumber
u\in H_0^1(\Omega),\quad  \Delta u + f \in\partial j(u)\quad \mbox{in}\ \ H^{-1}(\Omega),
\end{equation}
where $\partial j(u):=\{\xi\in H^{1}_0(\Omega) \mid j(v)-j(u) \geq  \langle \xi, v-u \rangle_{H^{-1}(\Omega), H_0^1(\Omega)}, \forall v\in H_0^1(\Omega)\}$ is the subdifferential of $j$ at $u$. 


Under some well-known additional regularity assumptions on $\Omega$, $f$, and $\psi$, we have $u\in H^2(\Omega)$ and $\Delta u + f\in L^2(\Omega)$, see e.g., \cite{grisvard1985ellptic, Ito1990augmented,ladyzhenskaya1985boundary}. Inspired by these results, we can choose $H = L^2(\Omega)$ and redefine $j(\cdot) = I_K(\cdot)$ with $K = \{ v\in H\mid \ v\geq \psi,\ \ \mbox{a.e.\ in}\ \Omega\}$. Then, we have
$$\text{Prox}_{\eta j}(w) = \max\{w,  \psi\}, ~\text{with}~w=(I+\eta \Delta)u+\eta f\in L^2(\Omega).$$
As a result, the equation \eqref{eq:oc} can be specified as
\begin{equation}\label{eq:VI_SF_exa}
\left\{
\begin{aligned}
u & = \max\{(I + \eta \Delta)u+\eta f, \psi\}\quad &&\mbox{in}\ \ \Omega;\\
u & = 0 &&\mbox{on}\ \ \partial\Omega,
\end{aligned}
\right.
\end{equation} 
which is a strong form of the EVI under consideration and coincides with the result presented in \cite{Ito1990augmented}.

\begin{remark}
For the current example, the proximal operator $\text{Prox}_{\eta j}(\cdot)$ in \eqref{eq:oc}, if defined on $H_0^1(\Omega)$,  lacks a closed-form expression. Consequently, obtaining a strong form of the underlying EVI is analytically intractable.
\end{remark}

\begin{remark}
In general cases, following \cite{glowinskinumerical1981, gong2025new,kinderlehrer1980introduction}, we introduce
$$ \Omega_0 = \{x\in \Omega\mid\ u(x) = \psi(x)\}\ \ \mbox{and}\ \ \Omega^+ = \{x\in \Omega\mid\ u(x) > \psi(x)\},$$
and have
$$-\Delta u - f = 0\ \ \mbox{a.e.\ in}\ \ \Omega^+,\ \  u= \psi\ \ \mbox{a.e.\ in}\ \ \Omega_0.$$
If $f\in L^2(\Omega)$ and $\psi\in H^2(\Omega)$, the interior regularity theory of second order elliptic partial differential equations (see e.g. \cite{evans2010PDEs}) gives that $u$ is in $H^2$ away from the set $\partial \Omega^+$. Therefore, the strong form (\ref{eq:VI_SF_exa}) can be treated as a regularization of (\ref{EVI}). 
\end{remark}

\subsection{The Prox-PINNs Framework}
Next, we elaborate on a neural network approach for solving the nonlinear equation \eqref{eq:oc} and present the Prox-PINNs framework for (\ref{EVI}). To fix ideas, we consider the homogeneous Dirichlet boundary condition $u=0$ on $\partial \Omega$ and construct a neural network $\hat{u}(x;\bm{\theta}_u)$ verifying $\hat{u}(x;\bm{\theta}_u)=0$ for $x\in \partial\Omega$ to approximate $u$.  To this end, we first introduce a function $h: \overline{\Omega} \to \mathbb{R}$ satisfying
\begin{equation*}
	h \in C(\overline{\Omega}), \quad h(x) = 0 \text{~if~and~only~if~} x \in \partial \Omega.
\end{equation*}
 We then approximate $u$ by
\begin{equation}\label{eq:neural-form-hard-bcij}
	\hat{u}(x; \bm{\theta}_u) = h(x) \mathcal{N}_u(x;\bm{\theta}_u),
\end{equation}
where $\mathcal{N}_u(x; \bm{\theta}_u)$ is a neural network parameterized by $\bm{\theta}_u$.
It is easy to verify that
\begin{equation*}
	\hat{u}(x; \bm{\theta}_u) =  h(x)  \mathcal{N}_u (x; \bm{\theta}_u)= 0, \quad \forall x \in \partial \Omega.
\end{equation*}
Hence, the boundary condition $u|_{\partial \Omega} = 0$ is satisfied by $\hat{u}(x; \bm{\theta}_u)$. 

\begin{remark}
{\em
If the boundary $\partial \Omega$ admits an analytic form, it is usually easy to construct $h$ with analytic expressions, see e.g., \cite{lagaris1998artificial,lai2023hard,lu2021physics} and also Section \ref{se: numerical} for some concrete examples.
Otherwise, we can adopt the method in \cite{sheng2021pfnn} or construct $h$ by training a neural network.
For instance, we can train a  neural network $\hat{h}(x; \theta_h)$ with smooth activation functions (e.g. the sigmoid function or the hyperbolic tangent function) by minimizing the following loss function:
\begin{equation*}
	\frac{w_{1h}}{M_b}\sum_{i=1}^{M_b}|\hat{h}(x_b^i; \theta_{h})|^2+\frac{w_{2h}}{M}\sum_{i=1}^{M}|\hat{h}(x^i; \theta_{h})-\bar{h}(x^i)|^2,
\end{equation*}
where $w_{1h}, w_{2h}>0$ are the weights, $\{x^i\}_{i=1}^{M} \subset \Omega$ and $\{x_b^i\}_{i=1}^{M_b} \subset \partial \Omega$ are sampled points, and $\bar{h}(x)\in C(\Omega)$ is a known function satisfying $\bar{h}(x)\neq 0$ in $\Omega$, e.g. $\bar{h}(x) = \min_{\hat{x}\in \partial\Omega} \{\| x- \hat{x} \|_2^2 \}$.
}
\end{remark}

With the neural network $\hat{u}(x; \bm{\theta}_u)$ given in \eqref{eq:neural-form-hard-bcij}, we approximate the equation \eqref{eq:oc} by 
\begin{equation}\label{eq:ocNN}
	\text{Prox}_{\eta j}\Big((I-\eta A)\hat{u}(x;\bm{\theta}_u)+\eta f(x)\Big)=\hat{u}(x;\bm{\theta}_u),~\text{a.e. in}~\Omega.
\end{equation}
Given a set $\mathcal{T}\subset{\Omega}$, the residual of the equation \eqref{eq:ocNN} can be measured by
\begin{equation}\label{eq:loss}
\mathcal{L}(\bm{\theta}_u)=\frac{1}{|\mathcal{T}|}\sum_{x\in \mathcal{T}}\left|\text{Prox}_{\eta j}\Big((I-\eta A)\hat{u}(x;\bm{\theta}_u)+\eta f(x)\Big)-\hat{u}(x;\bm{\theta}_u)\right|^2.
\end{equation}
As a result, we can train the neural network $\hat{u}(x;\bm{\theta}_u)$ by minimizing the loss function \eqref{eq:loss} and obtain the Prox-PINNs framework for solving (\ref{EVI}), which is listed as Algorithm \ref{alg:dl_evi}. 

\begin{algorithm}[H]
	\caption{The Prox-PINNs Framework for (\ref{EVI}).}
	\label{alg:dl_evi}
	\begin{algorithmic}[1]
		\REQUIRE Parameter $\eta>0$,  auxiliary function $h(x)$.
		\STATE Initialize the neural network $\hat{u}(x;\bm{\theta}_u)$
		\STATE Sample a training set $\mathcal{T} = \{x^i\}_{i=1}^M \subset {\Omega}$ and compute the value of $f$ over $\mathcal{T}$. 
		\STATE Train the neural network $\hat{u}(x;\bm{\theta}_u)$ to identify the optimal parameter $\bm{\theta}_u^*$ by minimizing the loss function \eqref{eq:loss} via a stochastic optimization method.
	\end{algorithmic}
	 \textbf{Output:} An approximate solution $\hat{u}(x;\bm{\theta}_u^*)$ to (\ref{EVI}).
\end{algorithm}

\begin{remark}
{\em
Note that, with the homogeneous Dirichlet boundary condition, the problem \eqref{eq:oc} is
\begin{equation*}
	\quad	\text{Prox}_{\eta j}((I-\eta A)u+\eta f)=u~\text{in}~\Omega,\quad u=0~\text{on}~\partial\Omega,
\end{equation*}
and the boundary condition $u=0$ on $\partial \Omega$ can be treated as a soft constraint by penalizing it in the loss function like the vanilla PINNs \cite{raissi2019physics} . 
In this case, the loss function (\ref{eq:loss}) needs to be revised accordingly to 
\begin{equation}\label{eq:loss_soft}
	\mathcal{L}(\bm{\theta}_u)=\frac{1}{|\mathcal{T}|}\sum_{x\in \mathcal{T}}\left|\text{Prox}_{\eta j}\Big((I-\eta A)\hat{u}(x;\bm{\theta}_u)+\eta f(x)\Big)-\hat{u}(x;\bm{\theta}_u)\right|^2+\frac{w_b}{|\mathcal{T}_b|}\sum_{x\in \mathcal{T}_b}|\hat{u}(x;\bm{\theta}_u)|^2,
\end{equation}
where $w_b>0$ is a weight parameter and $\mathcal{T}_b\subset \partial\Omega$ is a set of randomly sampled points. Note that
the boundary condition $u=0$ on $\partial \Omega$ cannot be strictly enforced under the soft-constraint loss \eqref{eq:loss_soft}. This approach jointly trains the nonlinear equation and boundary condition, and hence its performance depends critically on heuristic weight choices in the loss. However, no systematic principles exist to guide these weights, and setting them manually by trial and error is challenging and time-consuming.}
\end{remark}

\begin{remark}
{\em
	When a non-homogeneous boundary condition $u=u_b (\neq 0)$ on $\partial \Omega$ is considered, one can approximate $u$ by the following neural network 
	\begin{equation}\label{eq: bd_non}
		\hat{u}(x; \theta_u) =g(x)+ h(x) \mathcal{N}_u(x; \bm{\theta}_u),
	\end{equation}
	where the function $g: \bar{\Omega}\rightarrow \mathbb{R}$ is prescribed and satisfies $g\in C(\bar{\Omega})$ and $g\mid_{\partial\Omega}=u_b$.  See Example \ref{ex:2d_ob} for a demonstration. }
\end{remark}

We reiterate that Algorithm \ref{alg:dl_evi} is a high-level framework that imposes no strict restrictions on the operators $A$ and $j$. Meanwhile, the abstract 
and general Algorithm \ref{alg:dl_evi}  becomes practical for a specific EVI only when the analytical formulation of the loss function \eqref{eq:loss} is available, which depends only on the property of the nonsmooth functional $j$. Next, we shall show that the loss function \eqref{eq:loss} admits an analytical form for many nonsmooth functionals $j$ of practical interest, and hence Algorithm \ref{alg:dl_evi} is
feasible for a wide range of EVIs modeled by \eqref{EVI}.

\section{Case Studies for the Implementation of Algorithm \ref{alg:dl_evi}}\label{se:case_studies}

In this section, we present the formal derivation of the loss function for implementing the Prox-PINNs (Algorithm \ref{alg:dl_evi}) to solve EVIs in the form of (\ref{EVI}). Specifically, we consider four distinct types of nonsmooth functionals $j$ in the context of (\ref{EVI}), which are of great practical interest and capture important applications in different fields. 

We begin by addressing the cases where the proximal operator $\text{Prox}_{\eta j}$ admits an explicit analytical form. In this scenario, the loss function can be directly constructed by computing $\text{Prox}_{\eta j}((I-\eta A)u+\eta f) - u$.

\noindent$\bullet$ \textbf{Case 1}. Let $j(\cdot)=I_K(\cdot)$ be the indicator functional of the set ${K}=\{ v\in H_0^1(\Omega)\mid v\geq \psi,~ \text{a.e. in}~ \Omega\}$ with  $\psi\in H^1(\Omega)\cap C^0(\bar{\Omega})$ verifying $\psi\leq 0$ on $\partial \Omega$. Then, the resulting EVI \eqref{EVI} covers a variety of obstacle problems \cite{glowinski2008lectures,glowinski2015variational}.
As discussed in Section \ref{se:obstacle}, we choose $H = L^2(\Omega)$ and redefine $K = \{ v\in L^2(\Omega)\mid v\geq \psi,~ \text{a.e. in}~ \Omega\}$. This gives
$$
\text{Prox}_{\eta j}(w)=\mathop{\arg\min}\limits_{v\in L^2(\Omega)} I_K(v)+\frac{1}{2\eta}\|v-w\|_{L^2(\Omega)}^2, \quad \forall w\in L^2(\Omega),
$$
which implies that
$$
\text{Prox}_{\eta j}(w)(x) = \max\{\psi(x),w(x)\},\quad x\in\Omega, ~ \forall w\in L^2(\Omega).
$$
Hence, the computation of $\text{Prox}_{\eta j}((I-\eta A)u+\eta f) - u$ can be explicitly written as
$$
\max\{\psi(x),(I-\eta A) u(x)+\eta f(x)\} - u(x),
$$
or equivalently
$$
\max\{0,(I-\eta A)u(x)+\eta f(x)-\psi(x)\}+\psi(x) - u(x).
$$
As a result, the loss function \eqref{eq:loss} turns out to be
\begin{equation}\label{eq:loss-1}
	\mathcal{L}(\bm{\theta}_u)=\frac{1}{|\mathcal{T}|}\sum_{x\in \mathcal{T}}\Big|\text {ReLU}\left\{(I-\eta A)\hat{u}(x;\bm{\theta}_u)+\eta f(x)-\psi(x)\right\}+\psi(x)-\hat{u}(x;\bm{\theta}_u)\Big|^2.
\end{equation}

\noindent$\bullet$ \textbf{Case 2}. In this case, we consider $j(v)=\tau\int_\Omega |v| dx$ with $\tau>0$ a constant, which captures important applications in simplified friction problems and image denoising problems \cite{glowinski2015variational}. The corresponding proximal operator is given by
$$
\text{Prox}_{\eta j}(w)=\mathop{\arg\min}\limits_{v\in L^2(\Omega)}\tau \int_\Omega |v| dx+\frac{1}{2\eta}\|v-w\|_{L^2(\Omega)}^2,\quad
 \forall w\in L^2(\Omega),
$$
and hence
$$
\text{Prox}_{\eta j}(w)(x)=\mathbb{S}_{\tau\eta}(w)(x):=\text{sgn}(w(x))\max \{|w(x)|-\tau\eta,0\}~\text{a.e. in }~\Omega.
$$
The above result implies that the loss function \eqref{eq:loss} can be explicitly specified as
{\small
\begin{equation*}
	\mathcal{L}(\bm{\theta}_u)=\frac{1}{|\mathcal{T}|}\sum_{x\in \mathcal{T}}\Big|\text{sgn}\big((I-\eta A)\hat{u}(x;\bm{\theta}_u)+\eta f(x)\big)\text{ReLU}\left \{\big|(I-\eta A)\hat{u}(x;\bm{\theta}_u)+\eta f(x)\big|-\tau\eta\right\}-\hat{u}(x;\bm{\theta}_u)\Big|^2.
\end{equation*}}

Next, we consider a more complex scenario where the nonsmooth functional $j$ in (\ref{EVI}) has a composite structure. Specifically, let $j(u) = g(Bu)$, where the functional $g: \mathcal{H}\rightarrow \mathbb{R}\cup\{+\infty\}$, with $\mathcal{H}$ a Hilbert space, is non-smooth, convex, proper and lower semi-continuous, and the operator $B: V\rightarrow \mathcal{H}$ is assumed to be linear and continuous. The proximal operator of $g$ is thus defined as
$$\text{Prox}_{\eta g}(w)=\mathop{\arg\min}\limits_{v\in \mathcal{H}} g(v)+\frac{1}{2\eta}\|v-w\|_{\mathcal{H}}^2,\forall w\in \mathcal{H}.$$
As we shown subsequently, $\text{Prox}_{\eta g}(w)$ typically admits an explicit analytical form. Nevertheless,  the proximal operator of $j$ generally exhibits no closed-form expression and is computationally expensive to evaluate.

To address this issue, we note that, as shown in \cite[Proposition 6.19, Corollary 16.42]{bauschkeconvex}, if $0 \in \mathrm{int} (\mathrm{dom}\,g - \mathrm{ran}\, B)$, then it holds that $\partial j(u) = B^*\partial g( Bu)$ with $B^*$ the adjoint operator of $B$. By introducing an auxiliary variable $\bm{\lambda} \in  \mathcal{H}$, we can reformulate equation (\ref{eq:inclusion}) as follows
\[
u\in V, \quad \bm{\lambda} \in  \mathcal{H}, \quad -\bm{\lambda} \in \partial g( Bu), \quad (Au-f) - B^*\bm{\lambda} = 0.
\]
Let $\eta\in \mathbb{R}$ be a positive constant. We can further reformulate the above equation as
\[
u\in V, \quad \bm{\lambda} \in  \mathcal{H}, \quad Bu - \eta\bm{\lambda} \in Bu + \eta\partial g(Bu), \quad (Au-f) - B^*\bm{\lambda}= 0, 
\]
which, using $\text{Prox}_{\eta g} = (I + \eta \partial g)^{-1}$, can be rewritten  as
\begin{equation}\label{eq:inclusion2}
	u\in V, \quad \bm{\lambda} \in  \mathcal{H}, \quad \text{Prox}_{\eta g}( Bu - \eta\bm{\lambda}) = Bu, \quad (Au-f) - B^*\bm{\lambda} = 0.
\end{equation}

To approximate $\bm{\lambda}$, we introduce a neural network $\hat{\bm{\lambda}}(x;\bm{\theta}_{\lambda})$. The loss function is then determined based on the residual of equation (\ref{eq:inclusion2}) for training the neural networks $\hat{u}(x;\bm{\theta}_u)$ and $\hat{\bm{\lambda}}(x;\bm{\theta}_{\lambda})$,
\begin{equation}\label{eq:loss_2}
		\begin{aligned}
			\mathcal{L}(\bm{\theta}_u,\bm{\theta}_{\lambda})=\frac{1}{|\mathcal{T}|}\sum_{x\in \mathcal{T}}\Bigg\{w_1\left|\text{Prox}_{\eta g}\Big( B\hat{u}(x;\bm{\theta}_u) - \eta \hat{\bm{\lambda}}(x;\bm{\theta}_{\lambda})\Big)-B\hat{u}(x;\bm{\theta}_u)\right|^2\\
			+w_2\left|A\hat{u}(x;\bm{\theta}_u)- f(x)-B^*\hat{\bm{\lambda}}(x;\bm{\theta}_{\lambda})\right|^2
			\Bigg\}.
		\end{aligned}
\end{equation}

In the following, we consider two specific cases to further illustrate how to construct the loss function (\ref{eq:loss_2}) for training the neural networks $\hat{u}(x;\bm{\theta}_u)$ and $\hat{\bm{\lambda}}(x;\bm{\theta}_{\lambda})$.

\noindent$\bullet$ \textbf{Case 3}. We consider  $V=H_0^1(\Omega), j(\cdot)=I_K(\cdot)$ as the indicator functional of the set ${K}=\{v\in H_0^1(\Omega)\mid |\nabla  v(x)|\leq 1,~ \text{a.e. in}~ \Omega\}$, where $|\nabla  v(x)|=\left([\frac{\partial v}{\partial x_1}(x)]^2+\cdots+[\frac{\partial v}{\partial x_d}(x)]^2\right)^{\frac{1}{2}}$.  As a result, EVI \eqref{EVI} covers the elasto-plastic torsion 
problem, see e.g., \cite{duvant1976inequalities,glowinski2015variational} and references therein. 

In this case, we have $j(u) = g(Bu)$ with $B = \nabla$, $B^*=-\text{div}$, $ \mathcal{H}=[L^2(\Omega)]^d$, and $g(\cdot) = I_{\tilde{K}} (\cdot)$, where $\tilde{K}=\{\bm{q}| \bm{q}\in [L^2(\Omega)]^d, |\bm{q}(x)|\leq 1,~ \text{a.e. in}~ \Omega\}$ with $|\bm{q}(x)|=\left([\bm{q}_1(x)]^2+\cdots+[\bm{q}_d(x)]^2\right)^{\frac{1}{2}}$. Thus, we have
$$
\text{Prox}_{\eta g}(\bm{w})(x)= P_{\tilde{K}}(\bm{w}(x)):= \frac{\bm{w}(x)}{\max\{1, |\bm{w}(x)|\}}, \quad \forall \bm{w}\in [L^2(\Omega)]^d.
$$
This result implies that the loss function (\ref{eq:loss_2}) can be specified as 
\begin{equation}\label{eq:loss_3}
		\begin{aligned}
			\mathcal{L}(\bm{\theta}_u,\bm{\theta}_{\lambda})=\frac{1}{|\mathcal{T}|}\sum_{x\in \mathcal{T}}\Bigg\{w_1\left|\frac{\nabla\hat{u}(x;\bm{\theta}_u) - \eta\hat{\bm{\lambda}}(x;\bm{\theta}_{\lambda})}{\text{ReLU}\{|\nabla \hat{u}(x;\bm{\theta}_u) - \eta\hat{\bm{\lambda}}(x;\bm{\theta}_{\lambda})|-1\}+1} - \nabla\hat{u}(x;\bm{\theta}_u) \right|^2\\
			+w_2\left|A\hat{u}(x;\bm{\theta}_u)- f(x)+\nabla\cdot\hat{\bm{\lambda}}(x;\bm{\theta}_{\lambda})\right|^2
			\Bigg\}.
		\end{aligned}
\end{equation}

\noindent$\bullet$ \textbf{Case 4}. Finally, we consider $V=H_0^1(\Omega)$ and $j(v)=\tau\int_\Omega |\nabla v| dx$ with $\tau>0$ a constant, which is used in modeling Bingham visco-plastic flows \cite{Dean2007On} and  image restoration problems \cite{chan1999nonlinear}.

In this case, we have $j(u) = g(Bu)$ with $B = \nabla$, $B^*=-\text{div}$, $ \mathcal{H}=[L^2(\Omega)]^d$, and $g(\cdot) =\tau\int_\Omega |\cdot| dx$. Then, it is easy to show that
\begin{equation}\label{eq:prox_case4}
\text{Prox}_{\eta g}(\bm{w})(x)=\frac{\bm{w}(x)}{|\bm{w}(x)|}\max \{|\bm{w}(x)|-\tau\eta,0\}, ~\text{a.e. in }~\Omega,\quad \forall \bm{w}\in [L^2(\Omega)]^d.
\end{equation}
This implies that the loss function (\ref{eq:loss_2}) can be specified as 
\begin{equation}\label{eq:loss_4}
	{\small
	\begin{aligned}
		\mathcal{L}(\bm{\theta}_u,\bm{\theta}_{\lambda})=\frac{1}{|\mathcal{T}|}\sum_{x\in \mathcal{T}}\Bigg\{
		w_1	\left| \frac{ \nabla\hat{u}(x;\bm{\theta}_u) - \eta\hat{\bm{\lambda}}(x;\bm{\theta}_{\lambda})}{|\nabla\hat{u}(x;\bm{\theta}_u) - \eta\hat{\bm{\lambda}}(x;\bm{\theta}_{\lambda}) |}\text{ReLU}\left \{\big|\nabla\hat{u}(x;\bm{\theta}_u) - \eta\hat{\bm{\lambda}}(x;\bm{\theta}_{\lambda})\big|-\tau\eta\right\}- \nabla \hat{u}(x;\bm{\theta}_u) \right|^2\\
		+ w_2\left|A\hat{u}(x;\bm{\theta}_u)- f(x)+\nabla\cdot\hat{\bm{\lambda}}(x;\bm{\theta}_{\lambda})\big)\right|^2
		\Bigg\}.
	\end{aligned}
}
\end{equation}
 
 \begin{remark}
 {\em It follows from \eqref{eq:inclusion2} and \eqref{eq:prox_case4} that, for  \textbf{Case 4}, $u\in H_0^1(\Omega)$ and $\bm{\lambda} \in [L^2(\Omega)]^d $ satisfy
 \begin{numcases}
 		~\frac{\nabla u(x)- \eta\bm{\lambda}(x)}{|\nabla u(x) - \eta\bm{\lambda}(x)|}\max \{|\nabla u(x) - \eta\bm{\lambda}(x)|-\tau\eta,0\}=\nabla u(x),~\text{a.e.~in}~\Omega,\label{eq:o1_case4}\\
 	(Au-f) +\nabla\cdot\bm{\lambda} = 0,~\text{a.e.~in}~\Omega.\label{eq:o2_case4}
 \end{numcases}
 It is easy to show that the equation \eqref{eq:o1_case4} is equivalent to 
 \begin{equation}\label{eq:oc2_gradient_l1}
 	\left\{
 	\begin{aligned}
 			&\bm{\lambda}(x)\cdot \nabla u(x)=-\tau|\nabla u(x)|,\\
 			&|\bm{\lambda}(x)|\leq \tau ~(i.e.~ \bm{\lambda}(x)=\frac{\tau\bm{\lambda}(x)}{\max\{\tau,|\bm{\lambda}(x)|\}}).
 		\end{aligned}
 	\right.
 \end{equation}
 Indeed, if $\nabla u=0$, then we have $|\nabla u(x) - \eta\bm{\lambda}(x)|-\tau\eta\leq 0$ and thus $|\bm{\lambda}(x)|\leq \tau$. On the other hand, if $\nabla u\neq0$, then it holds that $|\nabla u(x) - \eta\bm{\lambda}(x)|-\tau\eta> 0$, which implies that 
 \begin{equation}\label{eq:eq1}
  \frac{\nabla u(x)- \eta\bm{\lambda}(x)}{|\nabla u(x) - \eta\bm{\lambda}(x)|}(|\nabla u(x) - \eta\bm{\lambda}(x)|-\tau\eta)=\nabla u(x),~\text{a.e.~in}~\Omega,
 \end{equation}
and hence
\begin{equation}\label{eq:eq2}
	\frac{\nabla u(x)- \eta\bm{\lambda}(x)}{|\nabla u(x) - \eta\bm{\lambda}(x)|}= \frac{\nabla u(x)}{|\nabla u(x)|}.
\end{equation}
It follows from \eqref{eq:eq1} and \eqref{eq:eq2} that
\begin{equation*}
\bm{\lambda}(x)=-\tau\frac{\nabla u(x)}{|\nabla u(x)|}.
\end{equation*}
We thus get the desired result. Note that the converse of the above arguments also holds.

To guarantee $|\bm{\lambda}(x)|\leq \tau $, we use 
	$$
	\hat{\bm{\lambda}}(x;\bm{\theta}_{\lambda})=\frac{\tau\mathcal{N}_{\lambda}(x;\bm{\theta}_{\lambda})}{\max\{\tau,|\mathcal{N}_{\lambda}(x;\bm{\theta}_{\lambda})|\}}
	$$ 
	with $\mathcal{N}_{\lambda}(x;\bm{\theta}_{\lambda})$ a neural network to approximate 
	$\bm{\lambda}$.
Then, one can use the residuals of \eqref{eq:o2_case4} and \eqref{eq:oc2_gradient_l1} to train $\hat{u}(x;\bm{\theta}_{u})$ and $\hat{\bm{\lambda}}(x;\bm{\theta}_{\lambda})$ and the resulting loss function reads as
 \begin{equation}\label{eq:loss_4_2}
 	{\small
 	\begin{aligned}
 			\mathcal{L}(\bm{\theta}_u,\bm{\theta}_{\lambda})=\frac{1}{|\mathcal{T}|}\sum_{x\in \mathcal{T}}\Bigg\{
 			w_1	\left|\hat{\bm{\lambda}}(x;\bm{\theta}_{\lambda})\cdot \nabla \hat{u}(x;\bm{\theta}_u)+\tau|\nabla \hat{u}(x;\bm{\theta}_u)|\right|^2
 			+w_2\left|A\hat{u}(x;\bm{\theta}_u)- f(x)+\nabla\cdot\hat{\bm{\lambda}}(x;\bm{\theta}_{\lambda})\right|^2
 			\Bigg\}
 		\end{aligned}}
 \end{equation}
 with $w_1>0$ and $w_2>0$ the weights for each component. Compared with  \eqref{eq:loss_4}, the loss function \eqref{eq:loss_4_2} employs fewer residual terms, thereby reducing training complexity. 
 Furthermore, the denominator  $|\nabla\hat{u}(x;\bm{\theta}_u) - \eta\hat{\bm{\lambda}}(x;\bm{\theta}_{\lambda}) |$  in \eqref{eq:loss_4} can approach zero during computation, potentially causing numerical instability. 
 This issue is avoided by adopting \eqref{eq:loss_4_2}, thereby substantially improving numerical stability.
}
 \end{remark}
 
\begin{figure}[htpb]
  \centering
 \begin{tikzpicture}[scale=0.4]
	\node at(1,8.5)    (1) [circle,fill=pink!60,draw=black,minimum width =20pt, minimum height =20pt,font=\fontsize{10}{10}\selectfont]{$x$};
	\node at(5,13)    (2) [circle,fill=blue!35,draw=black,minimum width =20pt, minimum height =20pt,font=\fontsize{25}{25}\selectfont]{};
	\node at(5,10)    (3) [circle,fill=blue!35,draw=black,minimum width =20pt, minimum height =20pt,font=\fontsize{25}{25}\selectfont]{};
	\node at(5,7)    (4) [circle,fill=blue!35,draw=black,minimum width =20pt, minimum height =20pt,font=\fontsize{25}{25}\selectfont]{};
	\node at(5,4)    (5) [circle,fill=blue!35,draw=black,minimum width =20pt, minimum height =20pt,font=\fontsize{25}{25}\selectfont]{};
	\node at(9,15)    (6) [circle,fill=blue!35,draw=black,minimum width =20pt, minimum height =20pt,font=\fontsize{25}{25}\selectfont]{};
	\node at(9,12)    (7) [circle,fill=blue!35,draw=black,minimum width =20pt, minimum height =20pt,font=\fontsize{25}{25}\selectfont]{};
	\node at(9,9)    (8) [circle,fill=blue!35,draw=black,minimum width =20pt, minimum height =20pt,font=\fontsize{25}{25}\selectfont]{};
	\node at(9,6)    (9) [circle,fill=blue!35,draw=black,minimum width =20pt, minimum height =20pt,font=\fontsize{25}{25}\selectfont]{};
	\node at(9,3)    (10) [circle,fill=blue!35,draw=black,minimum width =20pt, minimum height =20pt,font=\fontsize{25}{25}\selectfont]{};
	\node at(13,13)    (11) [circle,fill=orange!35,draw=black,minimum width =20pt, minimum height =20pt,font=\fontsize{25}{25}\selectfont]{};
		\node at(13,7)    (14)[box, fill=gray!20,minimum width =25pt, minimum height =88pt,draw=black]  {}; 
	\node at(13,9)    (12) [circle,fill=orange!35,draw=black,minimum width =20pt, minimum height =20pt,font=\fontsize{25}{25}\selectfont]{};
	\node at(13,5)    (13) [circle,fill=orange!35,draw=black,minimum width =20pt, minimum height =20pt,font=\fontsize{25}{25}\selectfont]{};
	\node at(16.5,13)  (15) [font=\fontsize{10}{10}\selectfont]{$\mathcal{N}_{u}(x;\bm{\theta}_u)$};
	\node at(16.7,7)  (16) [font=\fontsize{10}{10}\selectfont]{$\mathcal{N}_{\lambda}(x;\bm{\theta}_\lambda)$};
	\node at(24.5,13)    (17)[rectangle, minimum width =80pt, minimum height =20pt, inner sep=3pt,font=\fontsize{9}{9}\selectfont,line width=0.6 pt,draw=black]  {$\hat{u}(x;\bm{\theta}_u)=h(x)\mathcal{N}_{u}(x;\bm{\theta}_u)$}; 
	\node at(28,7)    (18)[rectangle, minimum width =120pt, minimum height =40pt, inner sep=3pt,font=\fontsize{9}{9}\selectfont,line width=0.6 pt,draw=black]  {$\begin{aligned}
			&\hat{\bm{\lambda}}(x; \bm{\theta}_{\lambda})=\mathcal{N}_{\lambda}(x; \bm{\theta}_{\lambda})&\text{for}~\textbf{Case~3}\\
			&\hat{\bm{\lambda}}(x; \bm{\theta}_{\lambda})=\frac{\tau \mathcal{N}_{\lambda}(x; \bm{\theta}_{\lambda})}{\max\{\tau,|\mathcal{N}_{\lambda}(x; \bm{\theta}_{\lambda})|\}}&\text{for}~\textbf{Case~4}
		\end{aligned}$}; 
	\draw(1)--(2);
	\draw (1)--(3);	
	\draw (1)--(4);	
	\draw (1)--(5);	
	\draw (2)--(6);	
	\draw (2)--(7);	
	\draw (2)--(8);	
	\draw (2)--(9);	
	\draw (2)--(10);	
	\draw (3)--(6);	
	\draw (3)--(7);	
	\draw (3)--(8);	
	\draw (3)--(9);	
	\draw (3)--(10);	
	\draw (4)--(6);	
	\draw (4)--(7);	
	\draw (4)--(8);	
	\draw (4)--(9);	
	\draw (4)--(10);	
	\draw (5)--(6);	
	\draw (5)--(7);	
	\draw (5)--(8);	
	\draw (5)--(9);	
	\draw (5)--(10);	
	\draw (6)--(11);	
	\draw (6)--(12);	
	\draw (6)--(13);	
	\draw (7)--(11);	
	\draw (7)--(12);	
	\draw (7)--(13);	
	\draw (8)--(11);	
	\draw (8)--(12);	
	\draw (8)--(13);	
	\draw (9)--(11);	
	\draw (9)--(12);	
	\draw (9)--(13);	
	\draw (10)--(11);	
	\draw (10)--(12);	
	\draw (10)--(13);	
    \draw[->](11)--(15);
	\draw[->] (14)--(16);
	\draw[->] (15)--(17);
	\draw[->] (16)--(18);
\end{tikzpicture}
 \caption{An illustrative example for the neural network architectures of $\hat{u}(x;\bm{\theta}_u)$ and $\hat{\bm{\lambda}}(x;\bm{\theta}_\lambda)$ when $A$ is linear and $d=2$.}\label{fig: nn}
\end{figure}
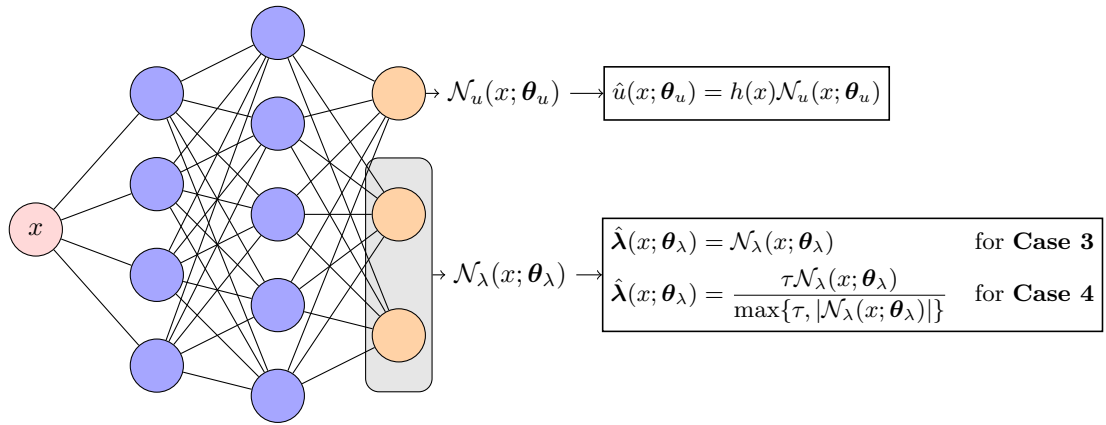

The results presented above demonstrate the broad applicability of Algorithm \ref{alg:dl_evi}. Next, we present some remarks on the neural network architecture in Algorithm \ref{alg:dl_evi} to complete the discussions on its implementation.  Suppose that we consider a neural network that takes spatial coordinates $x \in \mathbb{R}^d$ as input, where $d>0$ denotes the problem's dimension. For \textbf{Cases}~$\bm{1}$ and~$\bm{2}$, the neural network outputs a scalar $h(x)\mathcal{N}_u(x;\bm{\theta}_u)$ to approximate the solution $u(x)$. In \textbf{Cases}~$\bm{3}$ and~$\bm{4}$, we introduce the variable $\bm{\lambda}\in [L^2(\Omega)]^d$ to derive the explicit formulation of \eqref{eq:oc}, which would conceptually necessitate a separate network $\hat{\bm{\lambda}}(x,\bm{\theta}_\lambda)$, thereby increasing training complexity. To alleviate this issue, we observe from \eqref{eq:inclusion2}  that $u$ and $\bm{\lambda}$ share an affine relationship. Inspired by this insight, we expand the output dimension of the neural network to $d+1$ to incorporate $\bm{\lambda}(x)\in\mathbb{R}^d$, where the first component approximates $u(x)$ while the remaining $d$ components represent the approximation to $\bm{\lambda}(x)$, as illustrated in Figure \ref{fig: nn}. This unified architecture offers two significant computational and theoretical advantages:
 \begin{enumerate}
 	\item \textbf{Computational Efficiency}: 
 	Compared to using two separate neural networks (one for $u$ and the other for $\bm{\lambda}$), sharing parameters between $u$ and $\bm{\lambda}$ within a single neural network significantly reduces the computational costs and simplifies the training processes.
 	
 	\item \textbf{Mathematical Consistency}: 
The affine coupling between \( u \) and \( \bm{\lambda} \), inherent to the problem’s structure, aligns naturally with neural network design. Specifically, the output layer of a neural network is an affine transformation of its final hidden layer. As a result,  this architectural property inherently enforces the theoretical affine relationship between \( u \) and \( \bm{\lambda} \), ensuring consistency with the governing equations.
 \end{enumerate}
Notably, the above discussions can be similarly extended to the cases where $A$ is a nonlinear operator. For instance, one can take $\mathcal{N}_{\lambda}(x; \bm{\theta}_{\lambda})= \mathcal{N}(x;\bm{\theta})\circ\mathcal{N}_u(x;\bm{\theta}_u)$ with $\mathcal{N}(x;\bm{\theta})$ a (shallow) neural network parameterized by $\bm{\theta}$.
 
\section{Applications and Numerical Simulations}\label{se: numerical}
In this section, we implement Algorithm \ref{alg:dl_evi} to various concrete EVIs. To this end, we consider four classic and important EVIs, including obstacle problems, elasto-plastic torsion problems, Bingham visco-plastic flows, and simplified friction problems. For each EVI, numerical results of some benchmark examples are presented to validate the effectiveness, efficiency, accuracy, and robustness of Algorithm \ref{alg:dl_evi}. Some comparisons with the reference ones obtained by FEM-based traditional numerical methods and other deep learning methods are also included. All codes in the numerical experiments were written in Python and PyTorch, and are publicly available on GitHub at: \url{https://github.com/yugaomath/Prox-PINNs}.  The numerical experiments were conducted on a MacBook Pro with mac OS Monterey, Intel(R) Core(TM) i7-9570h (2.60 GHz), and 16 GB RAM.

Throughout, all the neural networks (outlined in Figure \ref{fig: nn}) are set as fully connected neural networks equipped with \texttt{tanh} activation functions. Unless otherwise specified,  the neural networks are initialized by the default PyTorch settings and trained by an ADAM optimizer with a learning rate $10^{-3}$. All the weights in the loss functions are set to be 1. Other parameter settings for different test problems are summarized in Table \ref{tab:para}.

\begin{table}[H]
					\centering
	{\small
	\begin{tabular}{|l|cc|c|cc|}
		\hline
		\multirow{2}{*}{Examples}  & \multicolumn{2}{c|}{size of data set}        & \multirow{2}{*}{training epochs} & \multicolumn{2}{c|}{neural networks}         \\ \cline{2-3} \cline{5-6} 
		& \multicolumn{1}{c|}{training} & test &                                  & \multicolumn{1}{c|}{hidden layers} & neurons \\ \hline
		1D obstacle problems               & \multicolumn{1}{c|}{$50$}         & $10^3$   & $1\times 10^4$                   & \multicolumn{1}{c|}{3}             & 100     \\ \hline
		2D obstacle problems               & \multicolumn{1}{c|}{$10^3$}       & $10^4$   & $1\times 10^4$                   & \multicolumn{1}{c|}{5}             & 100     \\ \hline
		2D elasto-plastic torsion problems & \multicolumn{1}{c|}{$10^3$}       & $10^4$   & $1\times 10^4$                   & \multicolumn{1}{c|}{3}             & 100     \\ \hline
		2D Bingham visco-plastic flows  & \multicolumn{1}{c|}{$10^3$}       & $10^4$   & $2\times 10^4$                   & \multicolumn{1}{c|}{10}            & 50      \\ \hline
		2D simplified friction problems    & \multicolumn{1}{c|}{$10^3$}       & $10^4$   & $1\times 10^4$                   & \multicolumn{1}{c|}{4}             & 50      \\ \hline
	\end{tabular}
}
	\caption{Parameter settings for different test problems}\label{tab:para}
\end{table}

\subsection{Obstacle Problems}
Let $\Omega$ be a bounded domain of $\mathbb{R}^d (d\geq 1)$ and $\partial\Omega$ its boundary. Suppose an elastic
membrane occupy $\Omega$ and  this membrane is fixed
along $\partial\Omega$. Obstacle problems aim to find the equilibrium position $u$ of the elastic membrane under the action of the vertical force $f$, which can be modeled by the EVI:
\begin{equation*}
	u\in {K},\quad\text{such that}\quad \int_\Omega Au(v-u) dx \geq \int_\Omega f(v-u)dx,\quad \forall v\in {K},
\end{equation*}
where $Av=-\alpha \Delta v+\bm{\beta}\cdot \nabla v+\gamma v$, $\forall v\in H_0^1(\Omega)$ with $\alpha,\gamma\in L^{\infty}(\Omega)$, $\alpha>0, \gamma\geq0$ a.e. in $\Omega$, $\bm{\beta}\in [L^{\infty}(\Omega)]^d$, $\nabla\cdot\bm{\beta}=0$, and $f\in L^{2}(\Omega)$. The set
${K}=\{v\mid v\in H_0^1(\Omega),v\geq \psi,~ \text{a.e. in}~ \Omega\}$ with  $\psi\in H^1(\Omega)\cap C^0(\bar{\Omega})$ verifying $\psi\leq 0$ on $\partial \Omega$.  Obstacle problems have numerous applications in diverse scientific areas, such as contact mechanics, processes in biological cells, ecology, fluid flow, and finance, see for example \cite{kinderlehrer1980introduction,rodrigues1987obstacle,sofonea2009variational}. Obstacle problems  have been extensively studied both numerically and theoretically in the literature, see e.g., \cite{glowinski2008lectures,glowinski2015variational} and references therein.

\begin{example}\label{ex:1d_ob_sy}
We first consider a one-dimensional problem with a symmetric elliptic operator, which has been intensively investigated in the literature, see \cite{bahja2023physics,cheng2023deep,liu2023fast,Tran2014An,zosso2017efficient}. 
Let \(\Omega = (0, 1)\), \(A v= -v_{xx},\forall v\in H^1_0(\Omega)\), and \(f \equiv 0\). The functions \(\psi(x)\) and \(u(x)\) are defined as follows:
\begin{equation*}
	{\small
	\psi(x) := \begin{cases} 
		100 x^2, & x \in [0, 0.25], \\ 
		100 x(1 - x) - 12.5, & x \in (0.25, 0.5], \\ 
		\psi(1 - x), & x \in (0.5, 1],
	\end{cases}
		u(x) = 
	\begin{cases} 
		(100 - 50 \sqrt{2}) x, &x\in [0, \frac{1}{2 \sqrt{2}}), \\ 
		100 x(1 - x) - 12.5, & x\in [\frac{1}{2 \sqrt{2}}, 1 - \frac{1}{2 \sqrt{2}}), \\ 
		u(1 - x), & x\in[1 - \frac{1}{2 \sqrt{2}}, 1].
	\end{cases}}
\end{equation*}
It is easy to verify that $u, f$, and $\psi$ satisfy the equation \eqref{eq:oc} with $j$ the indicator functional of $K$. Hence, $u$ is the exact solution to this example. 

To implement Algorithm \ref{alg:dl_evi}, we take $h(x)=x(1-x)$ and $\eta=10^{-3}$. The numerical results are presented in Figure \ref{fig:1d-obstacle}, where we plot the exact and the learned solutions, the point-wise error, the training trajectories for the loss function, and the test errors with respect to training epochs.  We observe that the numerical solution is in
good agreement with the exact one. Moreover, the results validate that Algorithm \ref{alg:dl_evi} can produce numerical solutions with low relative $L^2$- and $L^{\infty}$- errors.  

\begin{figure}[h!]\label{fig:1d-obstacle}
	\centering
	\includegraphics[width=.45\textwidth]{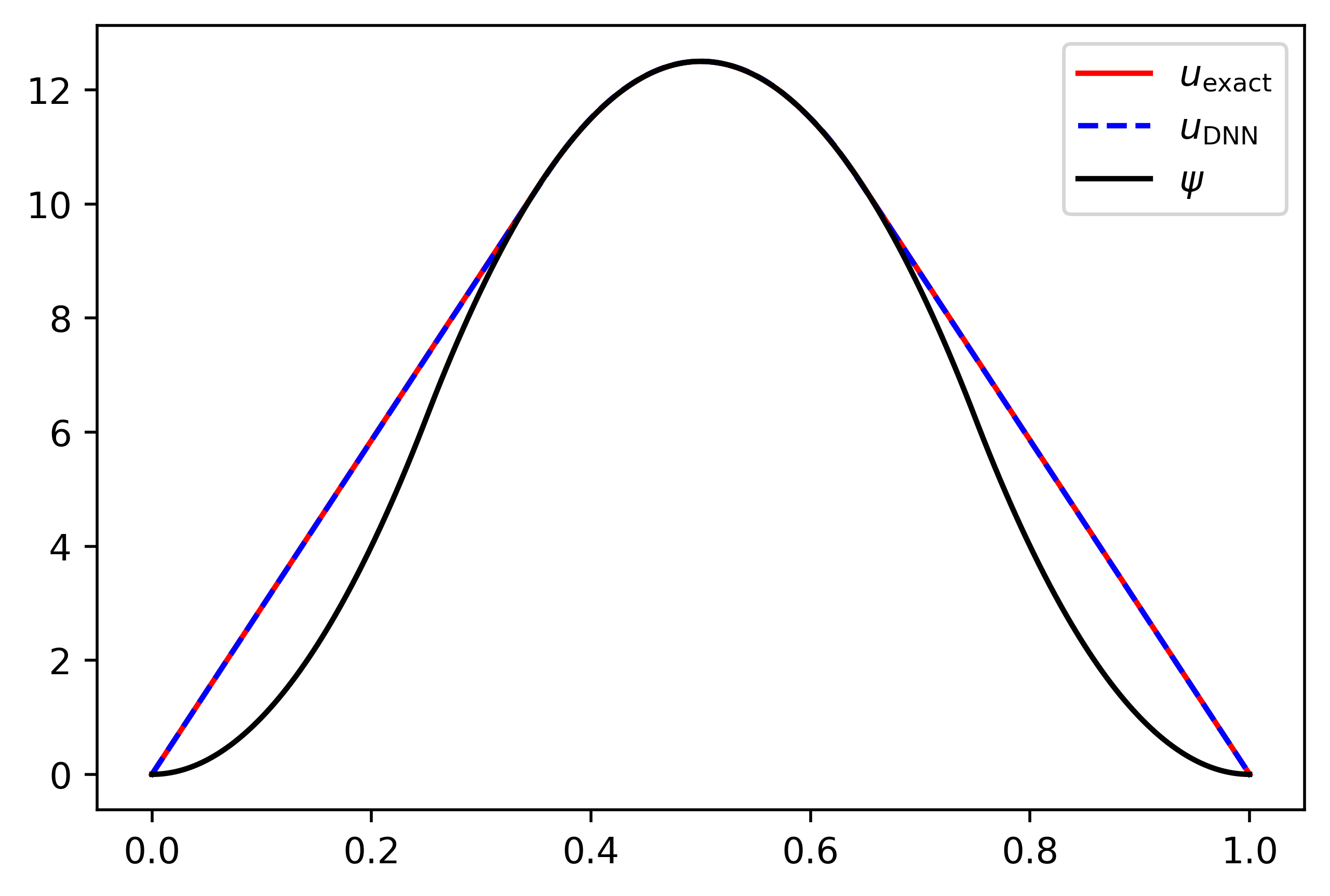}
	\includegraphics[width=.468\textwidth]{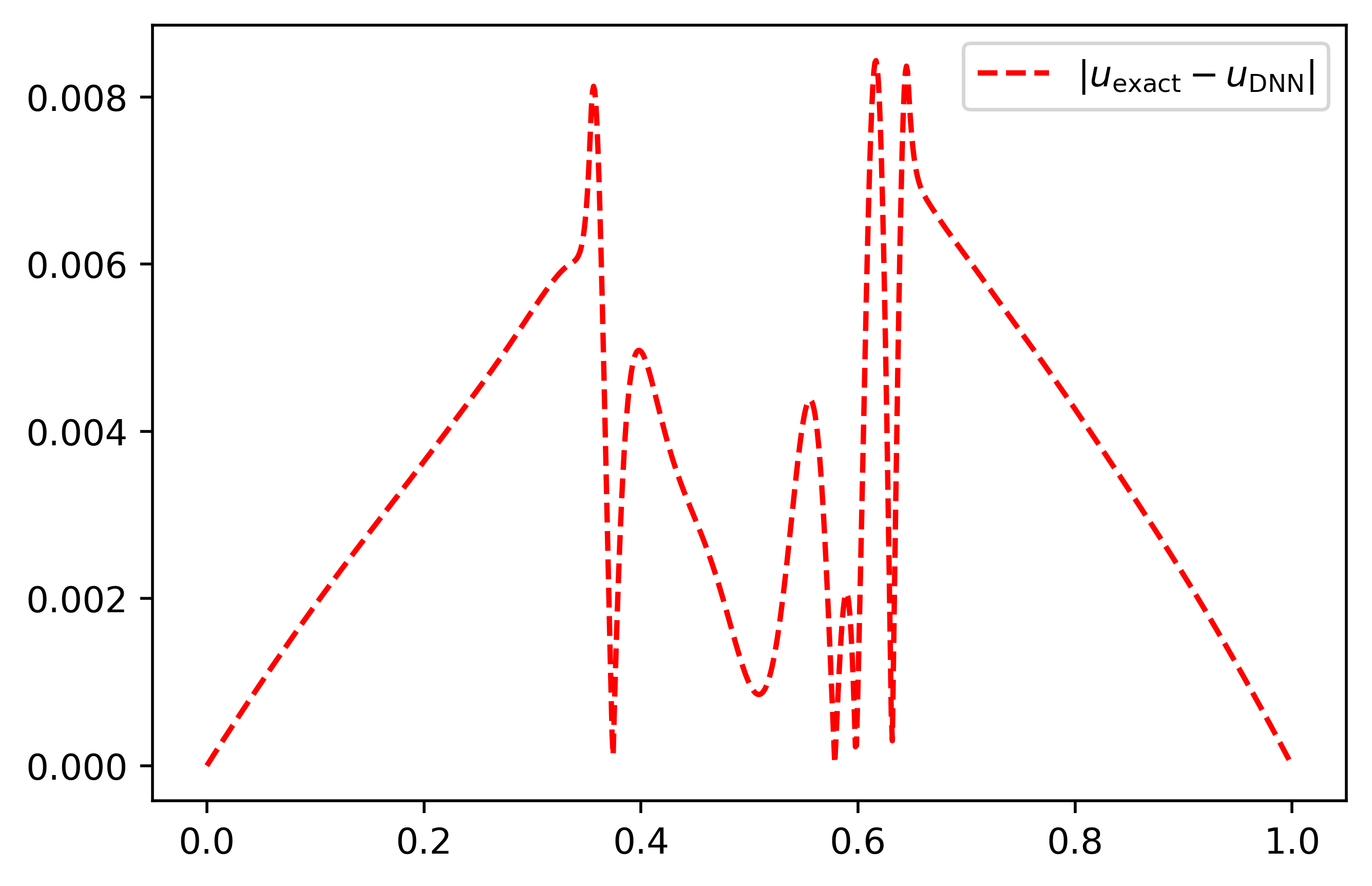}
	\includegraphics[width=.45\textwidth]{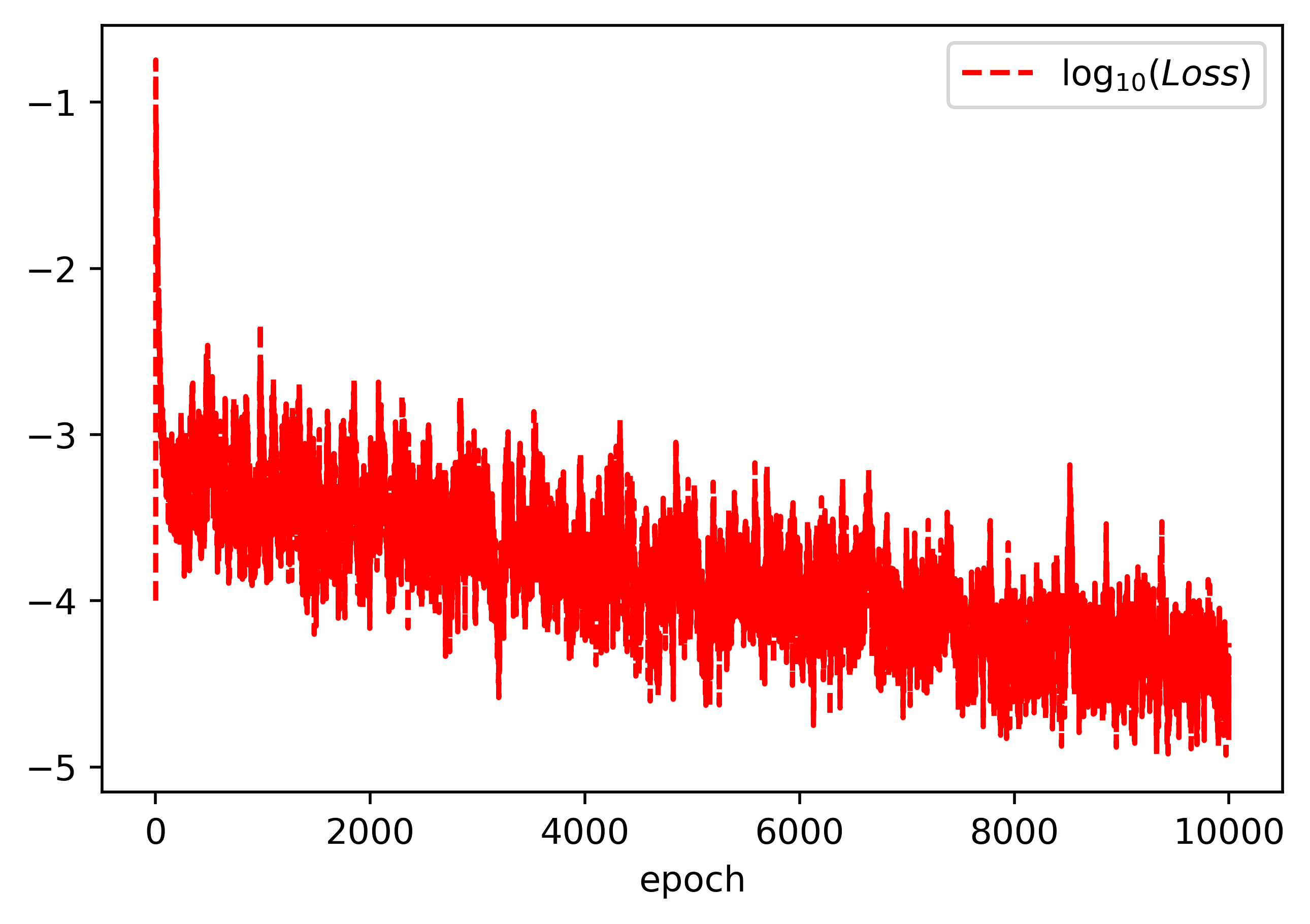}
	\includegraphics[width=.468\textwidth]{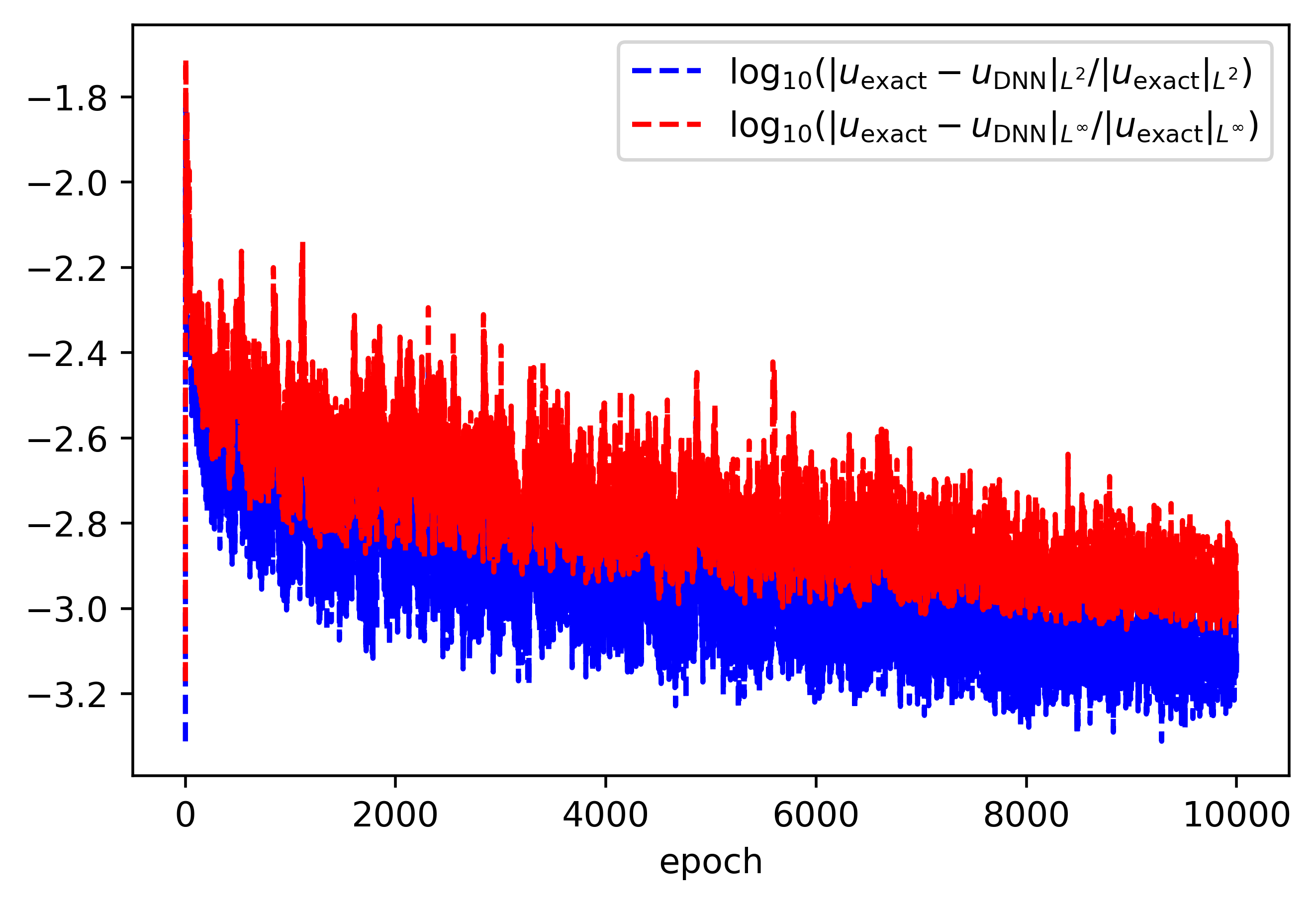}
	\caption{Numerical results for Example \ref{ex:1d_ob_sy} (Relative $L^2$-error: $6.998\times 10^{-4}$; Relative $L^{\infty}$-error: $1.042\times 10^{-3}$). 
	}
\end{figure}

To further validate the accuracy of Algorithm \ref{alg:dl_evi}, we compare it with the deep learning method in \cite{cheng2023deep}.
In \cite{cheng2023deep},  the neural networks are trained with different grid resolutions $N$ and tested on a uniform mesh with $10^3$ grids. Following the settings in \cite{cheng2023deep}, we train the neural networks for $5,000$ epochs and use the metric
$
\frac{1}{N}\sum_{i=1}^N\frac{|u(x_i)-\hat{u}(x_i; \bm{\theta}_u^*)|}{|u(x_i)|}
$
 to evaluate the numerical accuracy of the computed solutions. The deep learning method in \cite{cheng2023deep} is implemented using the  source code publicly available at \url{https://github.com/Xingbaji/Obstacle-Problem} with the parameters given in \cite{cheng2023deep}. The numerical comparisons are reported in Table \ref{tab: compare_1d_obs_sym}. 

	\begin{table}[H]
	\centering
	\scalebox{0.75}
	{\begin{tabular}{|c|c|c|c|c|c|c|c|}
			\hline
			$N$             & 20 & 50 & $10^2$ & $2\times 10^2$ & $5\times 10^2$ & $10^3$ & $10^4$ \\ \hline
			The deep learning method in \cite{cheng2023deep} & $4.7\times 10^{-1}$    & $3.5\times 10^{-1}$    &  $1.4\times 10^{-1}$     &   $6.7\times 10^{-2}$   &  $3.5\times 10^{-2}$   & $1.9\times 10^{-2}$     &  $4.2\times 10^{-3}$     \\ \hline
			Algorithm  \ref{alg:dl_evi}   & $1.8\times 10^{-3}$    & $1.1\times 10^{-3}$    & $1.1\times 10^{-3}$     &  $8.2\times 10^{-4}$    &    $1.0\times 10^{-3}$  &    $8.1\times 10^{-4}$   &  $8.1\times 10^{-4}$      \\ \hline
	\end{tabular}}
	\caption{Comparisons with the deep learning method in \cite{cheng2023deep} for Example \ref{ex:1d_ob_sy}.}\label{tab: compare_1d_obs_sym}
\end{table}

The results in Table \ref{tab: compare_1d_obs_sym} demonstrate that the numerical accuracy of the solutions computed by Algorithm \ref{alg:dl_evi} is significantly higher than that of the solutions produced by the deep learning method in \cite{cheng2023deep}, across varying values of $N$. Furthermore, Algorithm \ref{alg:dl_evi} exhibits strong robustness to the number of training points, indicating that high prediction accuracy can be maintained even with limited training data. These findings collectively validate the efficiency, accuracy, and robustness of Algorithm \ref{alg:dl_evi}, making it an attractive mesh-free method for obstacle problems.

	Furthermore, recall that the hyperparameter $\eta$ is introduced in Algorithm \ref{alg:dl_evi} (see (\ref{eq:loss-1}) for the details related to obstacle problems). To evaluate the impact of $\eta$ on the numerical accuracy of Algorithm \ref{alg:dl_evi}, we test $\eta=10^{-i}, i=2, 3, 4, 5$, while keeping other parameters unchanged. The numerical errors of Algorithm \ref{alg:dl_evi} with respect to different values of $\eta$ are reported in Table \ref{tab:eta}. We can see that Algorithm \ref{alg:dl_evi} achieves consistently low errors in all cases and thus is robust to the choice of $\eta$. 
	
	\begin{table}[H]
		\centering
		\begin{tabular}{|c|c|c|c|c|}
			\hline
			$\eta$&  $10^{-2}$ & $10^{-3}$  &  $10^{-4}$ & $10^{-5}$  \\ \hline
			Relative $L^2$-errors& $5.990\times 10^{-4}$  & $5.712\times 10^{-4}$ &$5.981\times 10^{-4}$  &  $5.641\times 10^{-4}$\\ \hline
			Relative $L^{\infty}$-errors&  $1.433\times 10^{-3}$ & $1.042\times 10^{-3}$ &$1.122\times 10^{-3}$  &$1.175\times 10^{-3}$  \\ \hline
		\end{tabular}
		\caption{Numerical errors with respect to different $\eta$ for Example \ref{ex:1d_ob_sy}.}\label{tab:eta}
	\end{table}
\end{example}

\begin{example}\label{ex:1d_ob_nsy}
In this example, we test a one-dimensional problem with a non-symmetric elliptic operator as that in \cite{alphonse2024neural,zhao2022two}.
Let \(\Omega = (-2, 2)\),  the operator \(A u = -u_{xx} + u_x\), and the obstacle function \(\psi(x) = 1 - x^2\). We define $f$ and the exact solution $u$ as
\begin{equation*}
	{\footnotesize
	f(x) := 
	\begin{cases}
		(4 - 2 \sqrt{3}), & x \in [-2, -2 + \sqrt{3}), \\
		-(2 \sqrt{3} - 2), & x \in [-2 + \sqrt{3}, 2 - \sqrt{3}], \\
		-(4 - 2 \sqrt{3}), & x \in (2 - \sqrt{3}, 2],
	\end{cases}
	~\text{and}~
		u(x) = 
	\begin{cases}
		(4 - 2 \sqrt{3})(x + 2), & x\in[-2, -2 + \sqrt{3}), \\
		1 - x^2, &x\in[ -2 +\sqrt{3}, 2 -\sqrt{3}), \\
		(4 - 2 \sqrt{3})(2 - x), & x\in [2 - \sqrt{3}, 2].
	\end{cases}}
\end{equation*}

In Algorithm \ref{alg:dl_evi}, we take $h(x)=\frac{1}{4}(x+2)(2-x)$ and $\eta=10^{-3}$. The numerical results are displayed in Figure \ref{fig:1d_obstable_non}. Visually, the learned solution aligns nearly perfectly with the exact one. The results further confirm that Algorithm \ref{alg:dl_evi} achieves high precision, as evidenced by the low relative $L^2$- and $L^{\infty}$- errors. Notably, the maximum point-wise error between the exact and predicted solutions is approximately 
 $1.2\times10^{-3}$, which is smaller than the one (around $5\times 10^{-3}$) reported in \cite{alphonse2024neural} and the one (around $3\times 10^{-3}$) obtained in \cite{zhao2022two}. These findings highlight the effectiveness of Algorithm \ref{alg:dl_evi} for solving obstacle problems with non-symmetric elliptic operators. 
\begin{figure}[h!]\label{fig:1d_obstable_non}
	\centering
	\includegraphics[width=.45\textwidth]{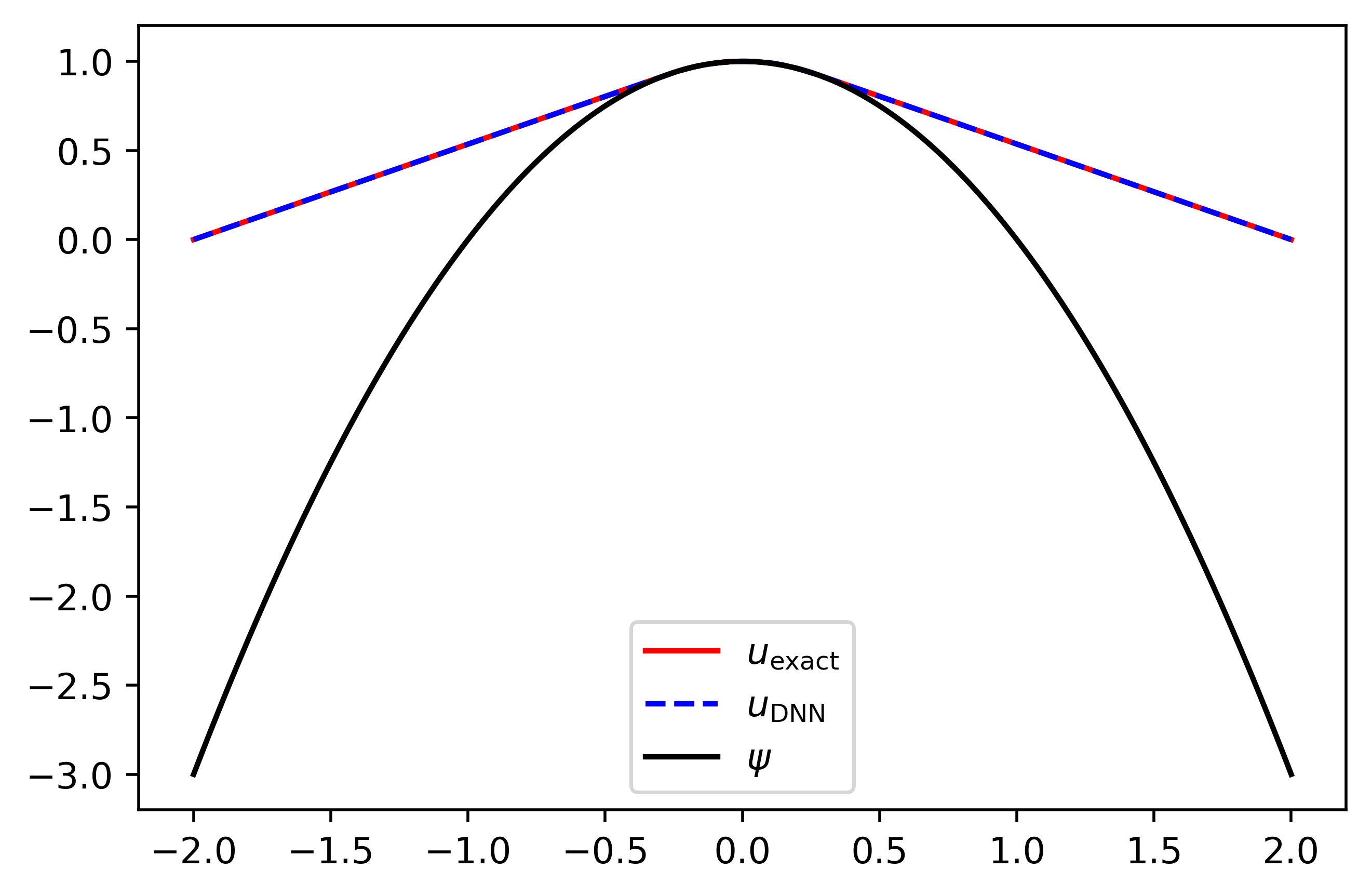}
	\includegraphics[width=.45\textwidth]{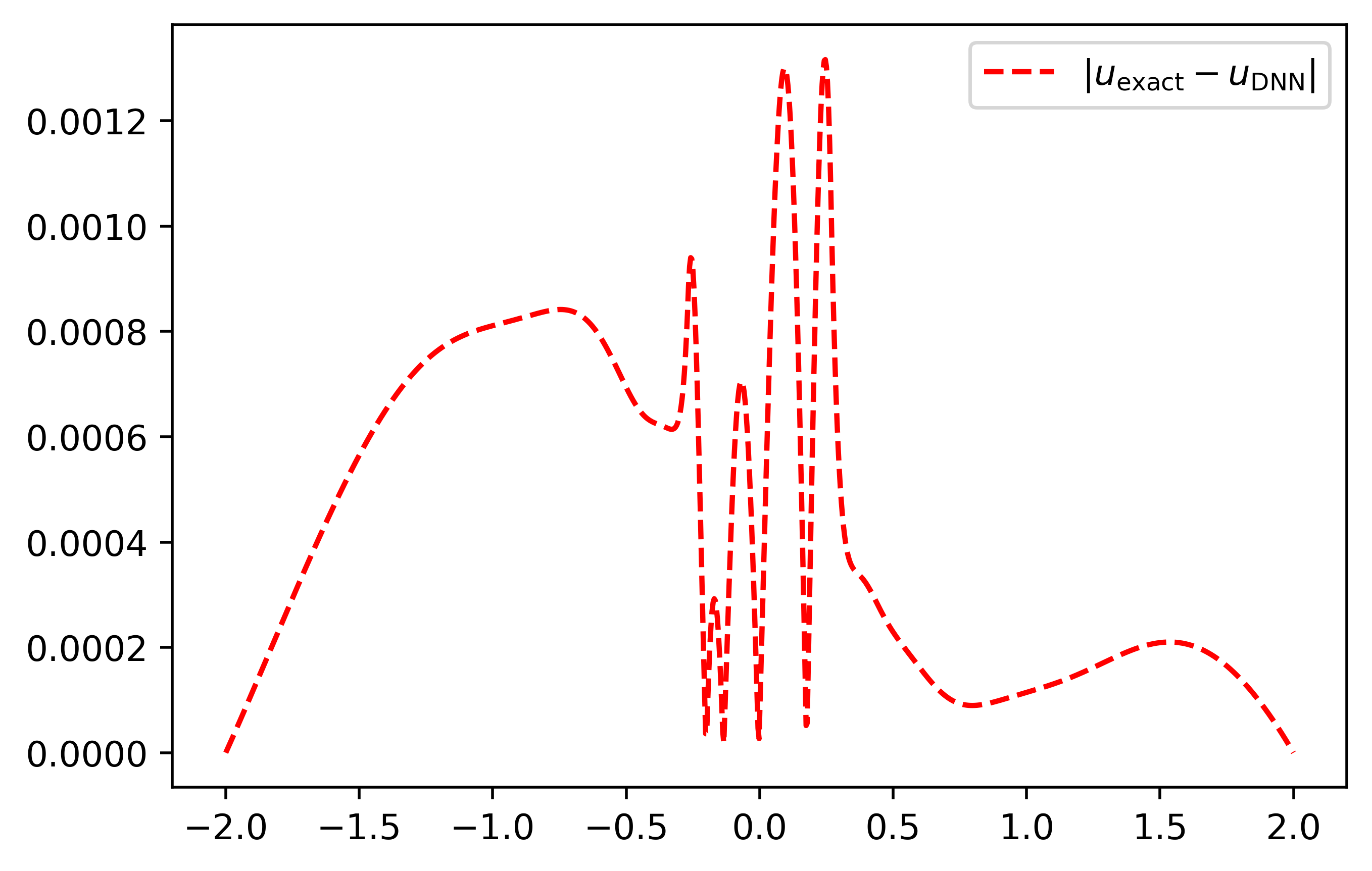}
	\includegraphics[width=.45\textwidth]{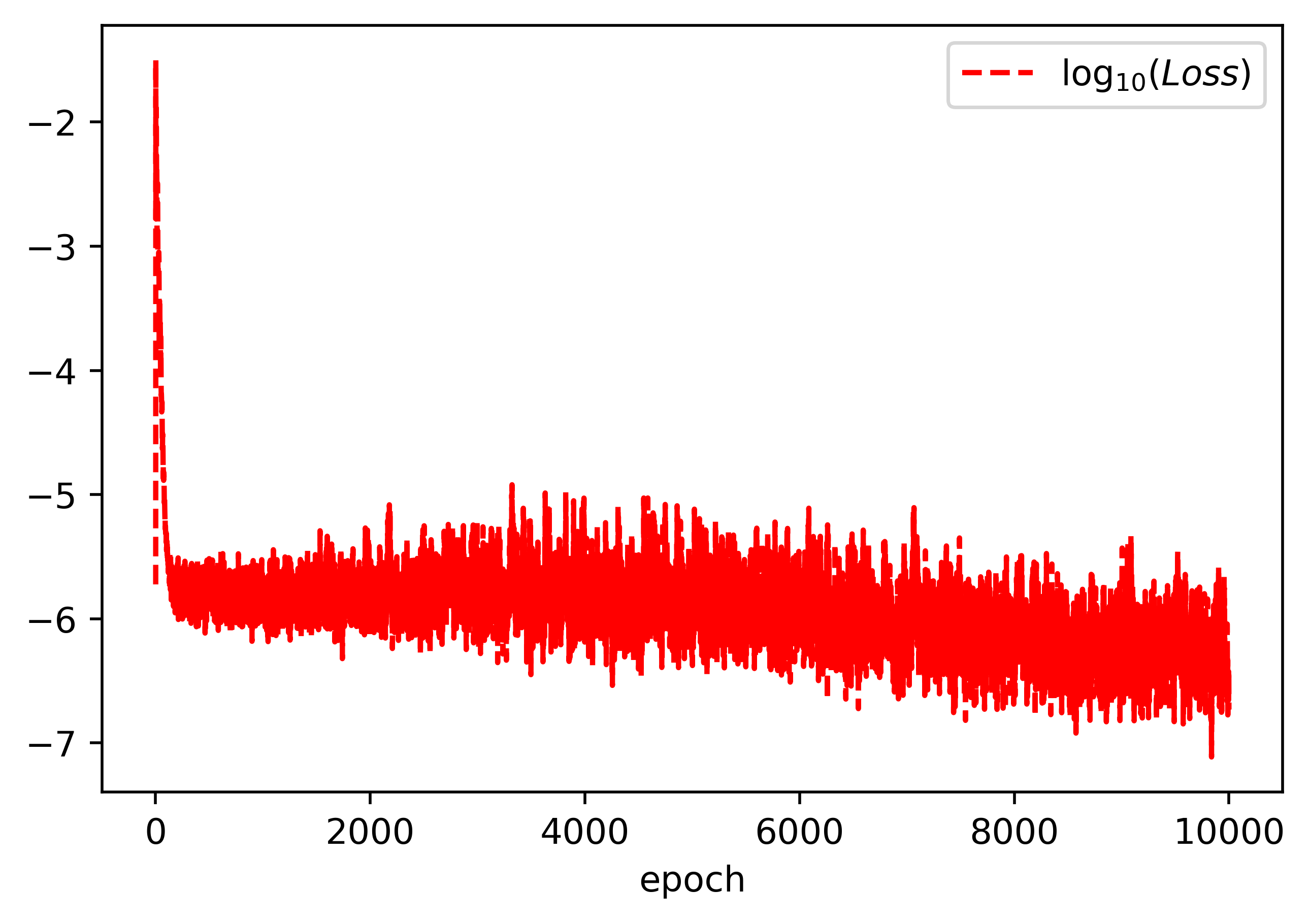}
	\includegraphics[width=.468\textwidth]{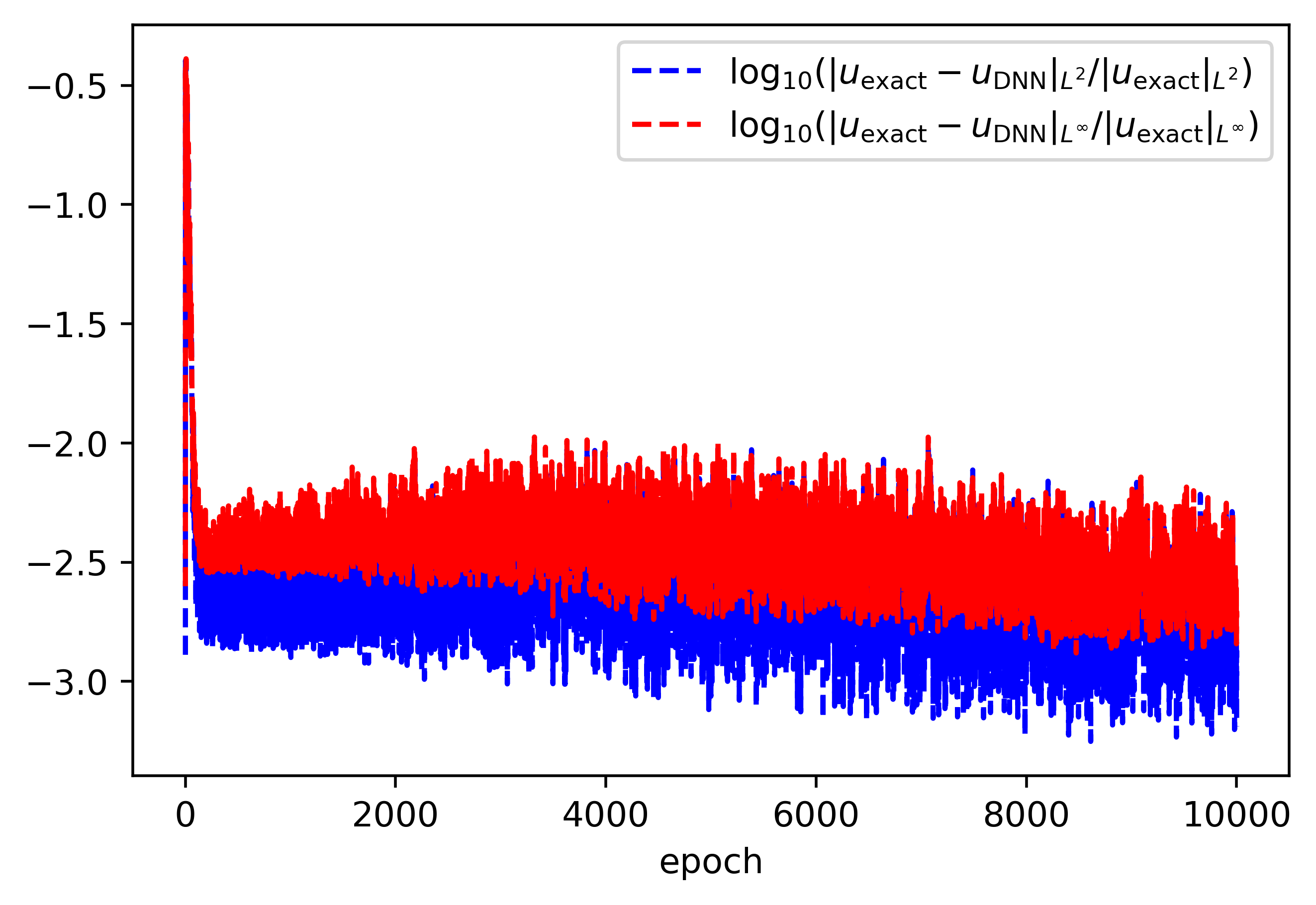}
	\caption{Numerical results for Example \ref{ex:1d_ob_nsy} ( Relative $L^2$-error: $6.472\times 10^{-4}$; Relative $L^{\infty}$-error: $1.863\times 10^{-3}$). 
	 }
\end{figure}

\end{example}

\begin{example}\label{ex:1d_ob_ps} We consider the one-dimensional obstacle problem with a piecewise smooth solution, previously investigated in \cite{alphonse2024neural}. 
	Let the domain $\Omega=(-1,1)$, the operator $A u=-u_{x x}$, and $f \equiv 0$. Let 
	\begin{equation*}
		\varphi(x):=\frac{\mu(0.4-|x|)}{\mu(|x|-0.3)+\mu(0.4-|x|)}, \quad \text { where } \quad \mu(x):= \begin{cases}\exp (-1 / x), &  x>0, \\ 
			0, &  x \leq 0.
		\end{cases}
	\end{equation*}
	Then the obstacle function is given by
	\begin{equation*}
	\psi(x):= \begin{cases}
			\begin{aligned}
			\varphi\left(x+\frac{1}{2}\right)\left(\frac{3}{2}-12\left|x+\frac{1}{2}\right|^{2-\alpha}\right)-\frac{1}{2},  &  \quad x \in(-1,0], \\ \varphi\left(x-\frac{1}{2}\right)\left(\frac{3}{2}-12\left|x-\frac{1}{2}\right|^{2-\alpha}\right)-\frac{1}{2}, &  \quad x \in(0,1),
		\end{aligned}
		\end{cases}
	\end{equation*}
	where $\alpha=0.4$. The exact solution is given by
	\begin{equation*}
		u(x)= \begin{cases}
			\begin{aligned}
				&\psi(-\beta-0.5) \frac{x+1}{0.5-\beta}, &  x \in(-1,-\beta-0.5), \\ 
					&\psi(x), &  x \in[-0.5-\beta,-0.5), \\ 
					&1, &  x \in[-0.5,0.5), \\ 
					&\psi(x), &  x \in[0.5,0.5+\beta), \\ 
					&\psi(\beta+0.5) \frac{x-1}{\beta-0.5}, &  x \in[\beta+0.5,1),
			\end{aligned}
		\end{cases}
	\end{equation*}
	where the constant $\beta$ is the unique solution of the equation
	$
		\psi(-\beta-0.5)=(0.5-\beta) \psi^{\prime}(-\beta-0.5), \beta \in(0,0.3).
	$
	In practice, we take $\beta=0.02376$ as an approximate solution of the equation.  Note that  the solution $u$ is composed
	of five separate pieces. 
	
	For this example, we take $h(x)=(x+1)(1-x)$ and $\eta=10^{-3}$ in the implementation of Algorithm \ref{alg:dl_evi}. The numerical results, summarized in Figure \ref{fig:1d-obstacle_ps}, include the exact and the learned solutions, the point-wise error, the training trajectories for the loss function, and the test errors versus training epochs.  We observe that the numerical solution is a good approximation to the exact one. Moreover, the results validate that Algorithm \ref{alg:dl_evi} can obtain accurate predictions with small relative $L^2$- and $L^{\infty}$- errors. Specifically, the maximum discrepancy between the learned and exact solutions is approximately $4\times10^{-3}$, which is significantly smaller than the one (approximately $2\times 10^{-2}$) reported in \cite{alphonse2024neural}. This result indicates the superiority of Algorithm \ref{alg:dl_evi} in terms of numerical accuracy. 
	
	\begin{figure}[h!]\label{fig:1d-obstacle_ps}
		\centering
		\includegraphics[width=.45\textwidth]{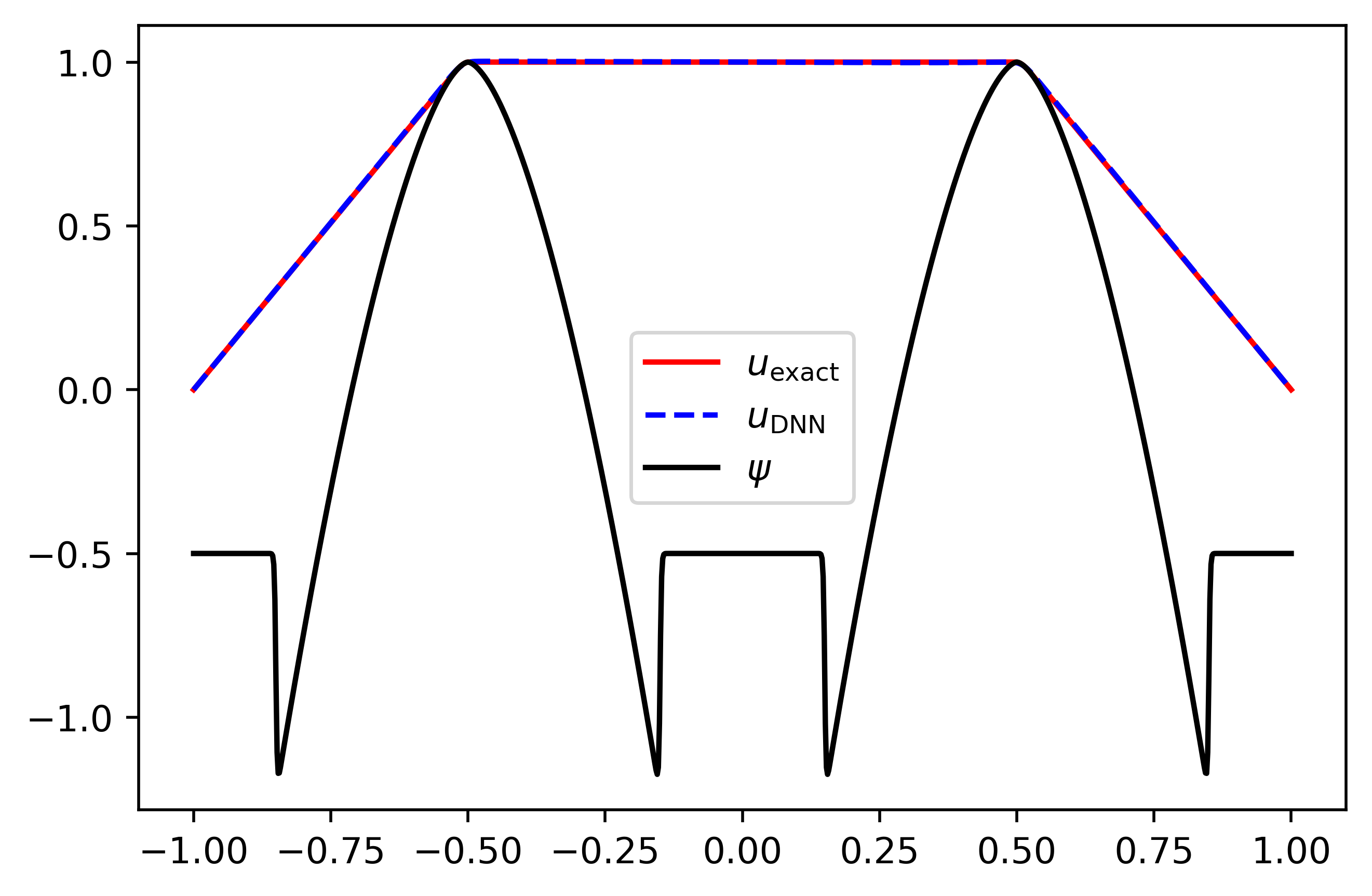}
		\includegraphics[width=.45\textwidth]{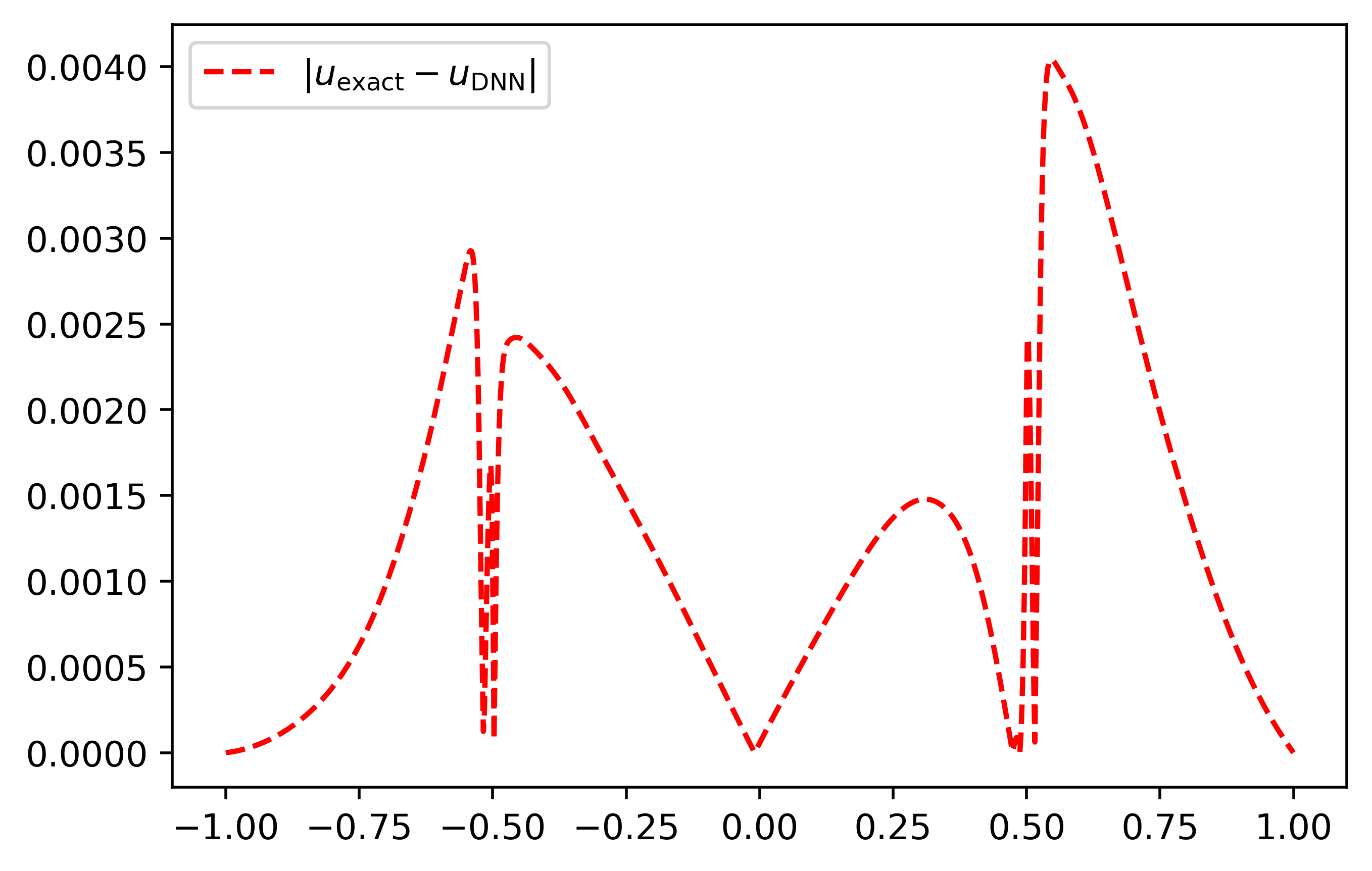}
		\includegraphics[width=.45\textwidth]{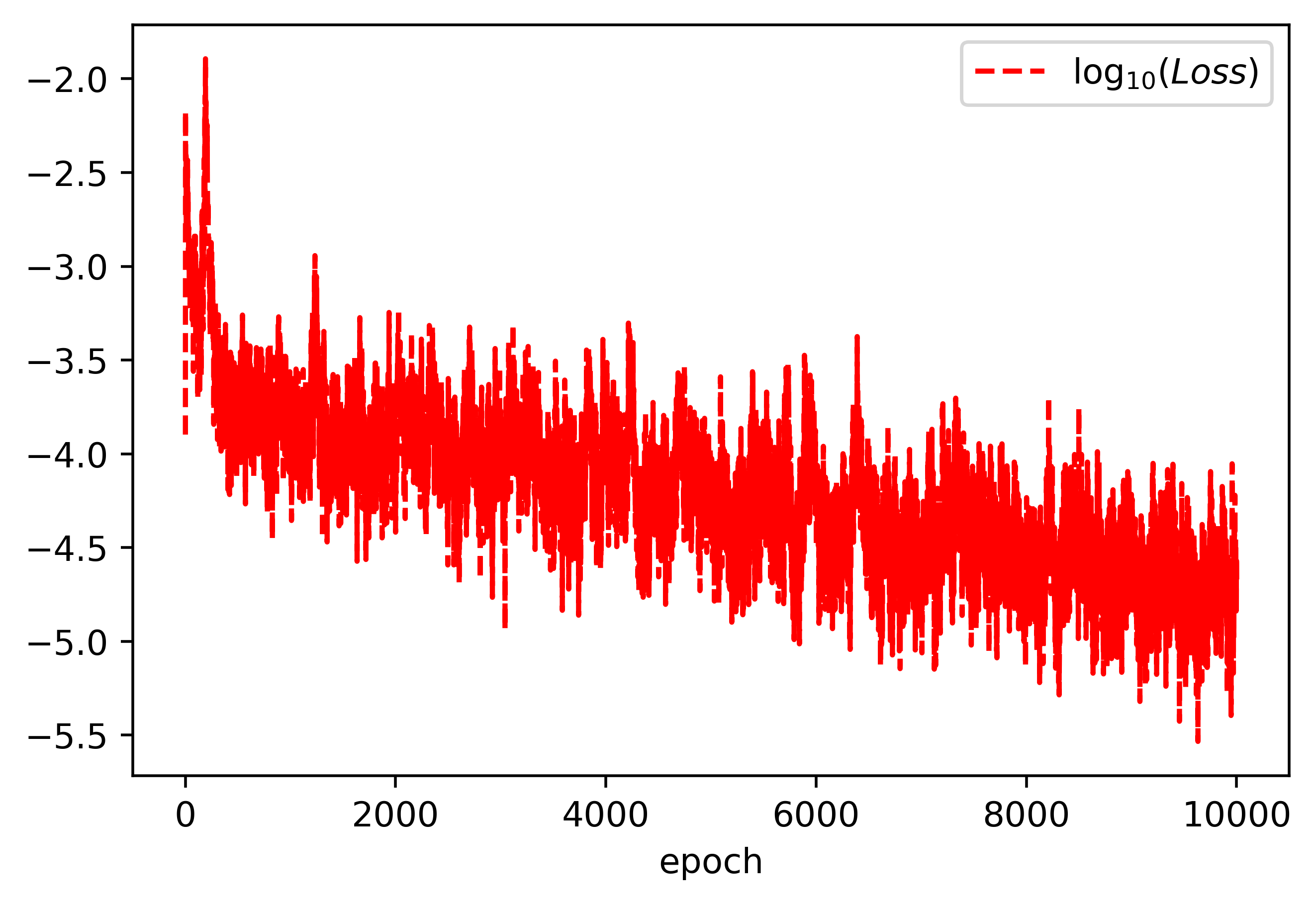}
		\includegraphics[width=.46\textwidth]{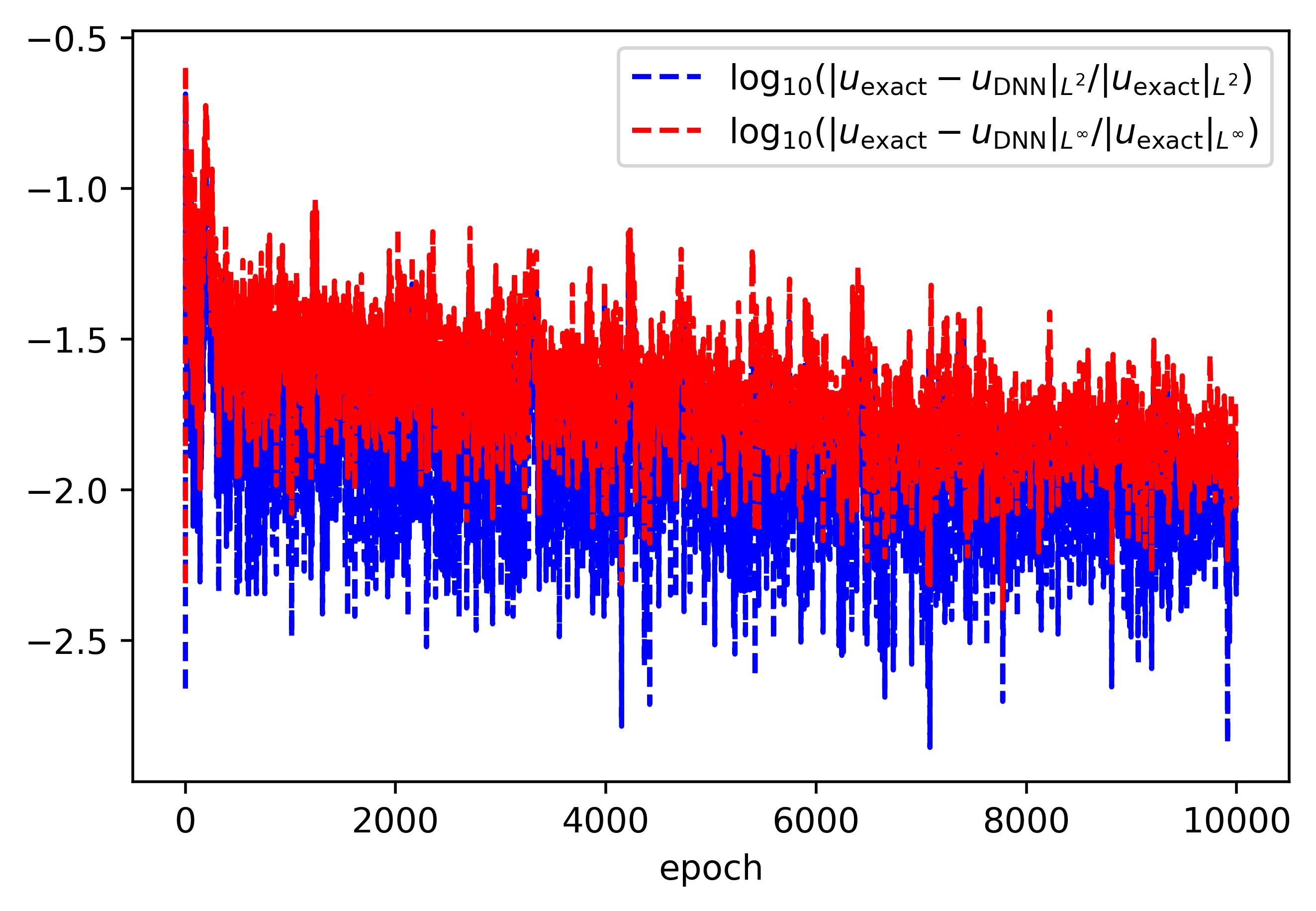}
		\caption{Numerical results for Example \ref{ex:1d_ob_ps} (Relative $L^2$-error: $5.045\times 10^{-3}$; Relative $L^{\infty}$-error: $9.167\times 10^{-3}$). 
			}
	\end{figure}

\end{example}

\begin{example}\label{ex:2d_ob}
	We test the two-dimensional problem with a piecewise smooth solution considered in \cite{cheng2023deep}. 
	Let $\Omega=(-2,2)^2$, and the operator $Au:=-\Delta u$, and $f=0$. The obstacle function \(\psi(x)\) and the exact solution \(u(x)\) are defined as follows:
	\begin{equation*}
		\psi(x)= \begin{cases}
			\sqrt{1-|x|^2}, & |x| \leq 1, \\ 
			-1, & \text { else where },
			\end{cases}
			~\text{and}~
					u(x)= \begin{cases}
				\sqrt{1-|x|^2}, & |x|^2 \leq r^*, \\ 
				-\left(r^*\right)^2 \ln (|x| / 2) / \sqrt{1-\left(r^*\right)^2}, & |x| \geq r^*,
			\end{cases}
	\end{equation*}
	where $x=(x_1,x_2)\in \overline{\Omega}$, $|x|=\sqrt{x_1^2+x_2^2}$, and $r^*$ satisfies $\left(r^*\right)^2\left(1-\ln \left(r^* / 2\right)\right)=1$. Here, we take $r^* \approx 0.6979651482$.
	
	Note that $u(x)\neq 0$ on $\partial\Omega$. Hence, to enforce the boundary condition, we construct a neural network in the form of \eqref{eq: bd_non} to approximate $u$. Next, we discuss the choices of $g$ and $h$. To this end, we let $a=c=-2$, $b=d=2$, and define
	\begin{equation*}
		 w(x_1)=\frac{x_1-a}{b-a} ~\text{and}~ w(x_2)=\frac{x_2-c}{d-c},
	\end{equation*}
    which satisfy 
    $w(a)=0, ~ w(b)=1, ~w(x_1) \in[0,1],~\text{and}~w(c)=0, ~ w(d)=1, ~w(x_2) \in[0,1].$
   Then, we let 
   	\begin{equation*}
   	\begin{aligned}
   		g(x_1, x_2)=&[1-w(x_1)]  u(a, x_2)+w(x_1)  u(b, x_2)+[1-w(x_2)]  u(x_1, c) +w(x_2)  u(x_1,d)\\
   		&-\Big\{[1-w(x_1)] [1-w(x_2)]  u(a, c)+[1-w(x_1)] w(x_2) u(a, d)\\
   		&\quad+
   		w(x_1) [1-w(x_2)]  u(b, c)+w(x_1)w(x_2)u(b, d)\Big\}.
   	\end{aligned}
   \end{equation*}
   It is straightforward to verify that $g\in C(\bar{\Omega})$ and
	\begin{equation*}
		g(x)=u(x),\quad \forall x=(x_1,x_2)\in\partial\Omega. 
	\end{equation*}
   For the choice of $h$, we define
   $\hat{h}(x_1,x_2)= (x_1-a)(b-x_1)(x_2-c)(d-x_2)$
   and then take 
   \begin{equation*}
   	 {h}(x_1,x_2)  = \hat{h}(x_1,x_2)/\|\hat{h}\|_{L^\infty} ~\text{with}~ \|\hat{h}\|_{L^\infty}=\frac{(b-a)^2(d-c)^2}{16}. 
   \end{equation*}

	The numerical results of Algorithm \ref{alg:dl_evi}, with the above constructed $g$, $h$, and $\eta=10^{-3}$, are presented in Figure \ref{fig:2d_obstacle}, which includes the exact and learned solutions, the point-wise error, the training trajectories for the loss function, and the test errors with respect to training epochs.  We observe that the numerical solution is in
	good agreement with the exact one. The results show that Algorithm \ref{alg:dl_evi} can obtain accurate predictions with low relative $L^2$-error $1.810\times 10^{-3}$ and $L^{\infty}$-error $5.523\times10^{-3}$. In particular, the maximal point-wise error between the learned and exact solutions is about $4.5\times10^{-3}$, which is much smaller than the one (around $1\times 10^{-2}$) reported in \cite{cheng2023deep}.
	\begin{figure}[h!]\label{fig:2d_obstacle}
		\centering
		\includegraphics[width=.32\textwidth]{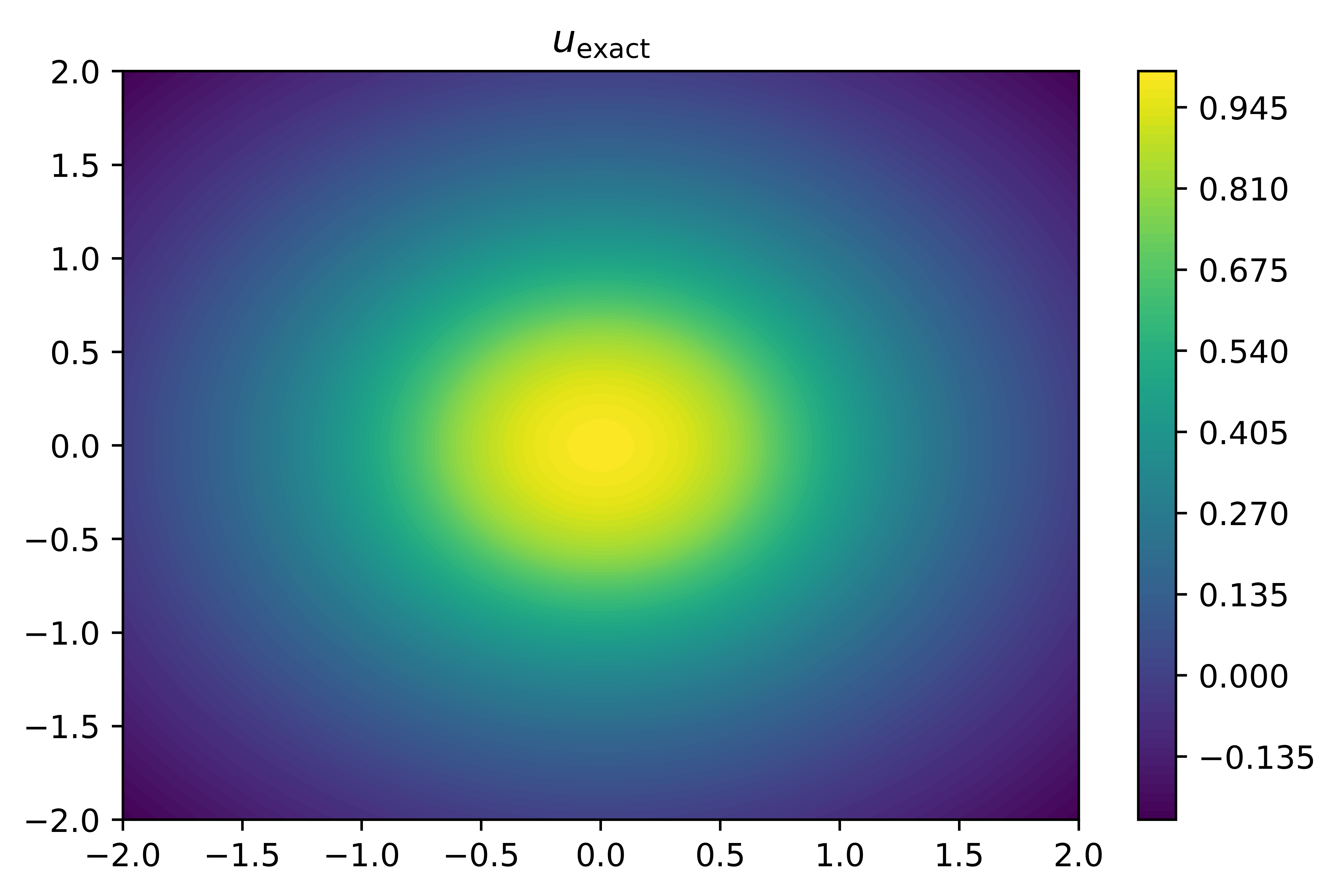}
		\includegraphics[width=.32\textwidth]{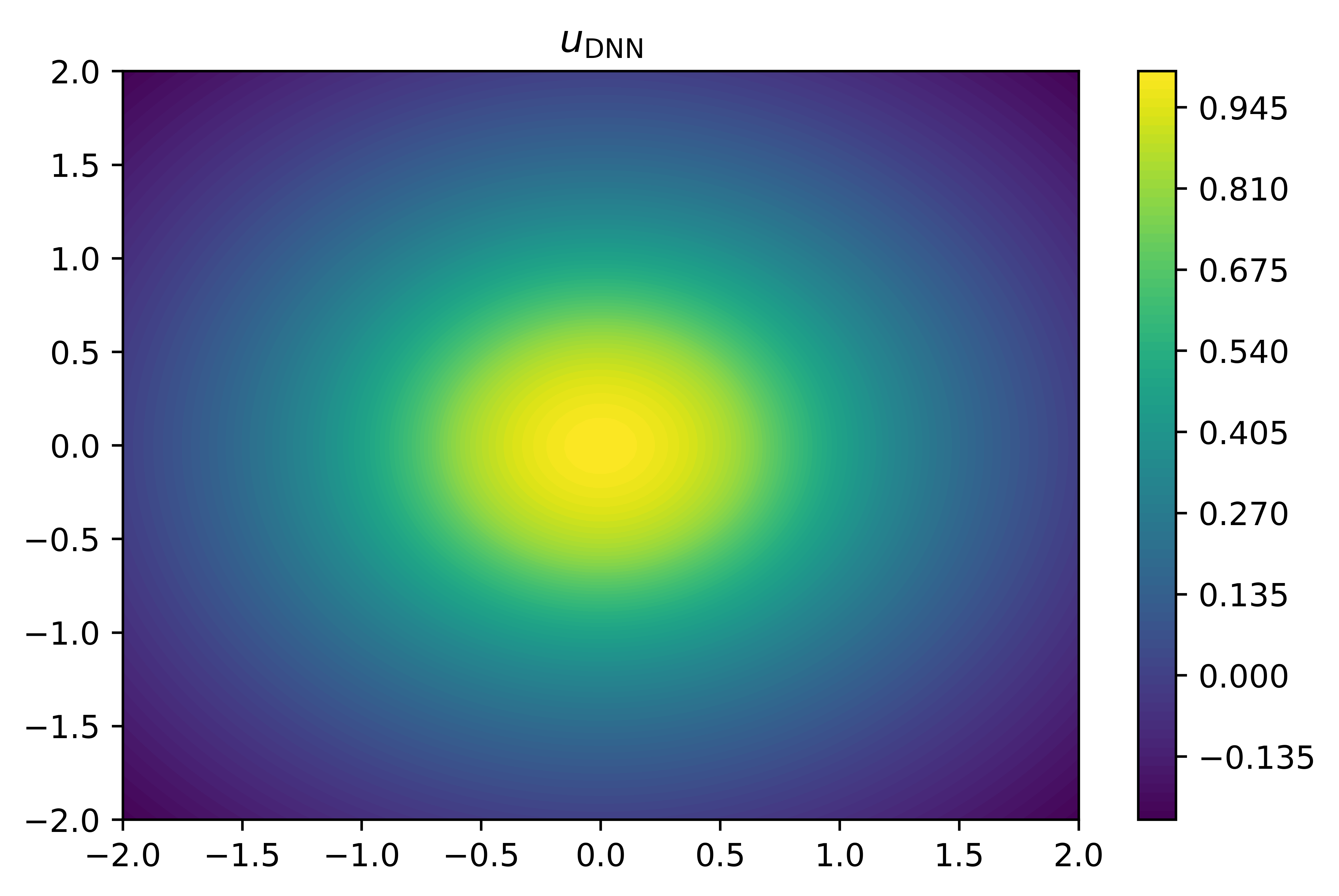}
		\includegraphics[width=.32\textwidth]{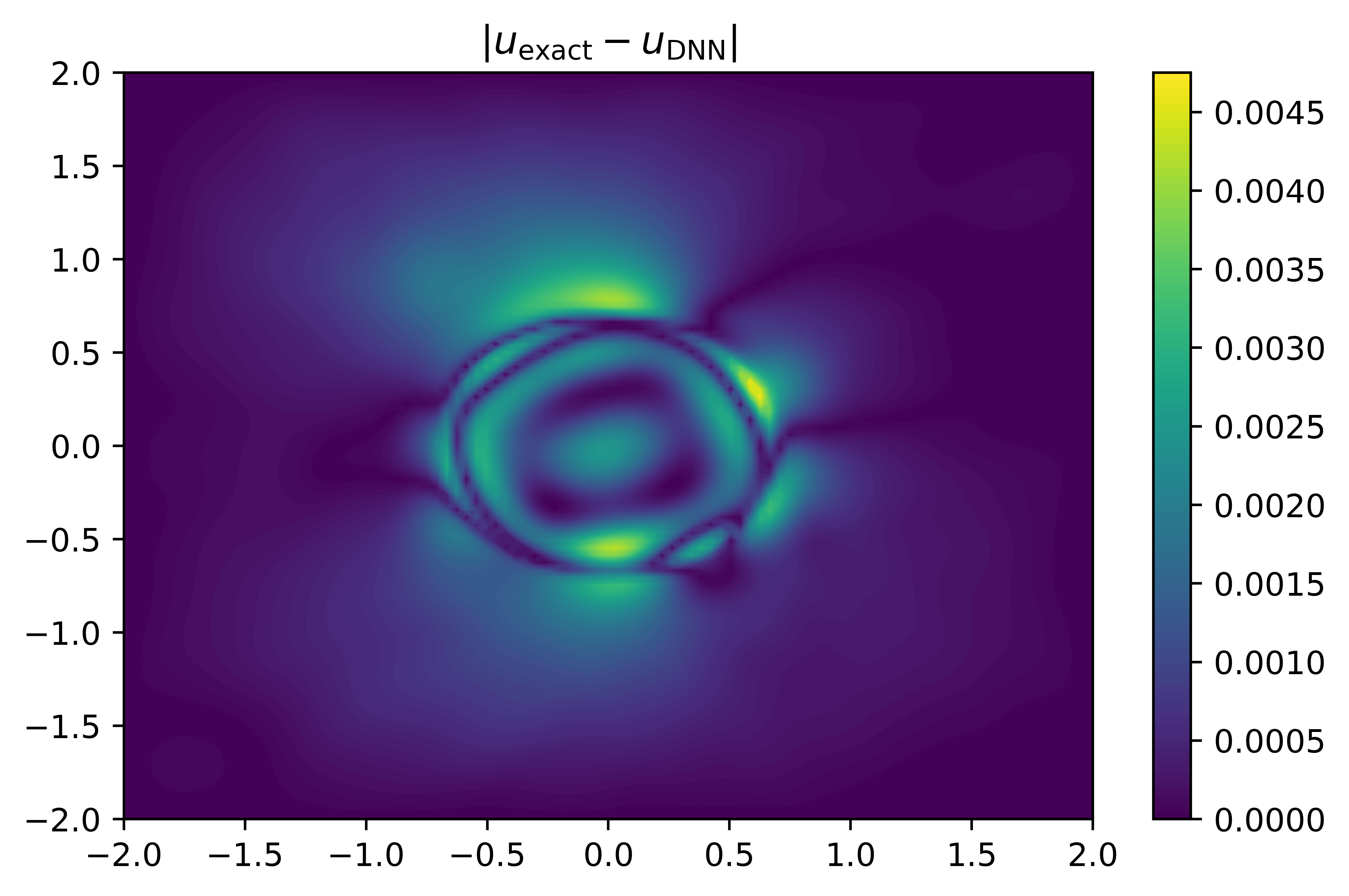}
			\includegraphics[width=.45\textwidth]{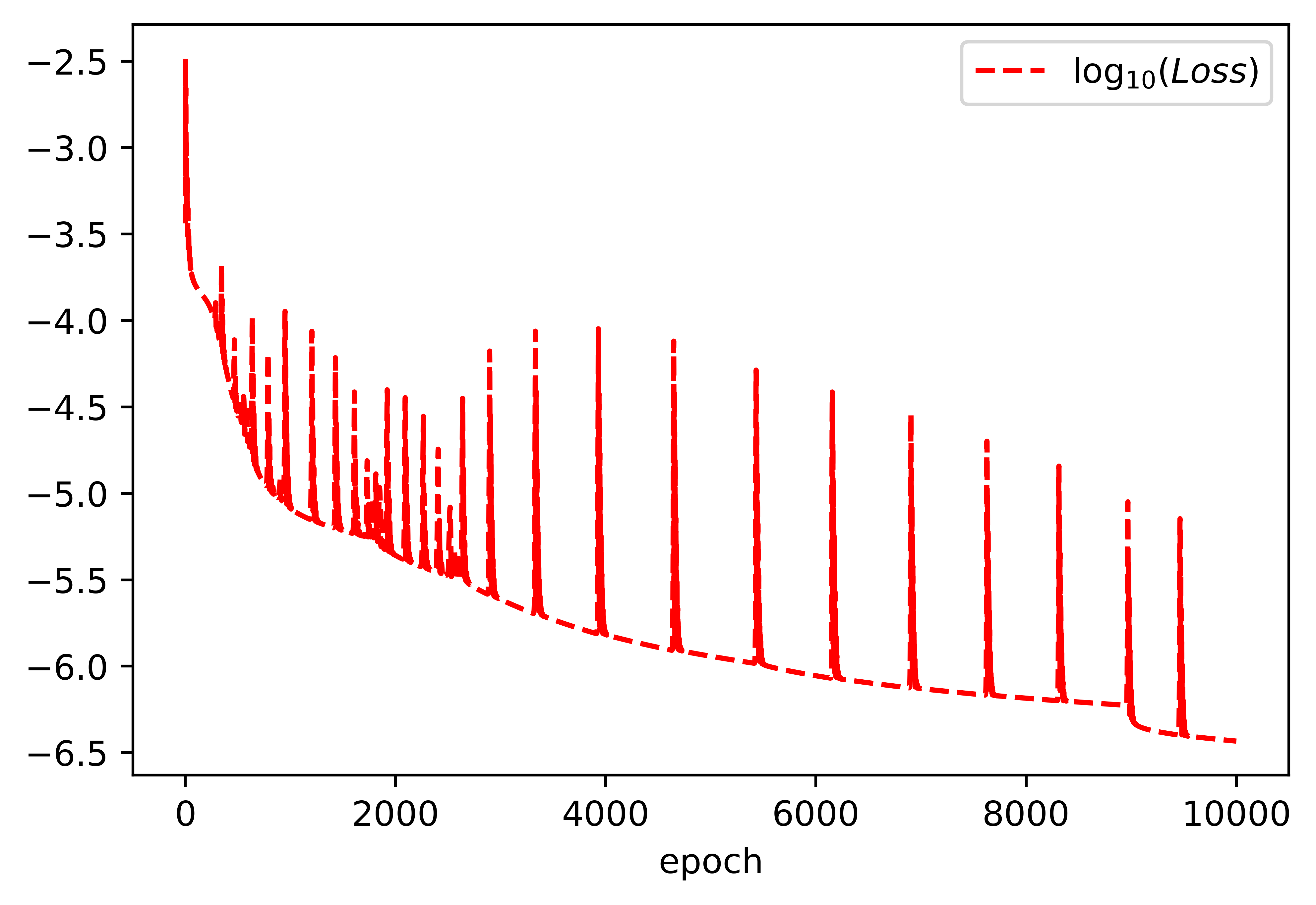}
		\includegraphics[width=.46\textwidth]{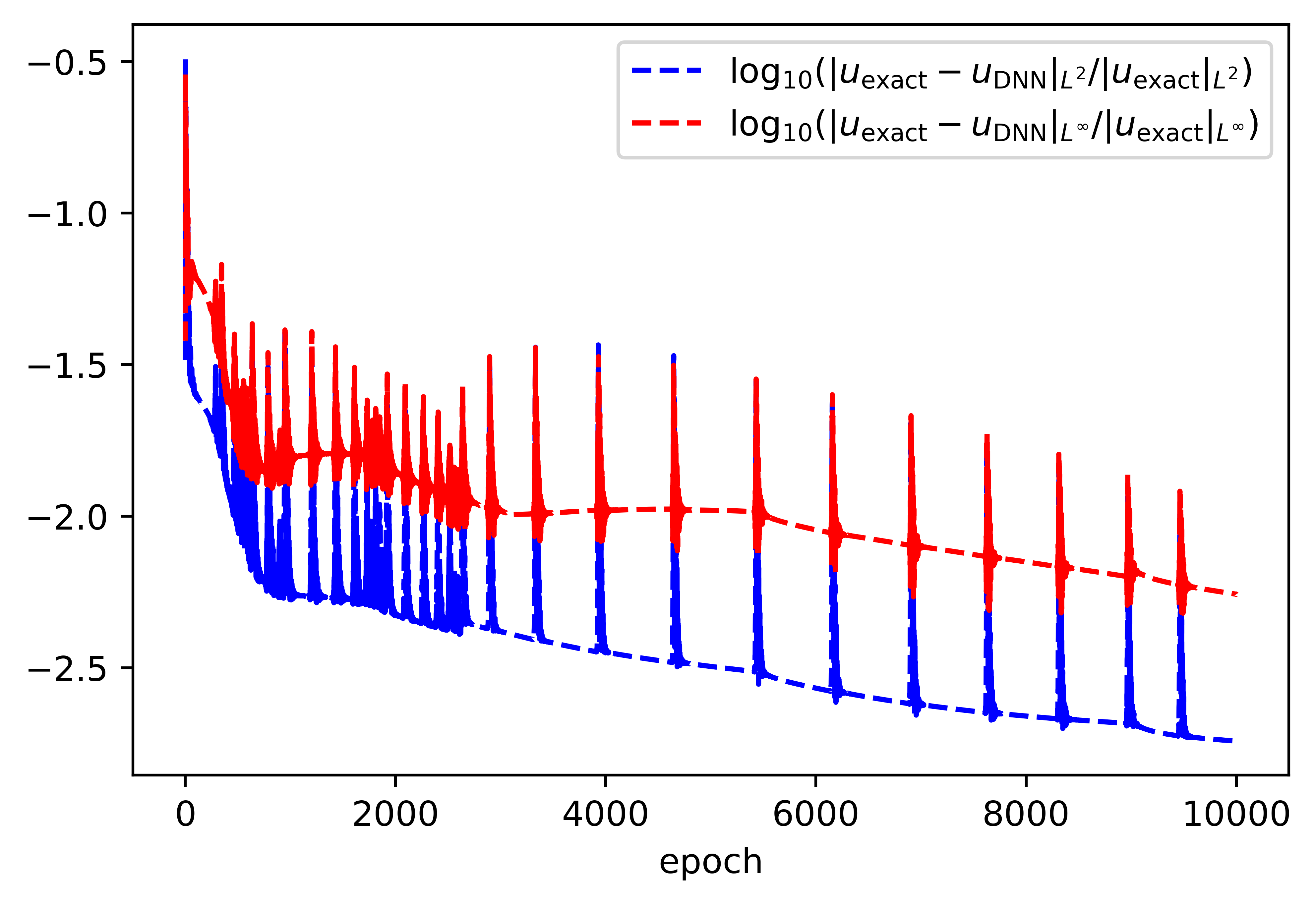}
		\caption{Numerical results for Example \ref{ex:2d_ob} (Relative $L^2$-error: $1.810\times 10^{-3}$; Relative $L^{\infty}$-error: $5.523\times 10^{-3}$). 
			}
	\end{figure}
	
\end{example}

\subsection{Elasto-Plastic Torsion Problems}

Let us consider an infinitely long cylindrical bar of cross-section $\Omega$, with $\Omega$ being bounded and simply connected. Assume that the bar is made of an isotropic elastic perfectly plastic material whose plasticity yield is given by the Von Mises criterion. Starting from a zero-stress initial state, an increasing torsion moment is applied to the bar. The torsion is characterized by $f$ (often set as a constraint), which is the torsion angle per unit length. Then, for all $f$, it follows from the Haar–K{\'{a}}rm{\'{a}}n principle that the determination of the stress field is equivalent (in a convenient system of physical units) to the solution of the following EVI:
\begin{equation}\label{elasto-plasticVI}
	u\in {K},\quad\text{such that}\quad \int_\Omega \nabla u\cdot\nabla (v-u) dx\geq \int_\Omega f(v-u) dx, \forall v\in {K},
\end{equation}
where 
	${K}=\{v| v\in H_0^1(\Omega), |\nabla v(x)|\leq 1,~ \text{a.e. in}~ \Omega\}.$
The existence and regularity of the solution of  \eqref{elasto-plasticVI} has been studied in \cite{Glowinski1984Numerical,glowinski2008lectures}. Moreover,  the discretization together with the iterative algorithms for solving \eqref{elasto-plasticVI} can be found in \cite{Glowinski1984Numerical,glowinski2008lectures,glowinskinumerical1981}.

\begin{example}\label{ex:2d_EPT}
	We test a two-dimensional problem constructed in \cite{Glowinski1984Numerical}.
	Let the domain $\Omega=\{x=(x_1, x_2) \mid |x|:=\sqrt{x_1^2+x_2^2}<R\}$ and $f(x)=c$, where $R$ and $c$ are given constants. The exact solution \(u(x)\) are defined as follows:
	  \begin{equation*}
	  	{\footnotesize
		\text { if } c R\leq 2, ~ u(x) = \frac{c}{4}\left(R^2-|x|^2\right); ~	\text { if } cR>2,~ u(x) =
		\begin{cases}
			\begin{aligned}
				& R-|x|, & \frac{2}{c} \leq |x| \leq R, \\
				& \frac{c}{4}\left[\left(R^2-|x|^2\right)-\left(R-\frac{2}{c}\right)^2\right], &0 \leq |x| \leq \frac{2}{c}. \\
			\end{aligned}
		\end{cases}}
	\end{equation*}
	
	  Note that the exact solution $u$ is determined by the constants $c$ and $R$. In our numerical experiments, we fix $R=1$ and test the example with $c=1$ and $c=4$. To rigorously enforce the boundary condition, we take $h(x)=R^2-(x_1^2+x_2^2)$ to implement Algorithm \ref{alg:dl_evi} with the loss function given in \eqref{eq:loss_3}. The numerical results for $c=1$ and $c=4$ are respectively presented in Figure \ref{fig:2d_elasto_plastic_1} and Figure \ref{fig:2d_elasto_plastic_4}. The results indicate that the numerical solutions for both configurations align closely with the exact ones. Specifically, the maximal point-wise errors between the learned and exact solutions reach magnitudes on the order of $10^{-3}$ and $10^{-4}$, which, together with the low relative $L^2$-errors and $L^{\infty}$-errors, validate that Algorithm \ref{alg:dl_evi} can produce solutions with high accuracy for this two-dimensional elasto-plastic torsion problem. 
	
			\begin{figure}[h!]
				\centering
				\includegraphics[width=.32\textwidth]{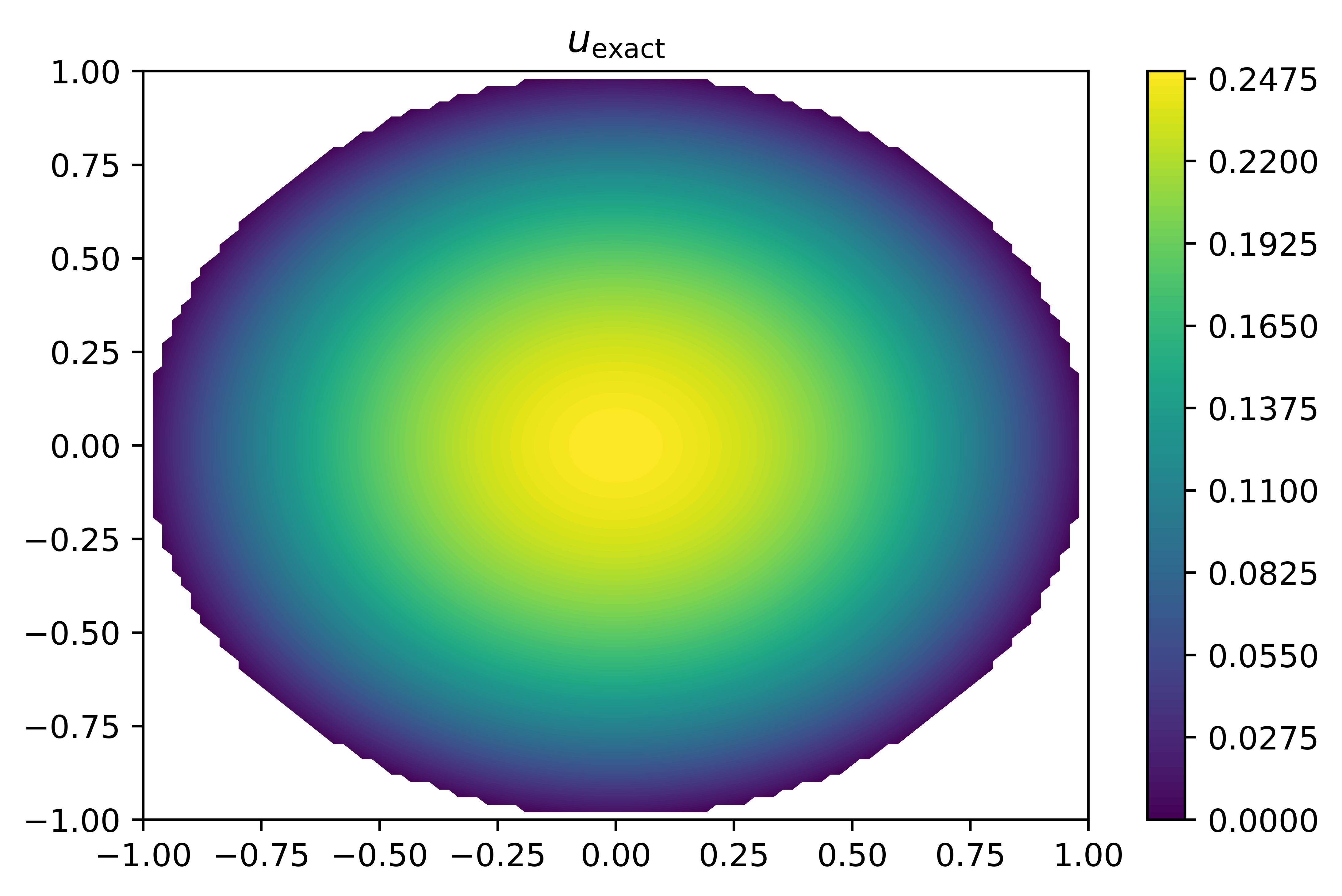}
				\includegraphics[width=.32\textwidth]{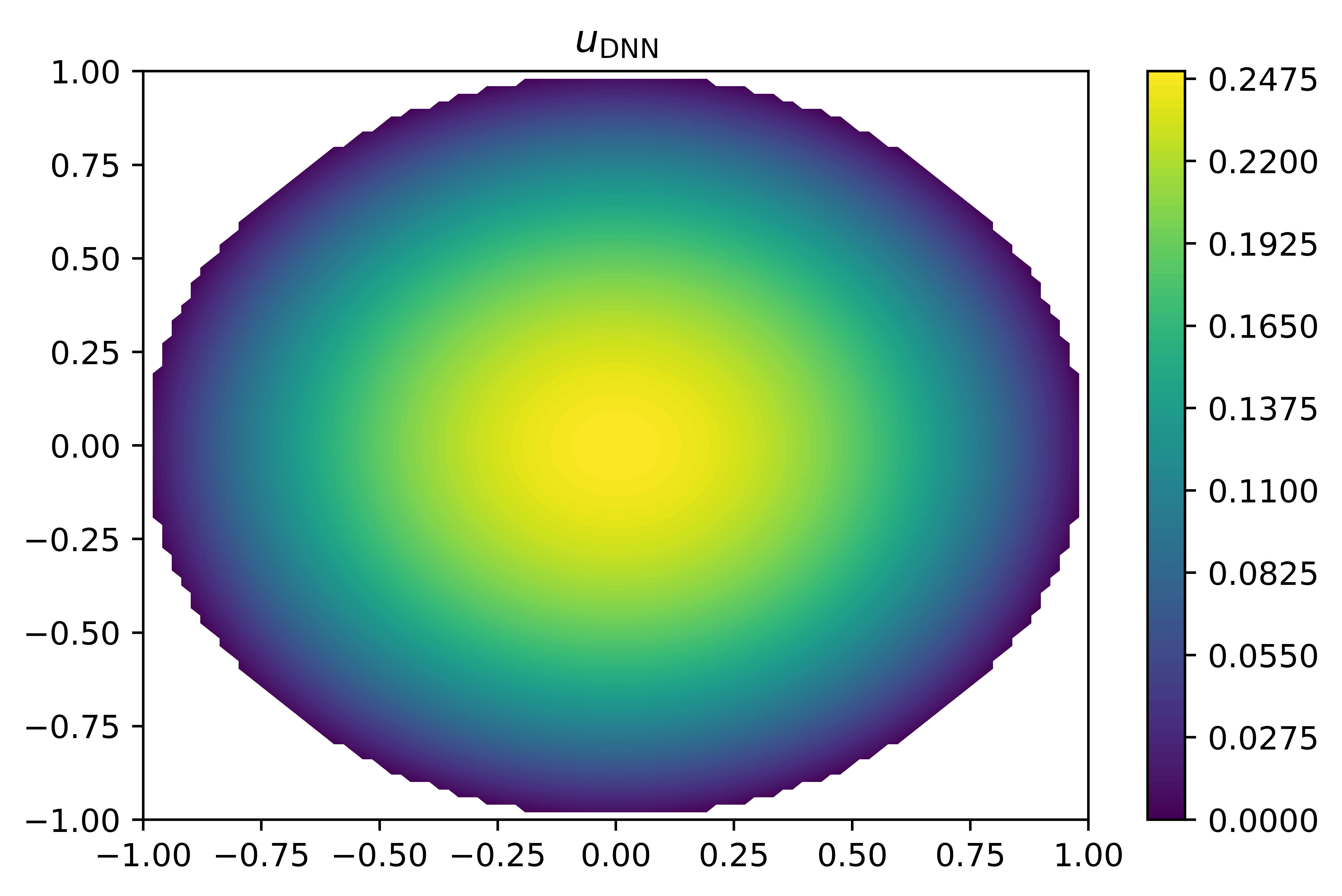}
				\includegraphics[width=.32\textwidth]{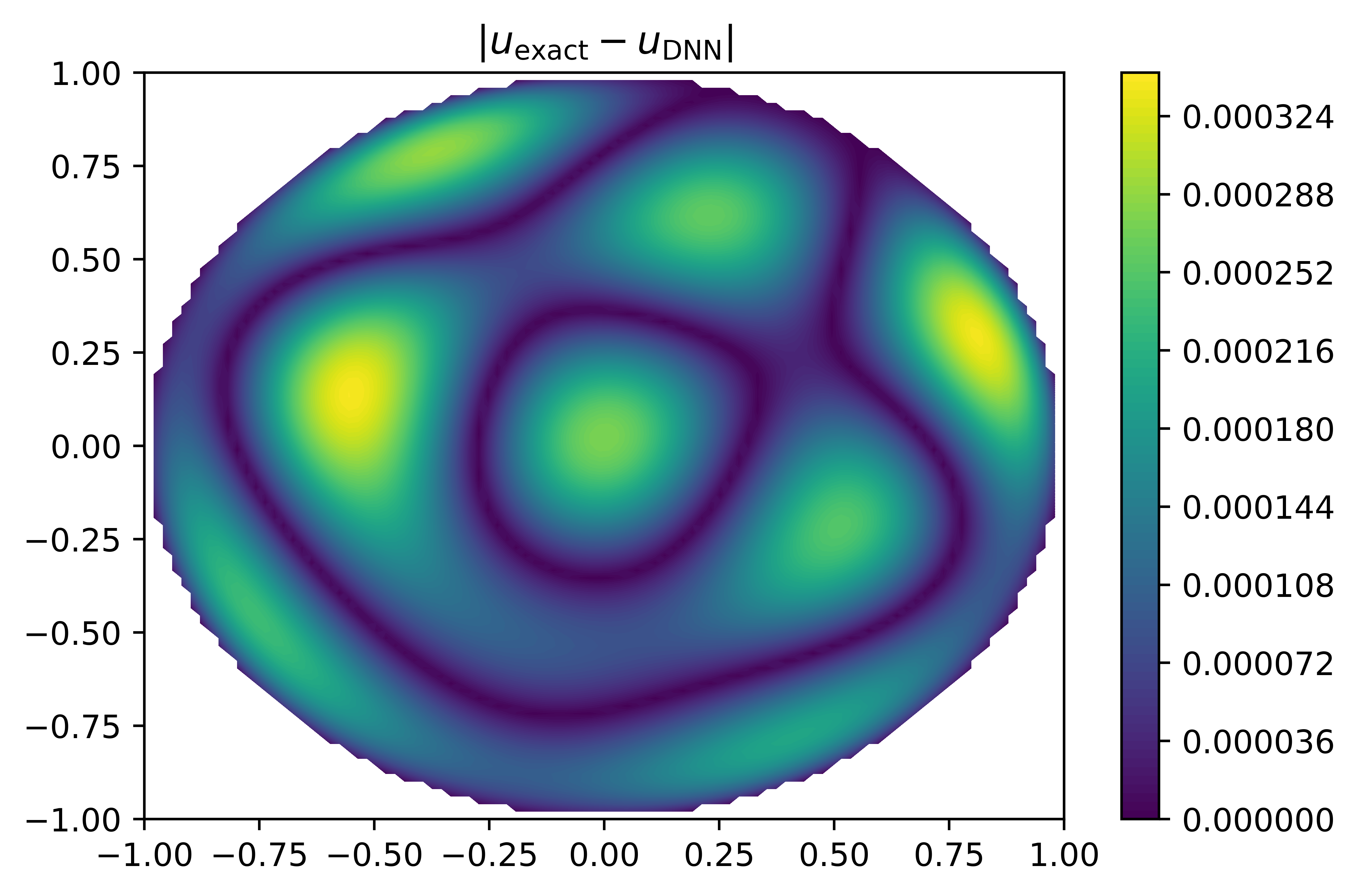}
				\includegraphics[width=.45\textwidth]{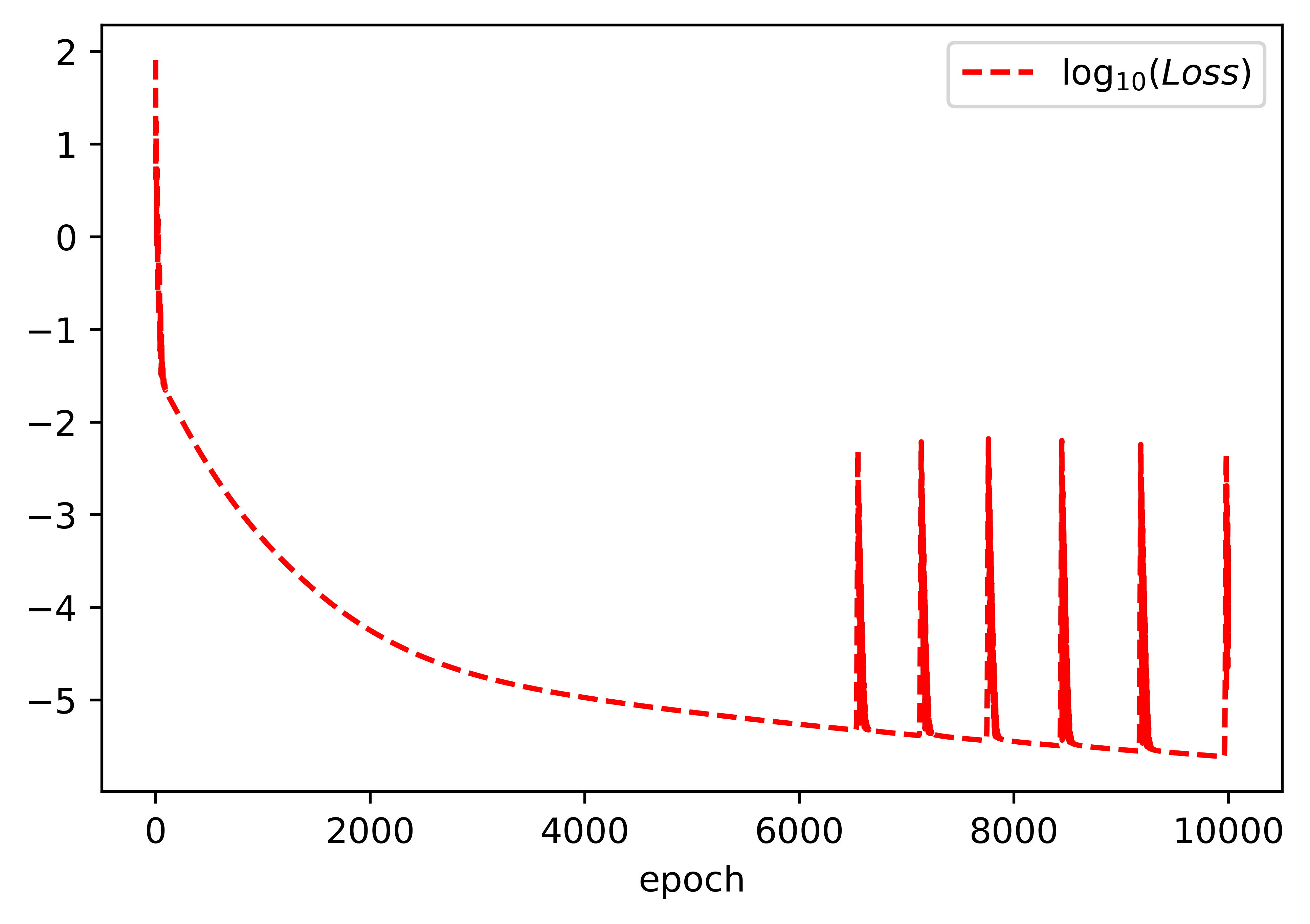}
				\includegraphics[width=.46\textwidth]{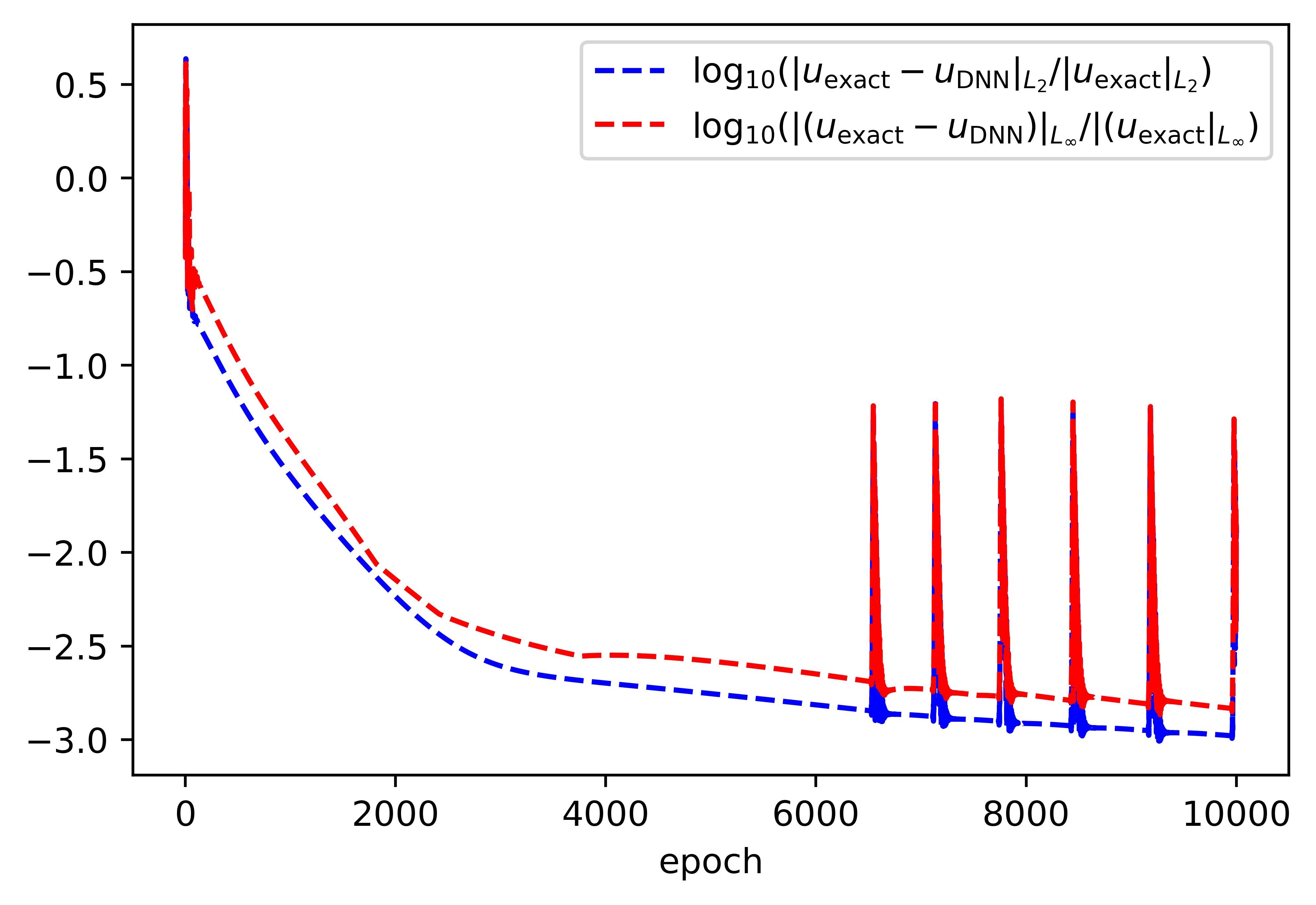}
				\caption{Numerical results for Example \ref{ex:2d_EPT} with $c=1, R=1$  (Relative $L^2$-error: $1.056\times 10^{-3}$; Relative $L^{\infty}$-error: $1.481\times 10^{-2}$). 
				}\label{fig:2d_elasto_plastic_1}
			\end{figure}
	
			\begin{figure}[h!]
			\centering
			\includegraphics[width=.32\textwidth]{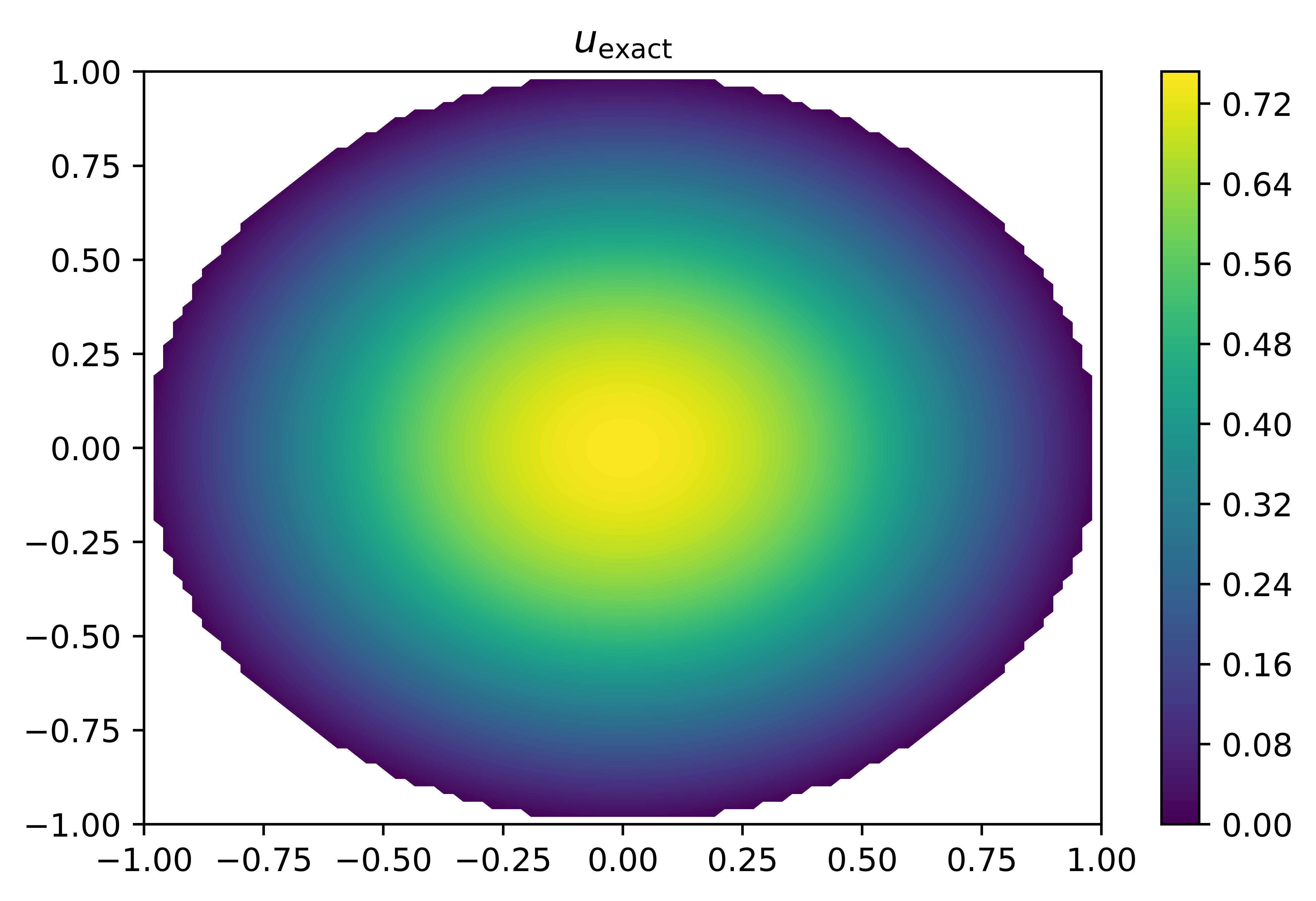}
			\includegraphics[width=.32\textwidth]{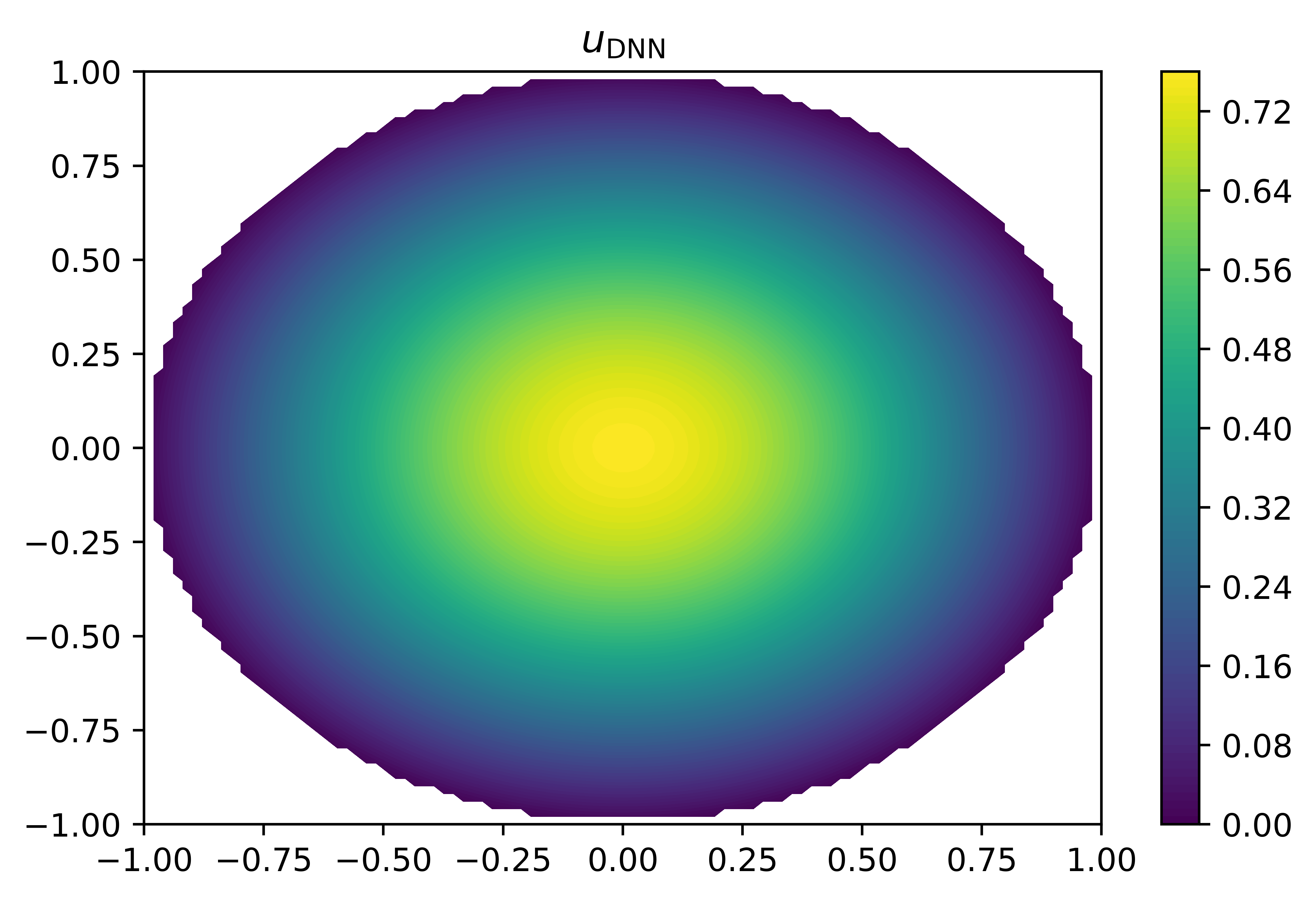}
			\includegraphics[width=.32\textwidth]{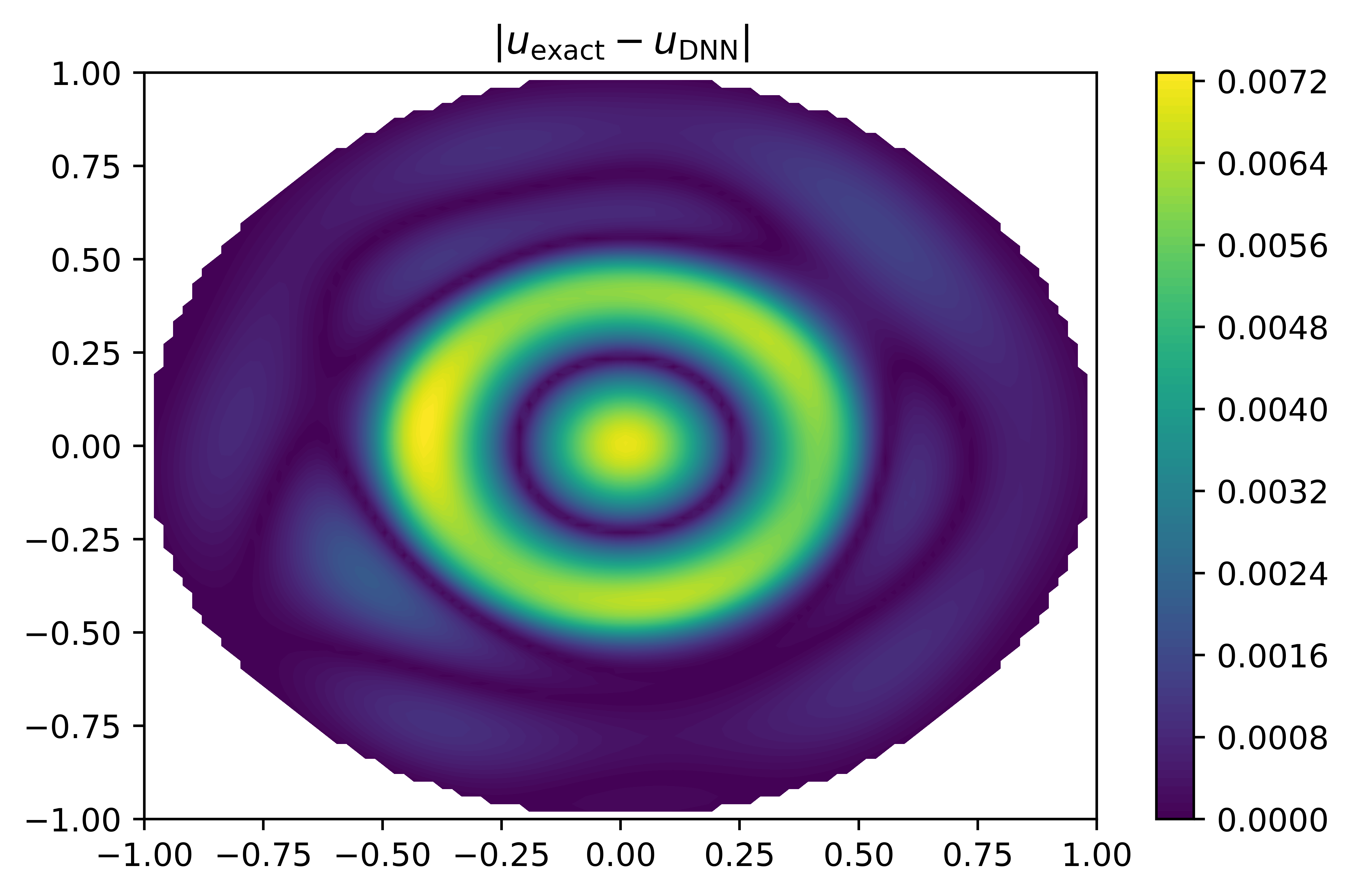}
				\includegraphics[width=.45\textwidth]{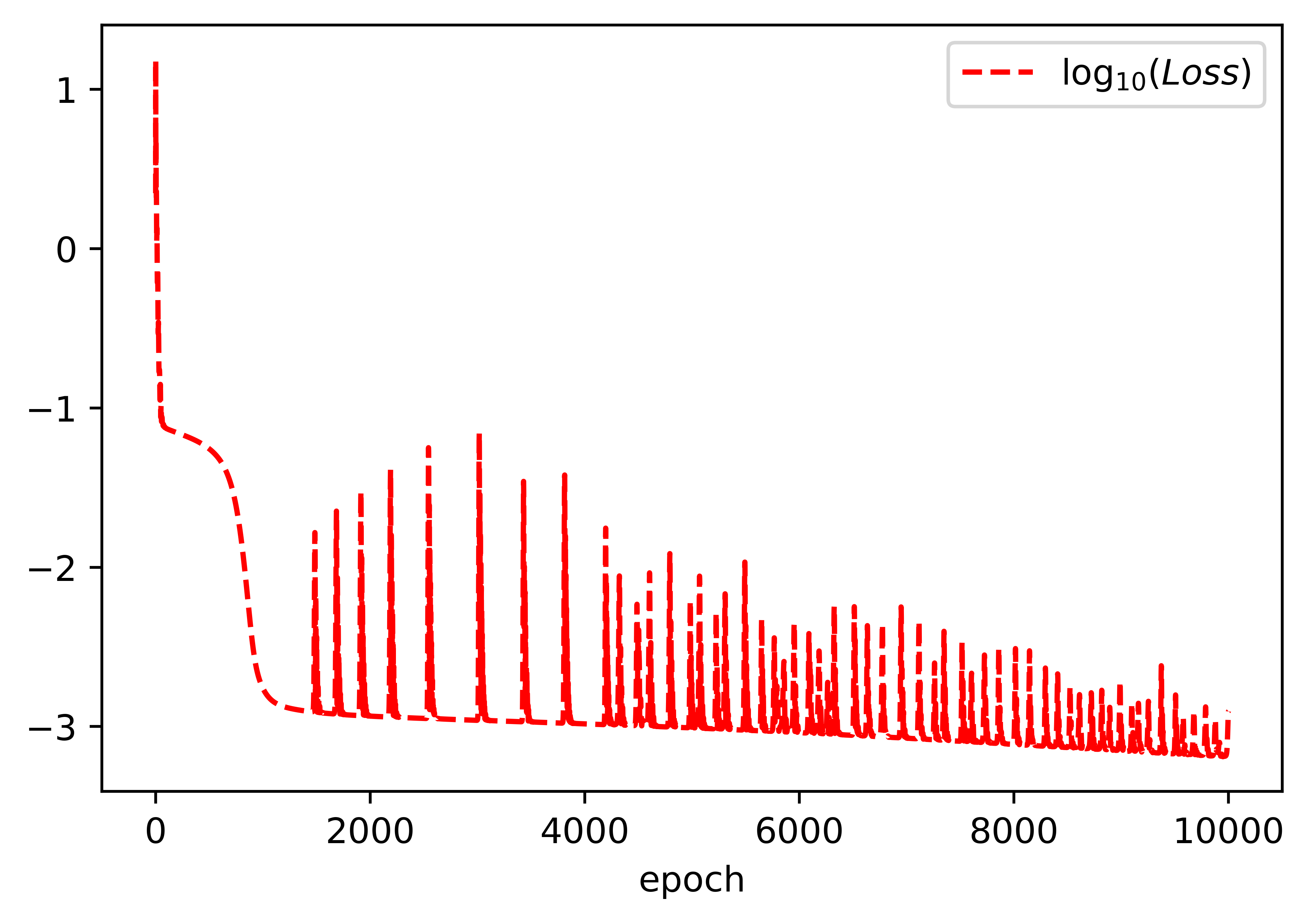}
			\includegraphics[width=.46\textwidth]{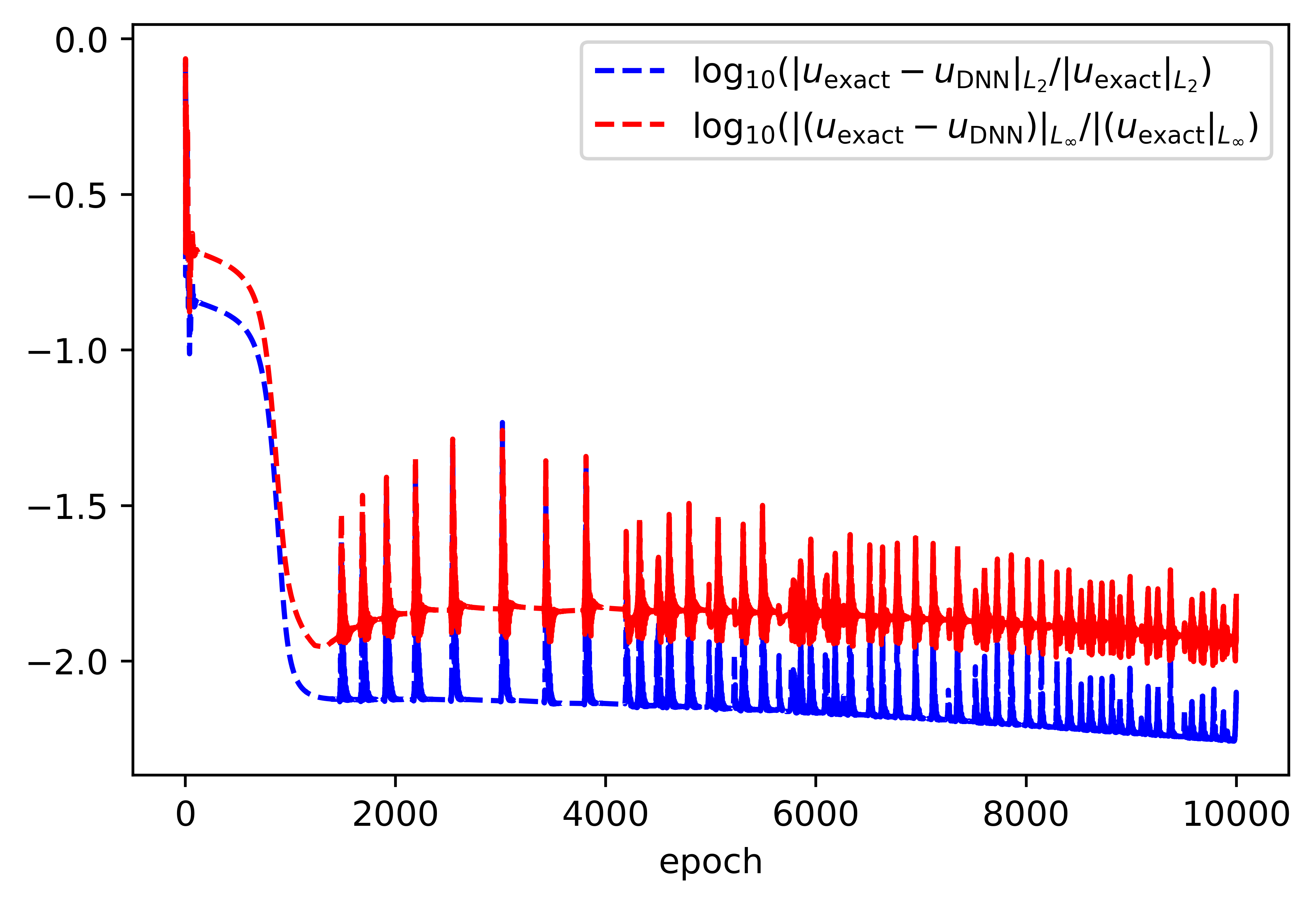}
			\caption{Numerical results for Example \ref{ex:2d_EPT} with $c=4, R=1$  (Relative $L^2$-error: $6.098\times 10^{-3}$; Relative $L^{\infty}$-error: $1.138\times 10^{-2}$). 
			}\label{fig:2d_elasto_plastic_4}
		\end{figure}

	To further validate the effectiveness of Algorithm \ref{alg:dl_evi} for solving elasto-plastic torsion problems, we compare the numerical results with the benchmark ones obtained by the ADMM method suggested in \cite{glowinski2015variational}. The ADMM decomposes the original problem into two simpler subproblems. One subproblem is to solve a convection equation and the other one requires to compute the projection onto $\{\bm{q}\in [L^2(\Omega)]^d\mid \|\bm{q}\|_{[L^2(\Omega)]^d}\leq 1\}$. To implement the ADMM, all the subproblems are discretized by a finite element method (FEM) with the mesh generated by the iFEM package \cite{Chen:2008ifem}.  An ADMM-FEM method is thus obtained. 
We test the ADMM-FEM on different grid resolutions $N$.
We train the neural network $\hat{u}(x; \bm{\theta}_u)$ using a set with fixed $10^2$ points randomly sampled from $\bar{\Omega}$ and then test $\hat{u}(x; \bm{\theta}_u)$ with different $N$.  We use the relative $L^2$-errors to evaluate and compare the numerical accuracy of the computed solutions. The numerical comparisons are reported in Table \ref{tab: compare_EP}.

\begin{table}[H]
	\centering
	\scalebox{0.9}{
		\begin{tabular}{|c|c|c|c|c|}
			\hline
			$N$   &$88$& $318$ & $1207$ & $4701$ \\ \hline
		ADMM-FEM \cite{glowinski2015variational}        &   $1.840\times 10^{-2}$                              &                  $1.119\times 10^{-2}$               &      $5.877\times 10^{-3}$                           &         $2.937\times 10^{-3}$                         \\ \hline
			Algorithm \ref{alg:dl_evi}              &        $4.782\times 10^{-3}$                          &           $3.857\times 10^{-3}$                          &                   $3.762\times 10^{-3}$               &              $3.669\times 10^{-3}$                             \\ \hline
	\end{tabular}}
	\caption{Comparisons with the ADMM-FEM \cite{glowinski2015variational} on different grid resolutions ($c=4$).}\label{tab: compare_EP}
\end{table}

The results in Table \ref{tab: compare_EP} demonstrate that when 
$N$ is small, the $L^2$-errors of solutions computed by Algorithm \ref{alg:dl_evi} are lower than those produced by the ADMM-FEM. Even as the grid resolution increases, Algorithm \ref{alg:dl_evi} remains competitive with the ADMM-FEM in accuracy. Notably, once the neural networks are trained on $10^2$ randomly sampled points, solve the problem for a new resolution requires only a forward pass of the pre-trained networks. In contrast, the ADMM-FEM must solve the problem from scratch for each resolution,  incurring significantly higher computational costs. These findings highlight the mesh-free nature and strong generalization capability of Algorithm \ref{alg:dl_evi}, establishing its effectiveness and numerical efficiency for solving elasto-plastic torsion problems.
		
\end{example}

\subsection{Bingham Visco-Plastic Flows}

We consider a visco-plastic medium of viscosity $\nu>0$ and plastic yield $\tau>0$ flowing in an infinitely long cylindrical pipe of bounded cross section $\Omega\subset \mathbb{R}^2$. Suppose that $\Omega$ is parallel to the horizontal plane, then in the steady state, the velocity of such a flow is given by $(0,0,u)$, where $u$ is characterized by
\begin{equation}\label{Bingham_VI}
	u\in H_0^1(\Omega),\quad\text{such that}\quad  \nu\int_\Omega \nabla u\cdot\nabla (v- u)dx+\tau\int_\Omega(|\nabla v|-|\nabla u|)dx\geq c\int_\Omega(v-u)dx,\forall v\in H_0^1(\Omega).
\end{equation}
The constant $c>0$ is the linear decay of pressure and
$\nu$, $\tau$ are, respectively, the viscosity and plasticity yield of the fluid. The above medium behaves like a viscous fluid (of viscosity $\nu$) in
$\Omega^+:=\{x\in \Omega\mid |\nabla u|>0\}$ and like a rigid medium in
$\Omega^0:=\{x\in \Omega\mid |\nabla u|=0\}$.   We refer to \cite{Glowinski1984Numerical,glowinski2008lectures,mosolov1966stagnant,mosolov1965variational} for a detailed study of the properties of \eqref{Bingham_VI}. A survey on the numerical methods for solving \eqref{Bingham_VI} can be found in \cite{Dean2007On}.

\begin{example}\label{ex:2d_VPF}
		We consider the two-dimensional problem with an exact solution given in \cite{Glowinski1984Numerical}.
	Let the domain $\Omega=\{x=(x_1, x_2) \mid |x|:=\sqrt{x_1^2+x_2^2}<R\}$. Let  $\displaystyle R^\prime = \frac{2\tau}{c}$  and then the exact solution \(u(x)\) is defined as follows:
	\begin{equation*}
		\text { if } c R\leq 2\tau, ~ u(x) = 0;\quad 	\text { if } cR>2\tau, ~ u(x) =
		\begin{cases}
			\begin{aligned}
				& \left(\frac{R-R^{\prime}}{2}\right)\left[\frac{c}{2}\left(R+R^{\prime}\right)-2 \tau\right], &0 \leq |x| \leq R^{\prime}, \\
				& \left(\frac{R-|x|}{2}\right)\left[\frac{c}{2}(R+|x|)-2 \tau\right] , & R^{\prime} \leq |x| \leq R.
			\end{aligned}
		\end{cases}
	\end{equation*}
	
 It is clear that the exact solution $u$ depends on the constants $R$, $c$, and $\tau$. In our numerical experiments, we set $R=1, c=10$ and take $\tau=1$ and $1.5$ to test Algorithm \ref{alg:dl_evi} with the loss function specified in \eqref{eq:loss_4_2}. To impose the boundary condition as a hard constraint, we take $h(x)=R^2-(x_1^2+x_2^2)$. The numerical results for this example with $\tau=1$ and $\tau=1.5$ are respectively reported in Figure \ref{fig:2d_viscous_plastic_4} and Figure \ref{fig:2d_viscous_plastic_10}. We observe that the numerical solutions are good approximations to the exact ones. In particular, for both cases, the low relative $L^2$-errors are of order $10^{-3}$, which indicates that Algorithm \ref{alg:dl_evi} can produce solutions with high accuracy for Bingham visco-plastic flows in different settings.

	\begin{figure}[h!]
		\centering
		\includegraphics[width=.32\textwidth]{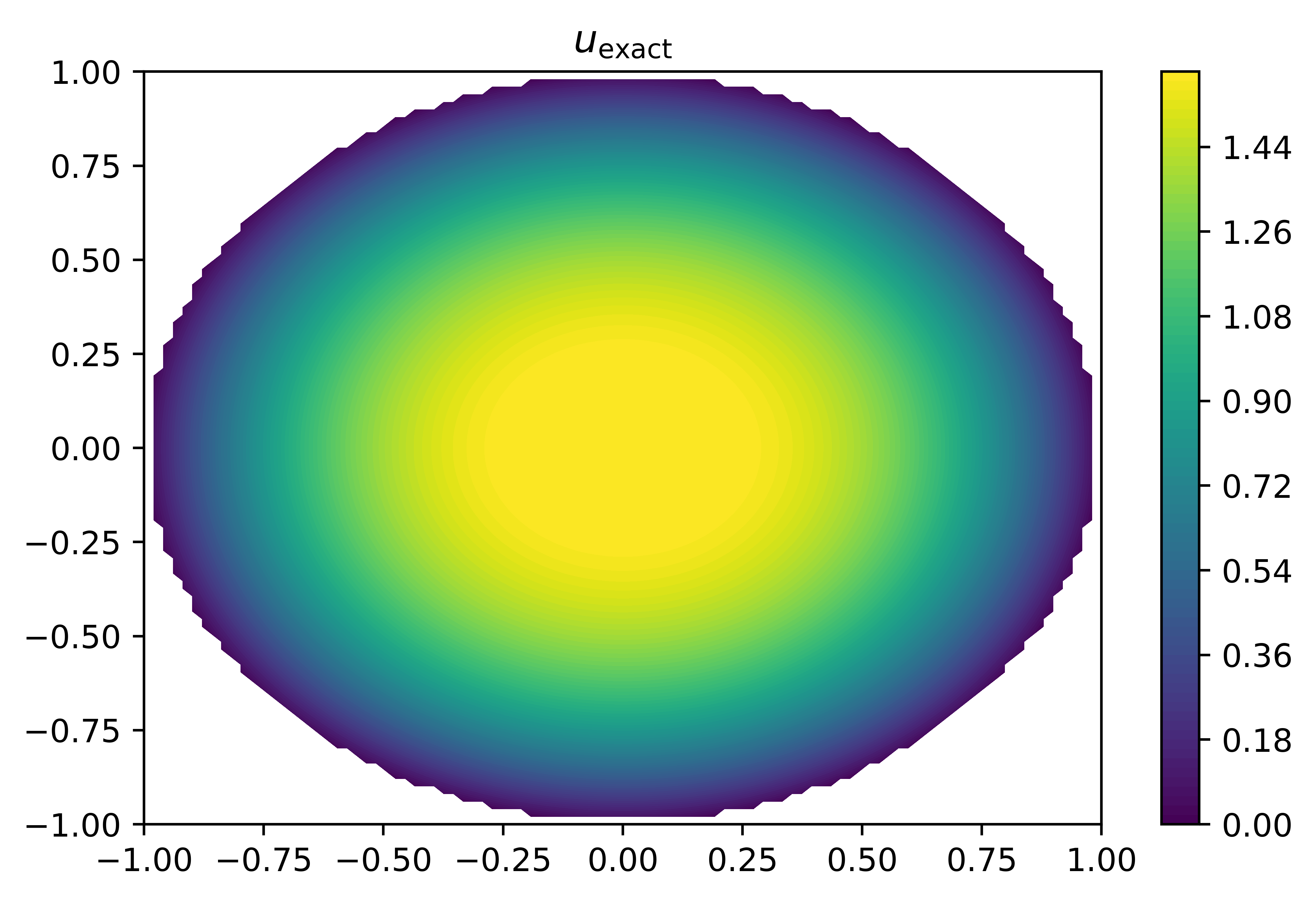}
		\includegraphics[width=.32\textwidth]{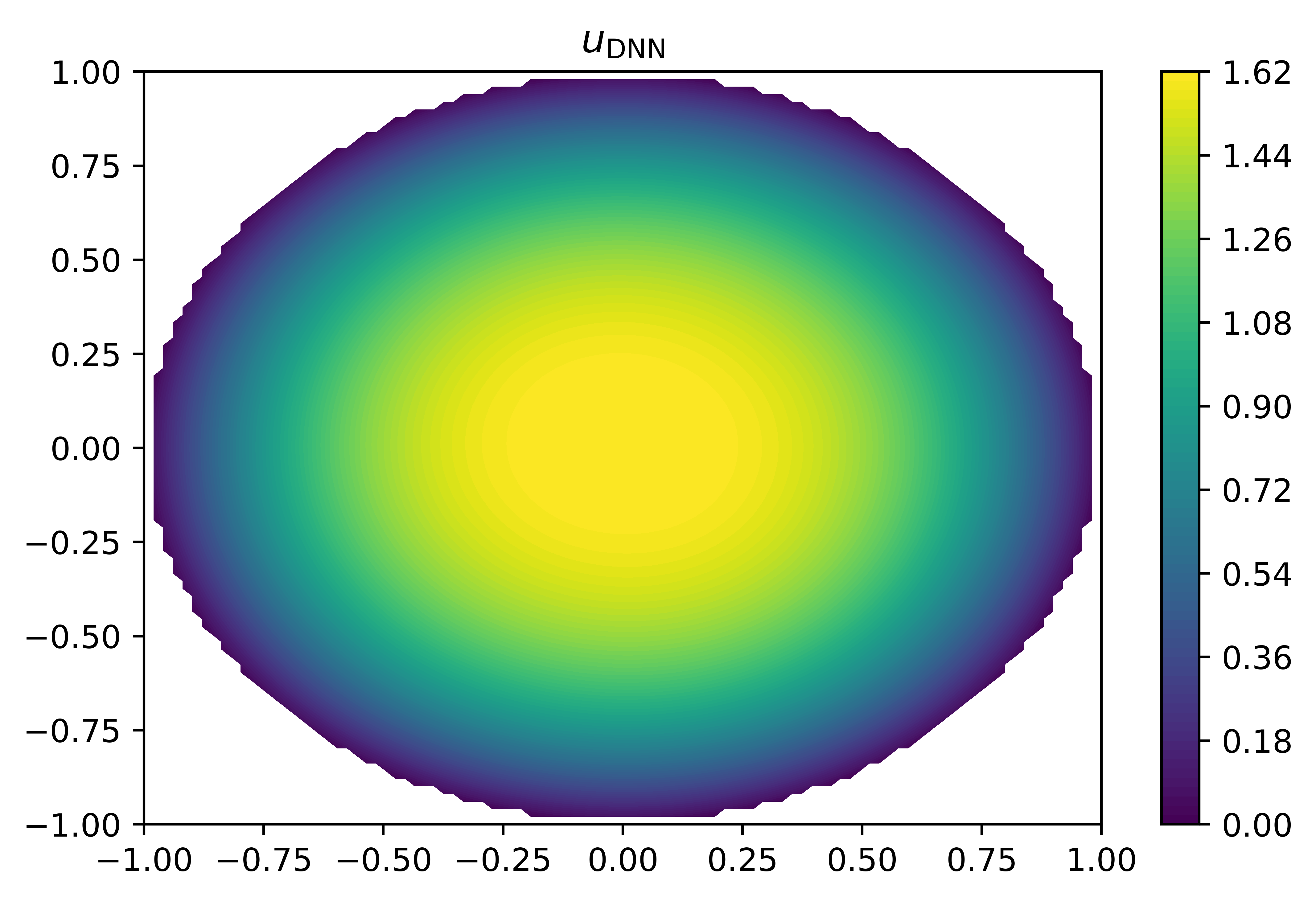}
		\includegraphics[width=.32\textwidth]{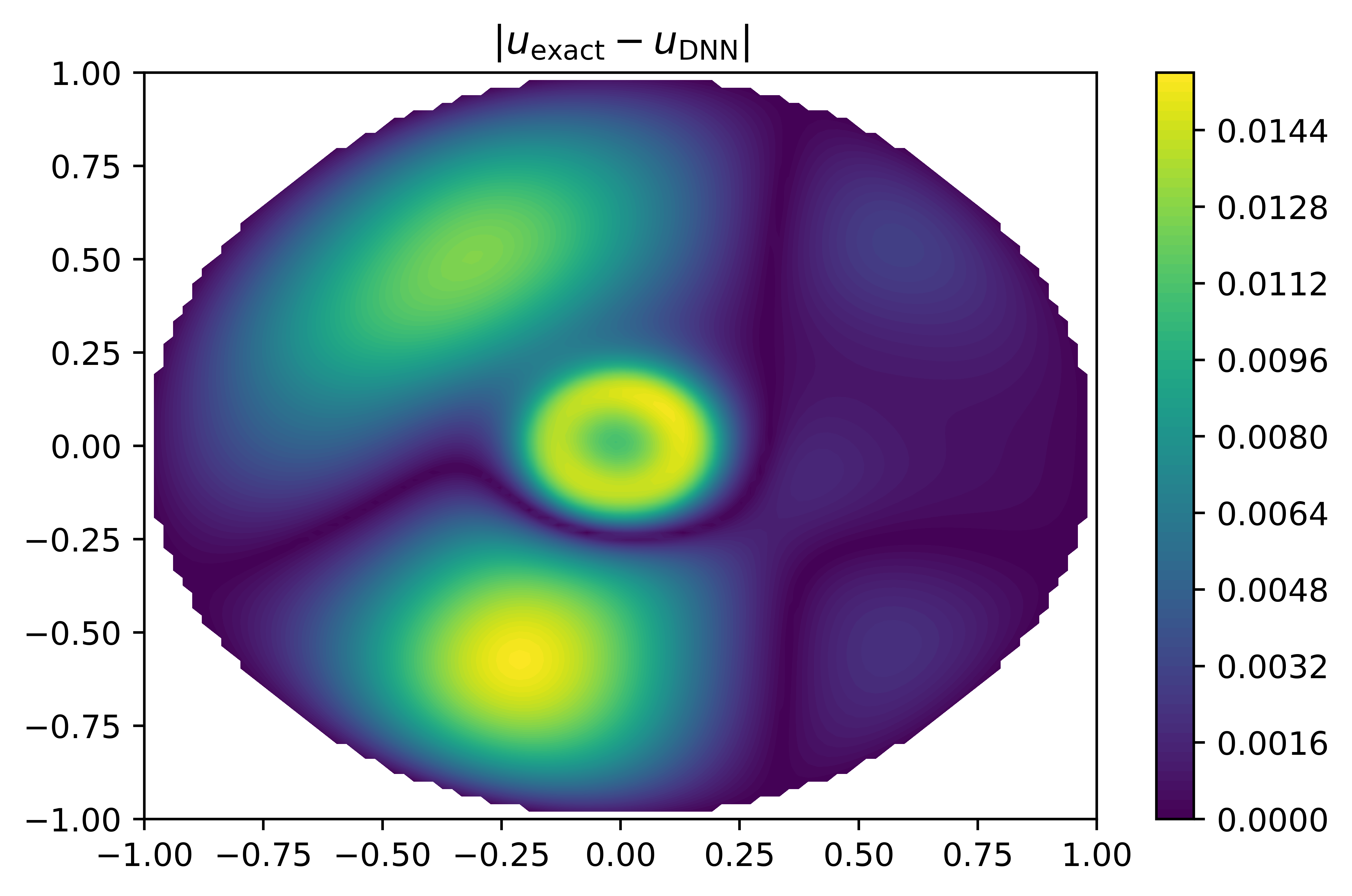}
		\includegraphics[width=.455\textwidth]{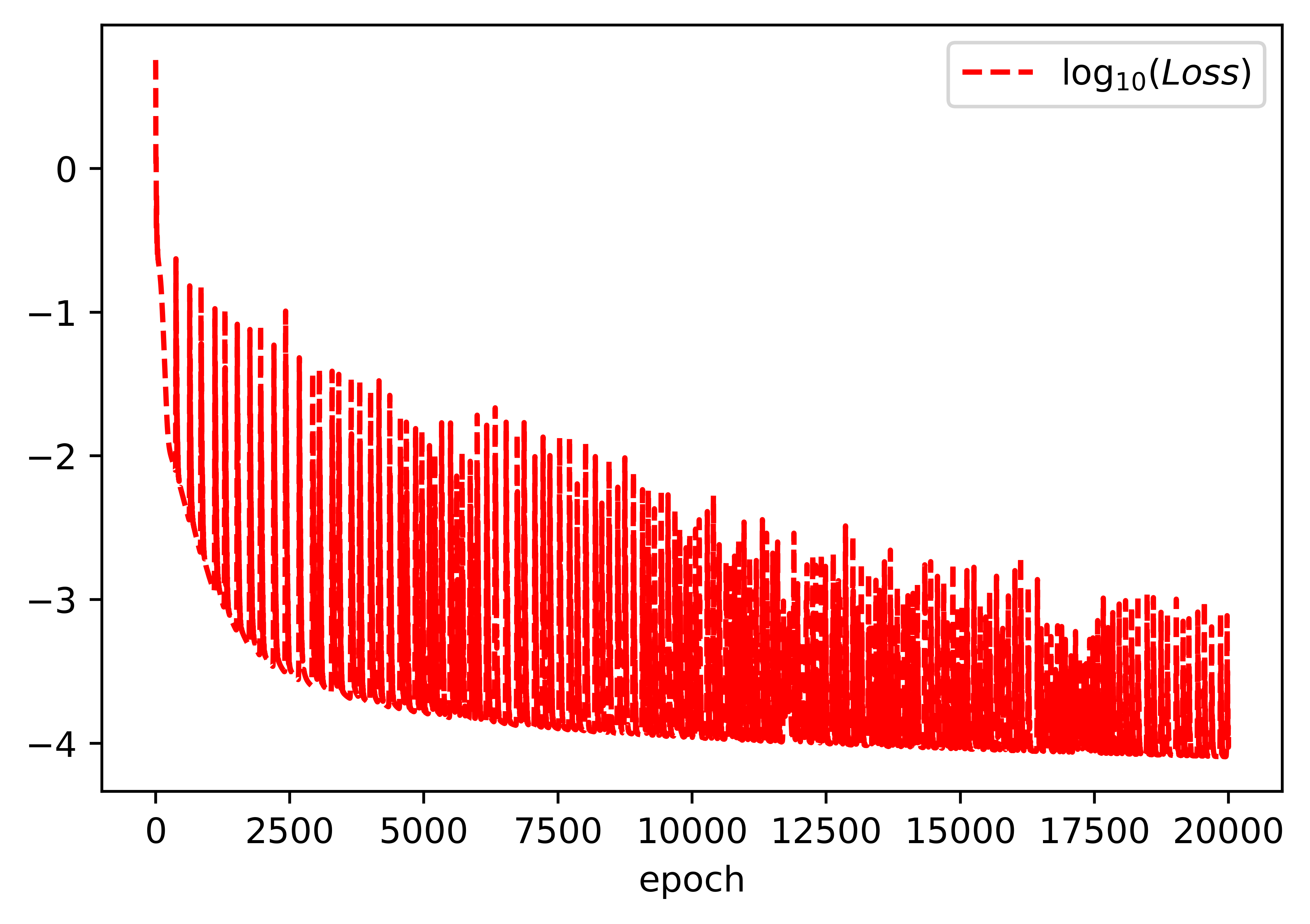}
		\includegraphics[width=.455\textwidth]{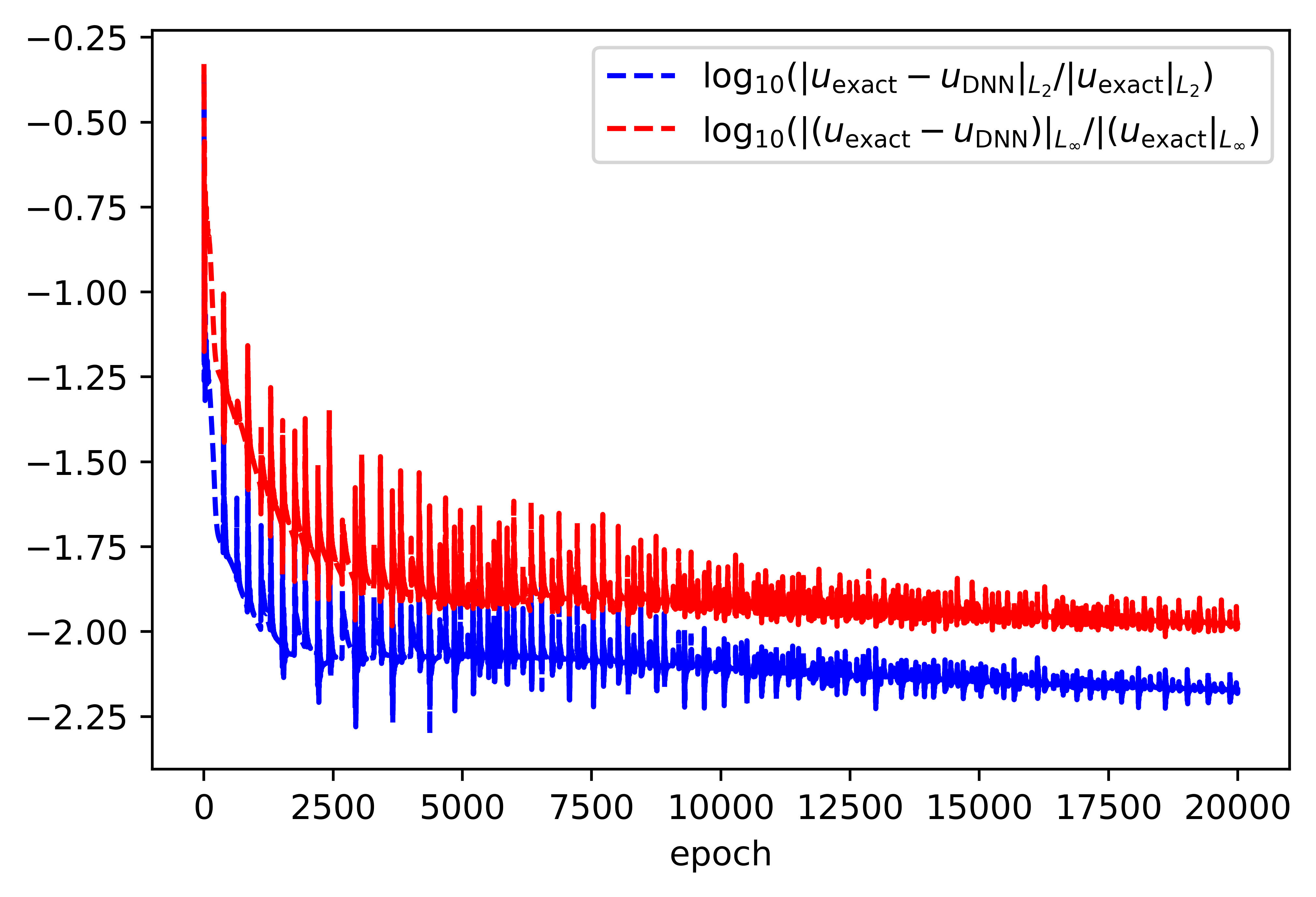}
		\caption{Numerical results for Example \ref{ex:2d_VPF} with $\tau=1$ and $c=10$ (Relative $L^2$-error: $6.819\times 10^{-3}$; Relative $L^{\infty}$-error: $1.060\times 10^{-2}$). 
		}\label{fig:2d_viscous_plastic_4}
	\end{figure}

	\begin{figure}[h!]
	\centering
	\includegraphics[width=.32\textwidth]{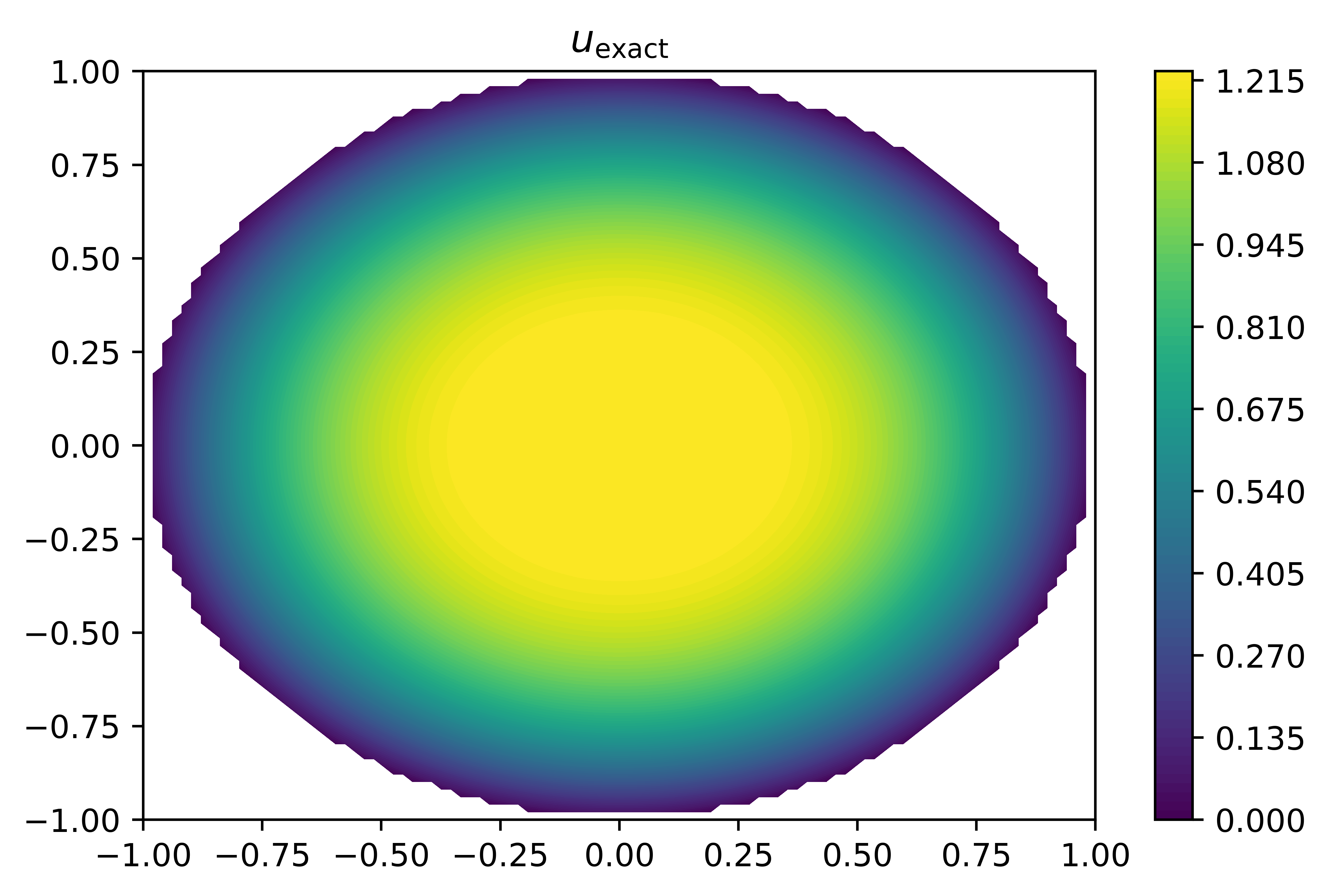}
	\includegraphics[width=.32\textwidth]{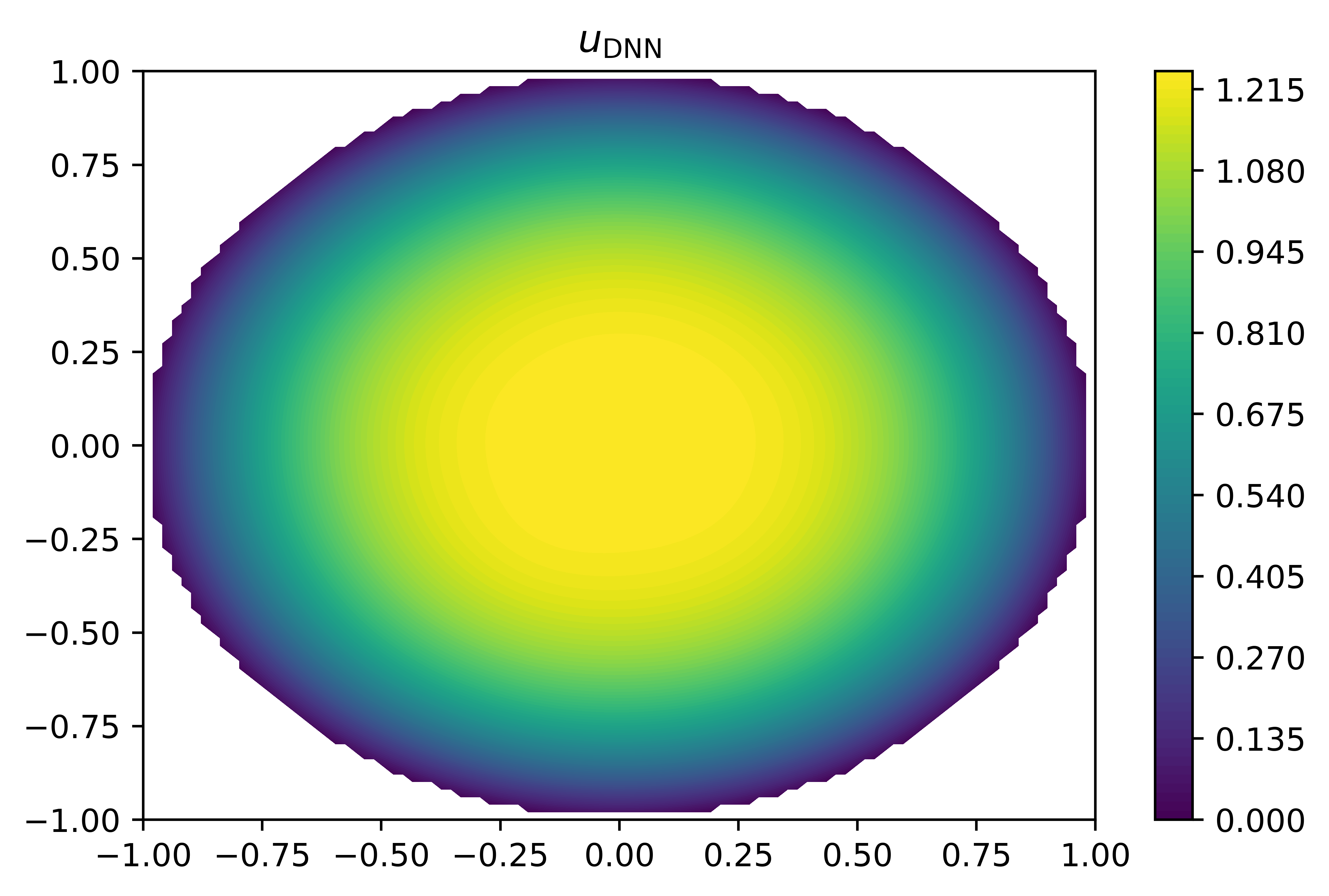}
	\includegraphics[width=.32\textwidth]{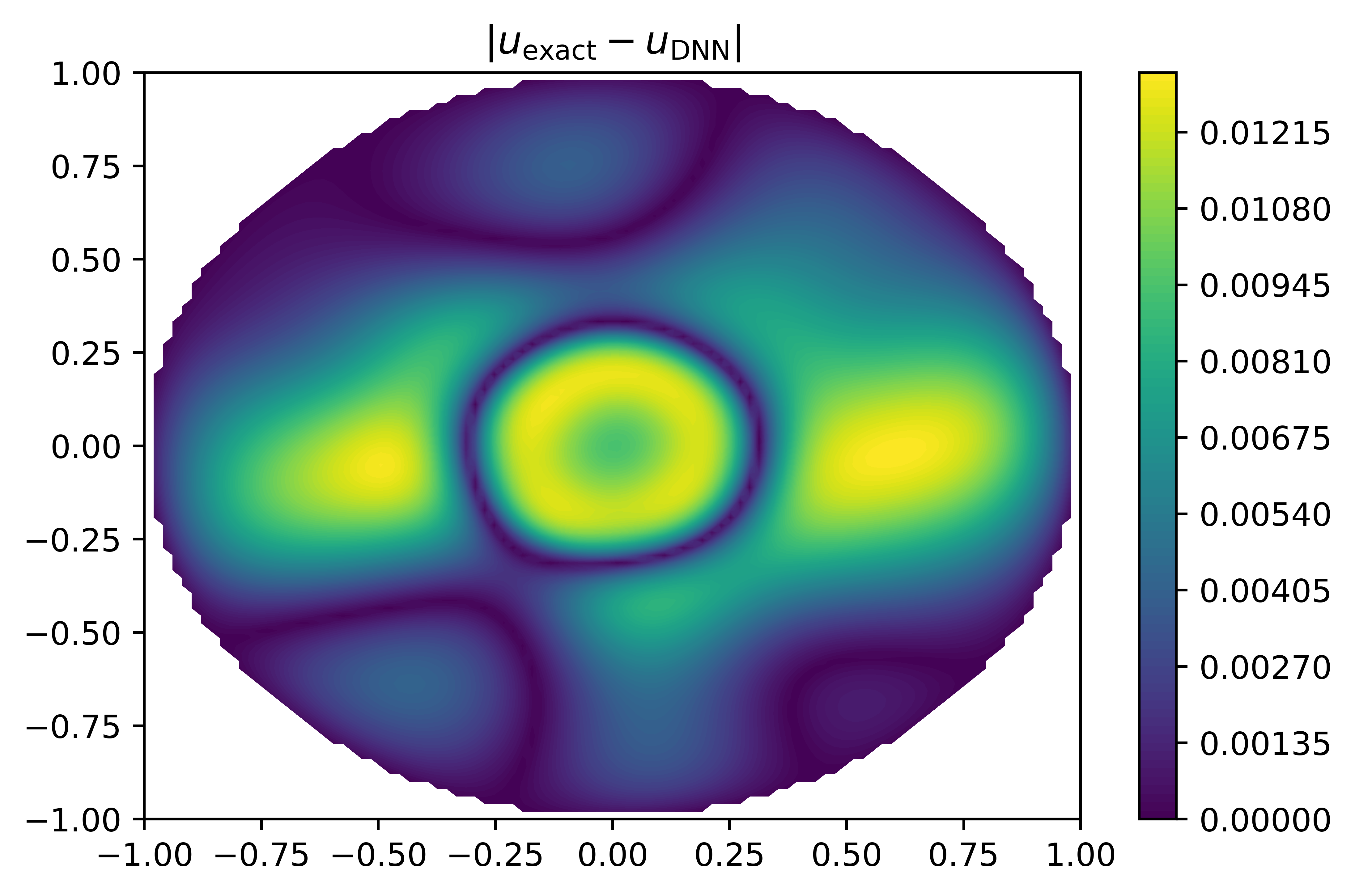}
	\includegraphics[width=.455\textwidth]{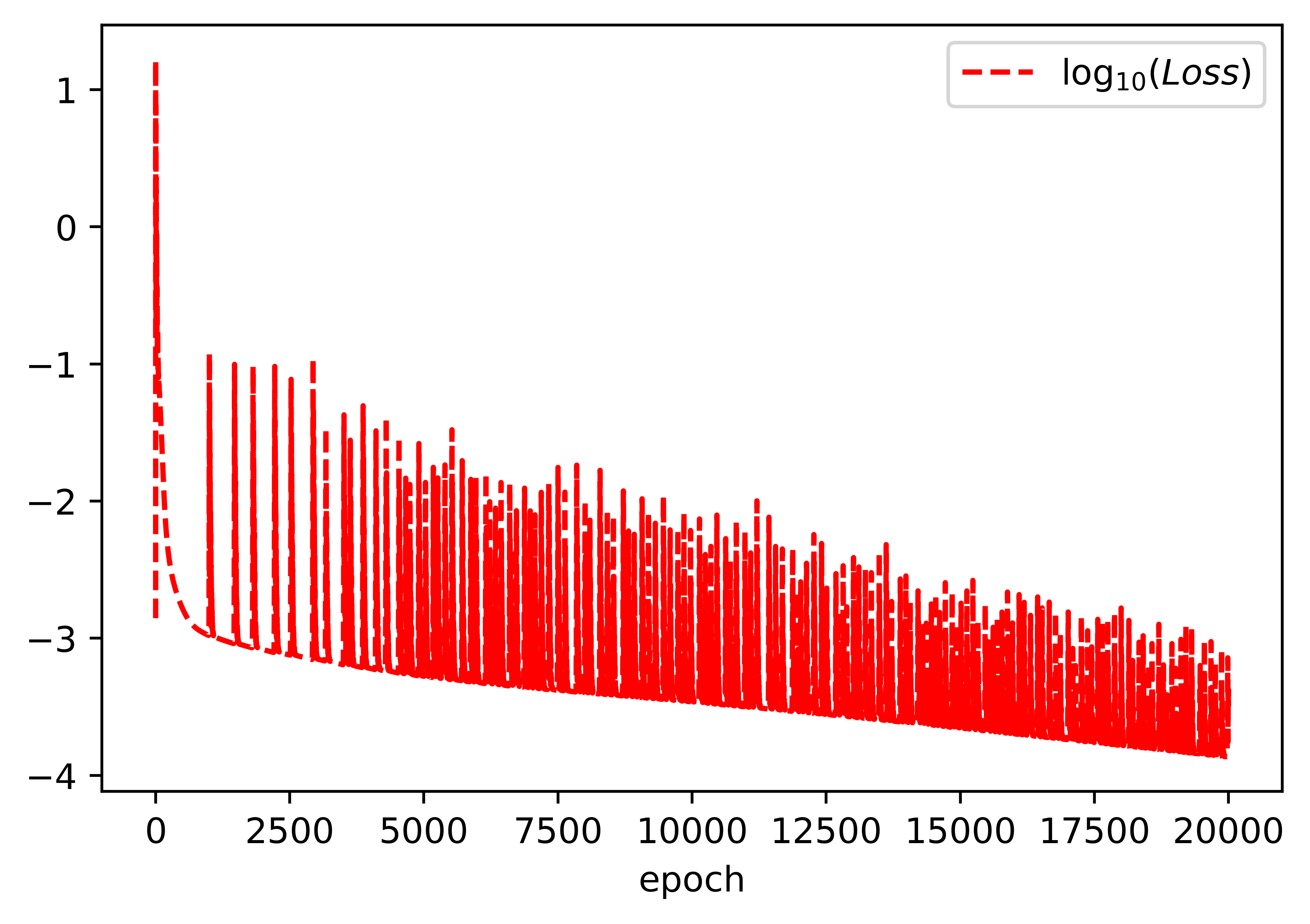}
	\includegraphics[width=.455\textwidth]{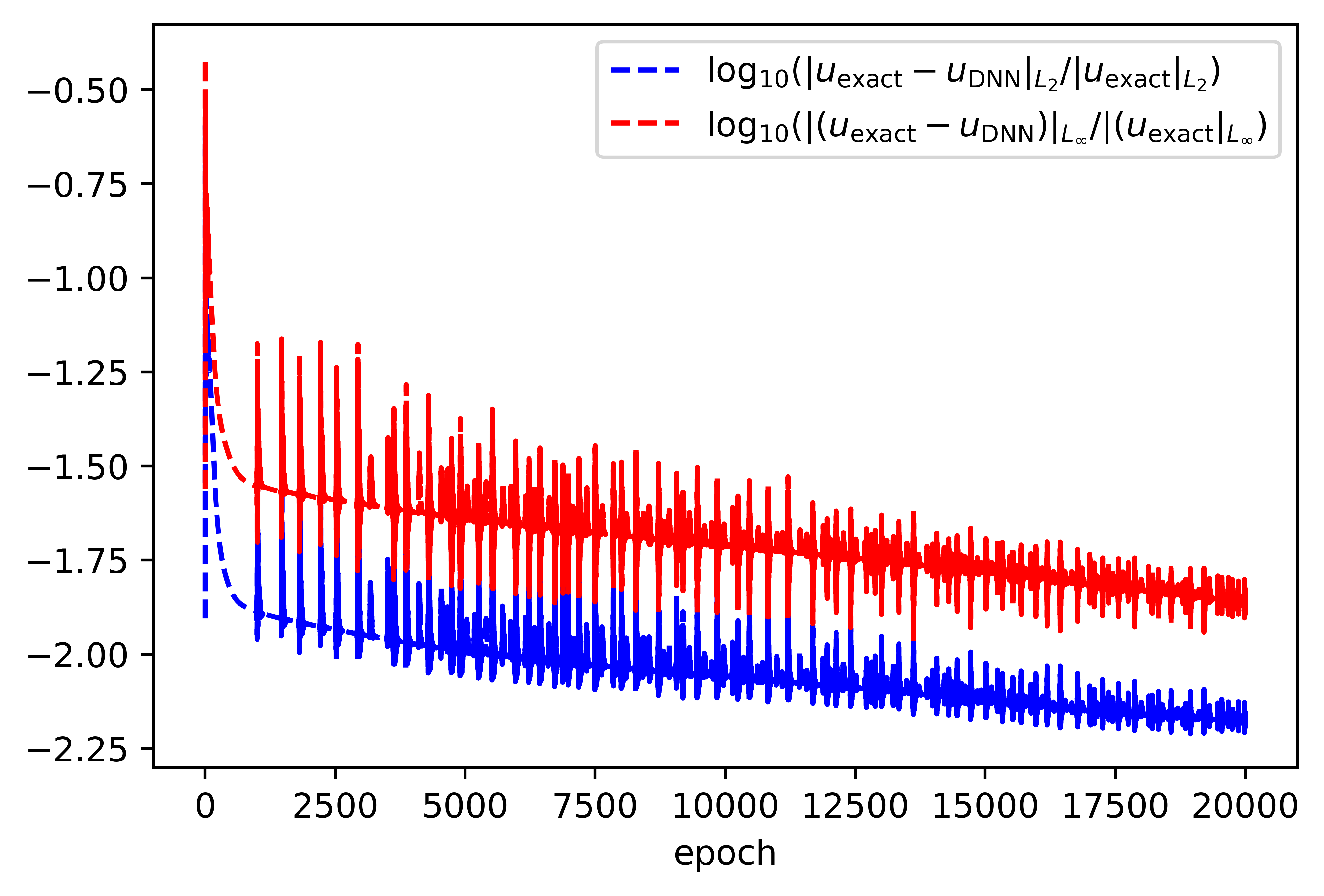}
	\caption{Numerical results for Example \ref{ex:2d_VPF} with $\tau=1.5$ and $c=10$ (Relative $L^2$-error: $6.350\times 10^{-3}$; Relative $L^{\infty}$-error: $1.267\times 10^{-2}$). 
		 }\label{fig:2d_viscous_plastic_10}
\end{figure}
	
\end{example}

\subsection{Simplified Friction Problems}\label{se: simplified_friction}
Friction phenomena between different bodied play an important role in structural and mechanical systems, see e.g, \cite{feng2019virtual,Glowinski1984Numerical,stadler2004semismooth}. Here, we consider the simplified friction problems \cite{Glowinski1984Numerical,glowinski2008lectures} that can be modeled by the EVI:
\begin{equation}\label{eq:simplified_friction_VI}
	u\in H^1_D(\Omega),~\text{such that}~  \int_\Omega Au(v- u)dx+\tau\int_{\Gamma_C}(|\gamma v|-|\gamma u|)dx\geq \int_\Omega f(v-u)dx,\forall v\in H^1_D(\Omega),
\end{equation}
where $\Omega$ is a bounded domain of $\mathbb{R}^d$ and $\partial\Omega$ is its boundary,  $Au=-\Delta u+u$, $H^1_D(\Omega)=\{v\in H^1(\Omega)\mid v=0~\text{on}~\Gamma_D\subset \partial\Omega\}$, $\Gamma_C=\partial\Omega/ \Gamma_D$, $\tau>0$,  and the trace operator $\gamma$ is defined by $\gamma v=v|_{\partial\Omega}$.

Let $j(v)=\tau\int_{\Gamma_C}|\gamma v|dx$. Then, following the similar arguments to those in Sections \ref{se: algorithmic_framework} and \ref{se:case_studies}, one can easily show that the solution $u\in H^1_D(\Omega)$ satisfies
$$\text{Prox}_{\eta j}((I-\eta A)u+\eta f)=u, ~\text{with}~\eta>0,$$
which, after introducing $\lambda^*\in L^2(\Gamma_C)$, can be reformulated as
\begin{equation*}
	\left\{
	\begin{aligned}
		&{\lambda^*}(x) \gamma u(x)=-\tau|\gamma u(x)| ~\text{on}~\Gamma_C,\\
		&|{\lambda^*}(x)|\leq \tau ~\left(i.e.~{\lambda^*}(x)=\frac{\tau{\lambda^*}(x)}{\max\{\tau,|{\lambda^*}(x)|\}}\right) ~\text{on}~\Gamma_C,	\\
		&Au=f~\text{in}~\Omega,~u=0~\text{on}~\Gamma_D,~ \frac{\partial u}{\partial n}-\lambda^*=0~\text{on}~\Gamma_C.\\
	\end{aligned}
	\right.
\end{equation*}

We construct neural networks $\hat{u}(x; \theta_u) = h(x) \mathcal{N}_u(x; \bm{\theta}_u)$ and $\hat{\bm{\lambda}}(x;\bm{\theta}_\lambda)=\frac{\tau{\mathcal{N}_{\lambda}(x; \bm{\theta}_u)}}{\max\{\tau,|\mathcal{N}_{\lambda}(x; \bm{\theta}_u)|\}}$ to respectively approximate $u$ and $\lambda$, where the function $h: \bar{\Omega}\rightarrow \mathbb{R}$ verifies $h(x)=0$ if and only if  $x\in \Gamma_D$, and $\mathcal{N}_u(x; \bm{\theta}_u)$ and $\mathcal{N}_\lambda(x; \bm{\theta}_\lambda)$ are neural networks parameterized by $\bm{\theta}_u$ and $\bm{\theta}_\lambda$, respectively.  We then implement Algorithm \ref{alg:dl_evi} to \eqref{eq:simplified_friction_VI} with the following loss function 
	\begin{equation*}
	\begin{aligned}
		\mathcal{L}(\bm{\theta}_u,\bm{\theta}_{\lambda})=\frac{1}{|\mathcal{T}_C|}\sum_{x\in \mathcal{T}_C}\Bigg\{
		w_1	\left|\hat{\bm{\lambda}}(x;\bm{\theta}_{\lambda}) \hat{u}(x;\bm{\theta}_u)+\tau| \hat{u}(x;\bm{\theta}_u)|\right|^2
		+w_2\left|\frac{\partial \hat{u}(x;\bm{\theta}_u)}{\partial n}-\hat{\bm{\lambda}}(x;\bm{\theta}_{\lambda})\right|^2\Bigg\}\\
		+w_3\frac{1}{|\mathcal{T}|}\sum_{x\in \mathcal{T}}\left|A\hat{u}(x;\bm{\theta}_u)- f(x)\right|^2,
	\end{aligned}
\end{equation*}
where $w_i>0, i=1,2,3$, are the weights, $\mathcal{T}\subset\Omega$ and $\mathcal{T}_C\subset\Gamma_C$ are sampled training sets. 

\begin{example}\label{ex:2d_SF}
	We consider an example that has been studied in \cite{feng2019virtual,huang2020int}. In particular, we let $\Omega=(0,1)\times (0,1)$, $\tau=1$, $\Gamma_C=\{1\}\times [0,1]$, and $\Gamma_D=\partial\Omega\backslash\Gamma_C$. The exact solution $u$ is given by 
	$
	u(x_1, x_2)=(\sin x_1-x_1\sin 1)\sin 2\pi x_2,
	$
	and the source term is
	$
	f(x_1,x_2)=((2+4\pi^2)\sin x_1-(1+4\pi^2)x_1\sin 1)\sin 2\pi x_2.
	$
	
	To impose the boundary condition $u=0$ on $\Gamma_D$ as a hard constraint, we take $h(x)=4x_1x_2(1-x_2)$.  The numerical results of Algorithm \ref{alg:dl_evi} for this example are presented in Figure \ref{fig:2d_simplified_friction}.  We observe that the numerical solutions are good approximations to the exact ones. In particular, the maximal point-wise error between the learned and exact solutions is of order $10^{-5}$, which, together with the low relative $L^2$-error $3.616\times 10^{-4}$ and $L^{\infty}$-error $5.281\times 10^{-4}$, validates that Algorithm \ref{alg:dl_evi} can produce a high-accurate solution for the simplified friction problem under investigation. 
	
		\begin{figure}[h!]
		\centering
		\includegraphics[width=.32\textwidth]{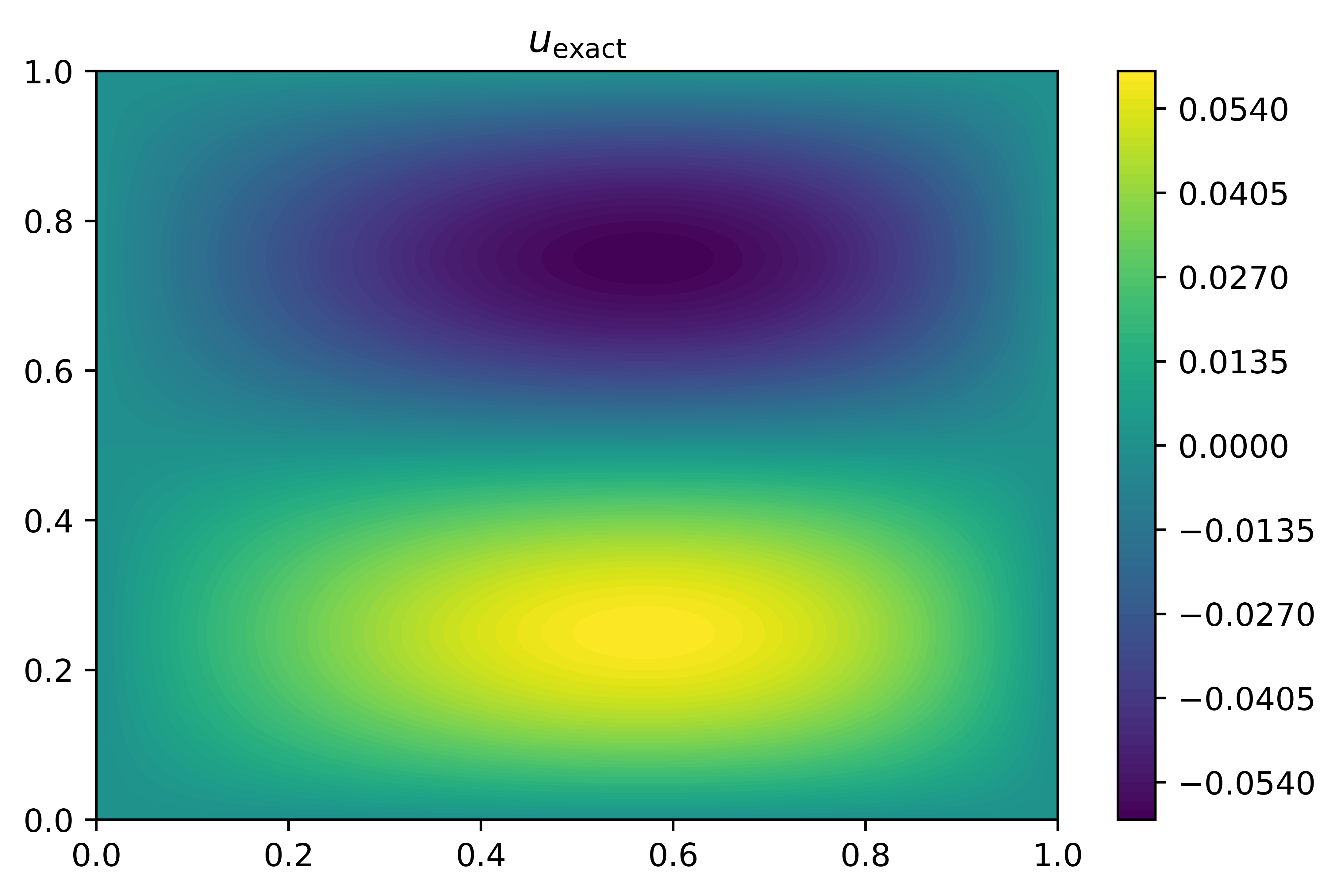}
		\includegraphics[width=.32\textwidth]{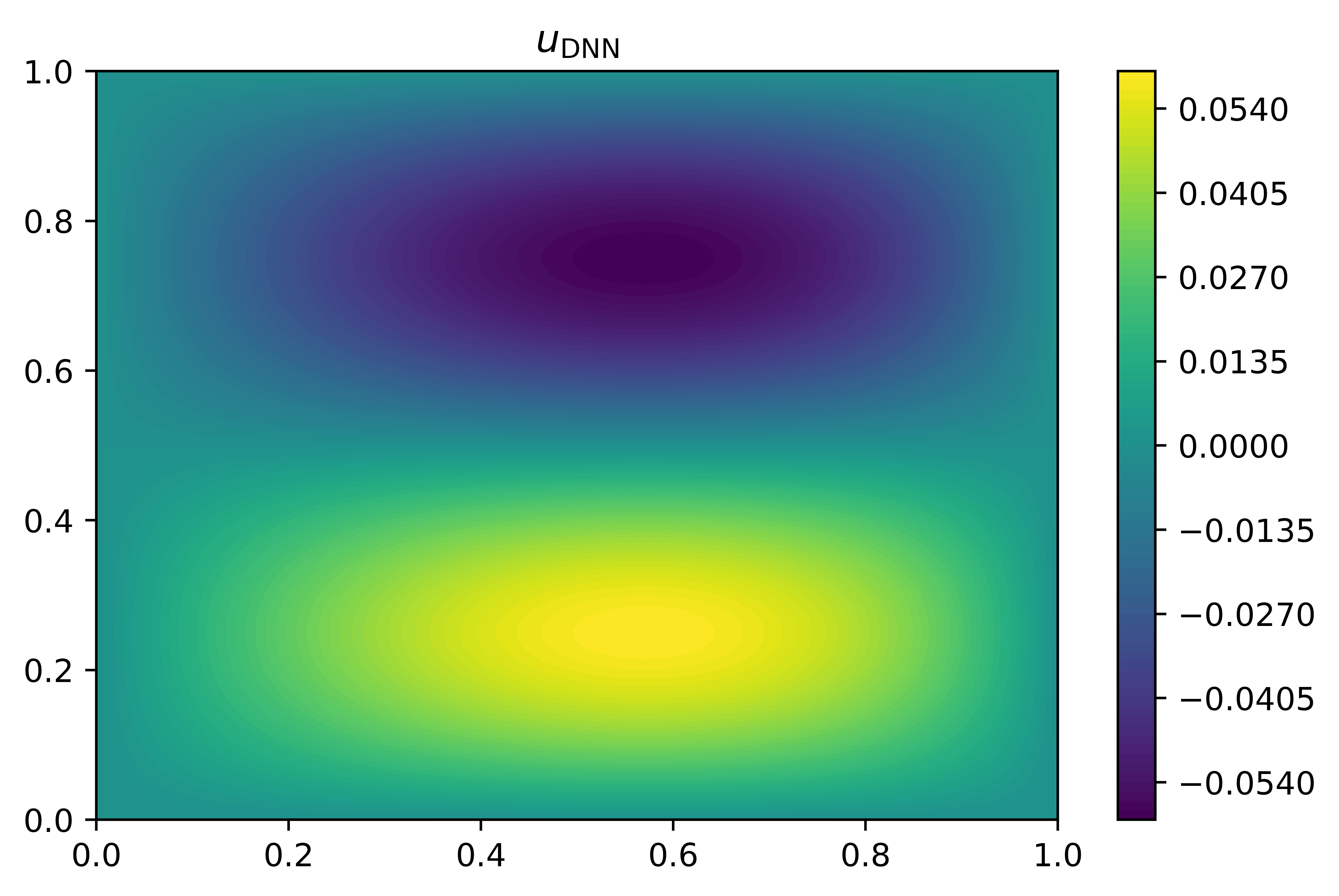}
		\includegraphics[width=.32\textwidth]{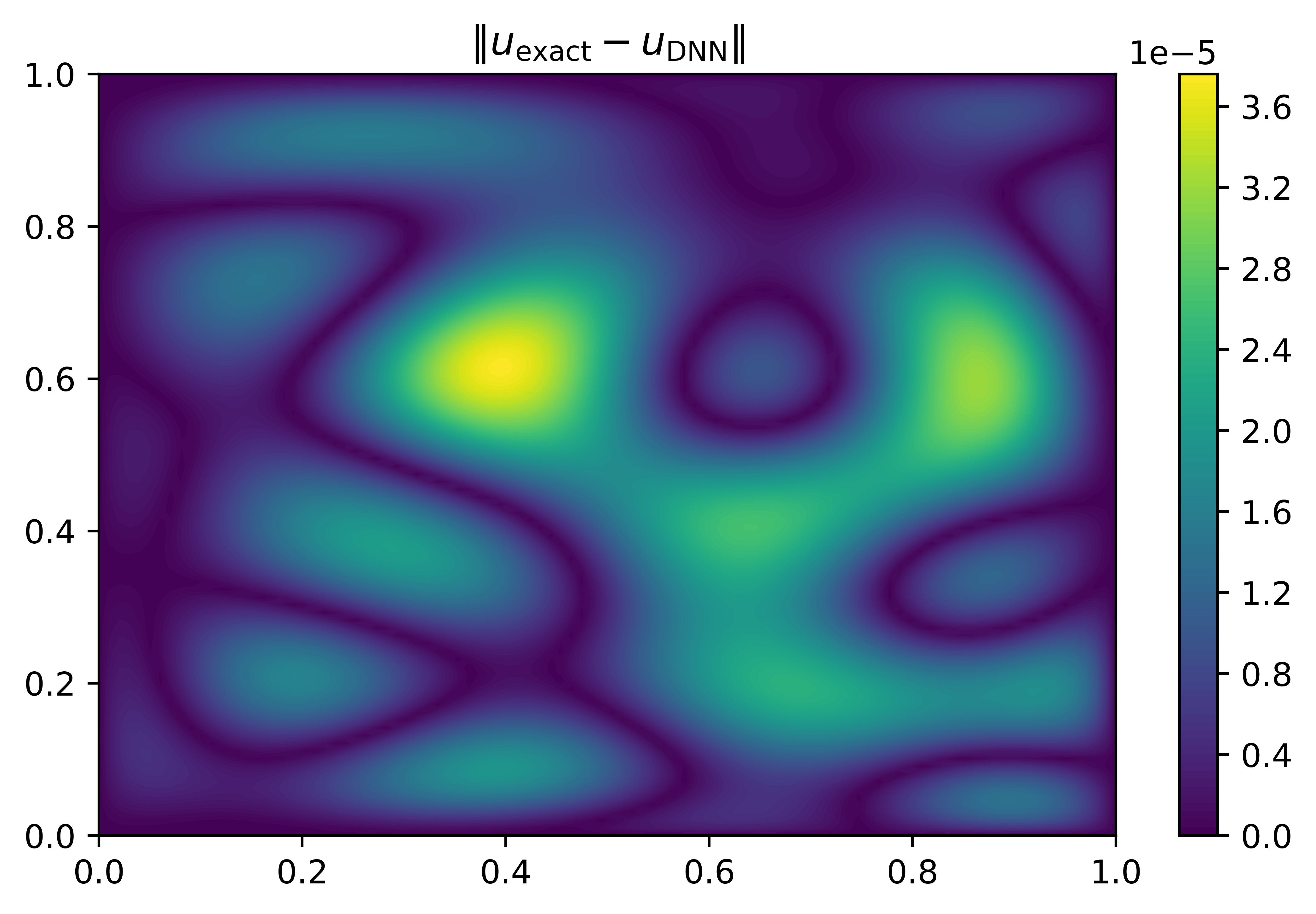}
		\includegraphics[width=.455\textwidth]{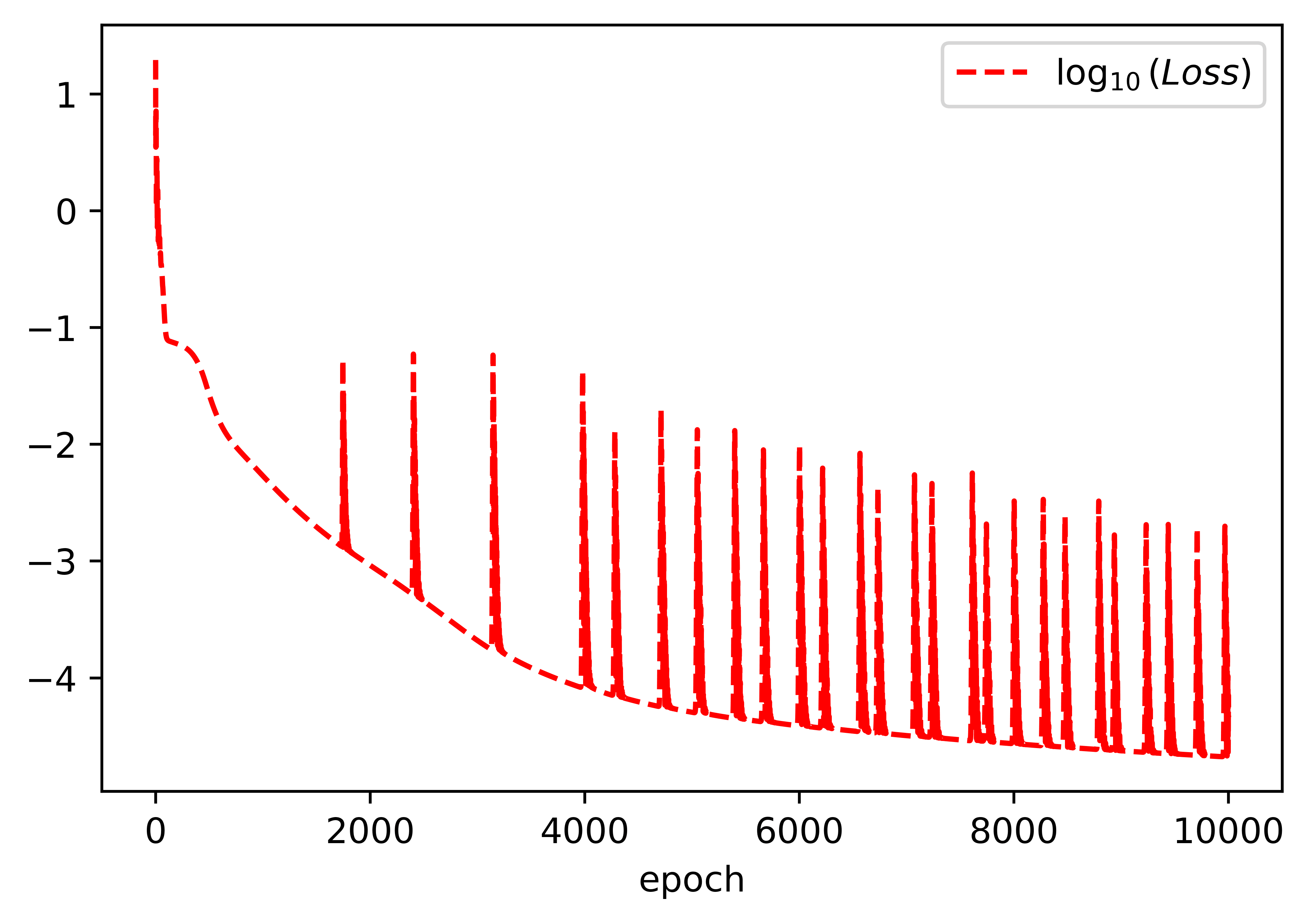}
		\includegraphics[width=.455\textwidth]{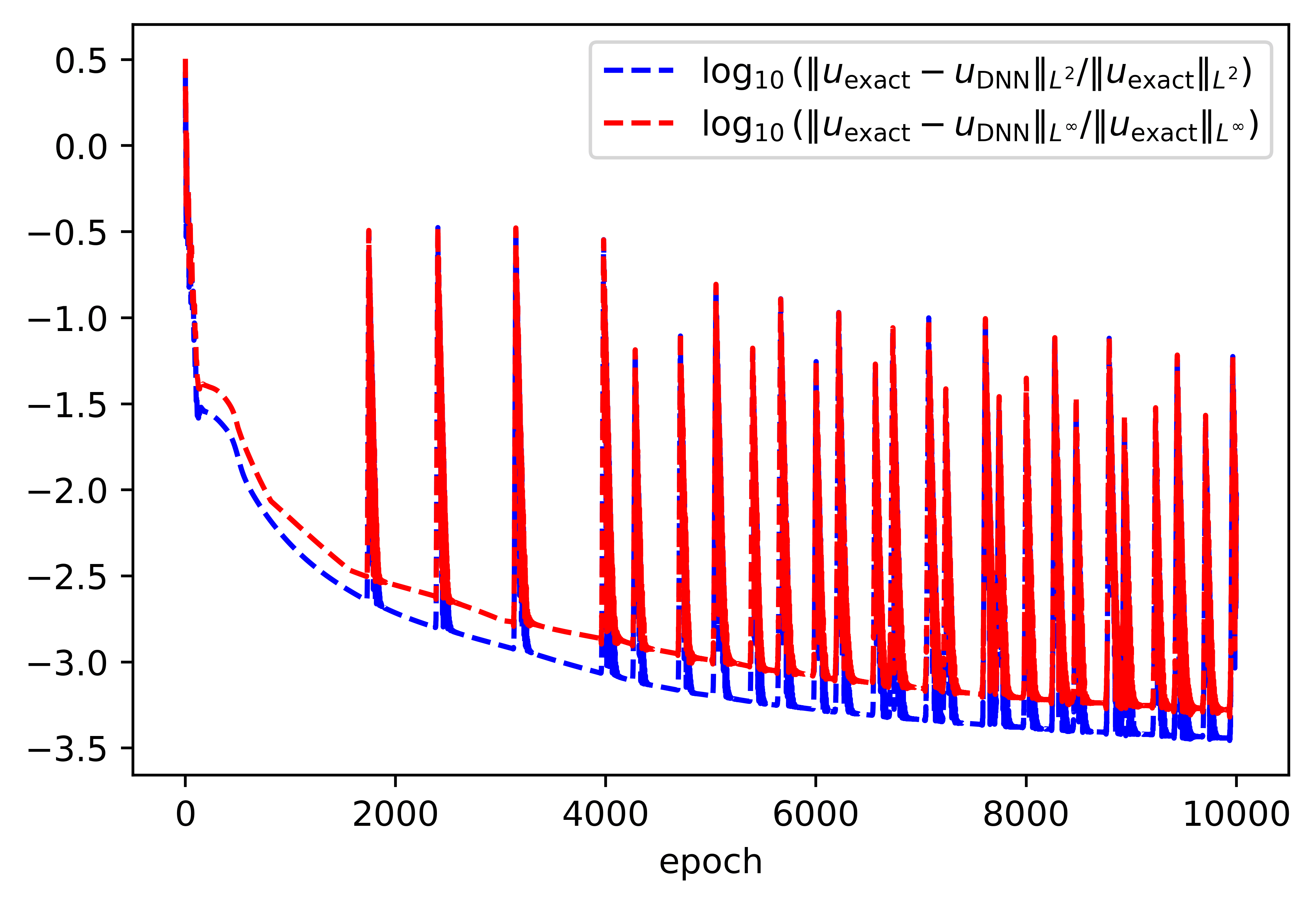}
		\caption{Numerical results for Example \ref{ex:2d_SF}  (Relative $L^2$-error: $3.616\times 10^{-4}$; Relative $L^{\infty}$-error: $5.281\times 10^{-4}$). 
			 }\label{fig:2d_simplified_friction}
	\end{figure}
	
	To further validate the effectiveness of Algorithm \ref{alg:dl_evi}, we compare it with the virtual finite element method in \cite{feng2019virtual}, which is a benchmark mesh-based traditional numerical algorithm for solving simplified friction problems.
	We train the neural network $\hat{u}(x; \bm{\theta}_u)$ using $10^3$ points randomly sampled from $\overline{\Omega}$ and then test $\hat{u}(x; \bm{\theta}_u)$ with different grid resolutions.  We use the absolute $L^\infty$-error used in \cite{feng2019virtual} to evaluate and compare the numerical accuracy of the computed solutions.
	Following \cite{feng2019virtual}, we evaluate the absolute $L^\infty$-errors on an $N \times N$ uniform grid over $\overline{\Omega}$ with $N = 8$, $16$, $32$, $64$, and $128$. The numerical comparisons are reported in Table \ref{tab: compare_SF}. 
\begin{table}[H]
			\centering
			\scalebox{0.9}{
	\begin{tabular}{|c|c|c|c|c|c|}
		\hline
		$N$   &8&${16}$& ${32}$ & ${64}$ & ${128}$ \\ \hline
		FEM \cite{feng2019virtual}      &     $2.70\times 10^{-3}$                           &   $4.22\times 10^{-4}$                              &                  $1.47\times 10^{-4}$               &      $4.66\times 10^{-5}$                           &         $1.20\times 10^{-5}$                         \\ \hline
		Algorithm \ref{alg:dl_evi}&           $2.41\times 10^{-5}$                     &        $3.42\times 10^{-5}$                          &           $3.22\times 10^{-5}$                          &                   $3.20\times 10^{-5}$               &              $3.35\times 10^{-5}$                             \\ \hline
	\end{tabular}}
	\caption{Comparison with the FEM \cite{feng2019virtual} on different grid resolutions.}\label{tab: compare_SF}
\end{table}

From the results in Table \ref{tab: compare_SF}, we can see that when $N \leq 32$, the $L^\infty$-errors of the computed solutions by Algorithm \ref{alg:dl_evi} are significantly lower than those by the FEM.
Even if the mesh resolution increases to $N = 128$,  Algorithm \ref{alg:dl_evi} is still comparable with the FEM. Moreover, note that after training the neural networks with $10^3$ randomly sampled points, the evaluation of  Algorithm \ref{alg:dl_evi} for a new resolution requires only a forward pass of these neural networks. In contrast, for each resolution, the FEM requires solving the simplified friction problem from scratch, which is more computationally expensive. 
These results validate the mesh-free nature and the generalization ability of Algorithm \ref{alg:dl_evi}, making it effective and numerically favorable for simplified friction problems.

\end{example}


\section{Conclusions and Perspectives}\label{se:conclusion}
This work presents the Prox-PINNs, a deep learning algorithmic framework that combines proximal operators and physics-informed neural networks (PINNs), for solving elliptic variational inequalities (EVIs). The Prox-PINNs framework reformulates EVIs as nonlinear equations through proximal operators, which are subsequently solved using hard-constraint PINNs. The Prox-PINNs framework is adaptable to various EVIs by leveraging analytical proximal operators for specific nonsmooth functionals. It thus alleviates the limitations of traditional mesh-based approaches and existing deep learning methods, which often lack generality or impose restrictive assumptions on the EVIs under investigation. The Prox-PINNs framework can be
used to develop efficient deep learning algorithms for diverse EVIs, including obstacle problems, elasto-plastic torsion problems, Bingham flows, and simplified friction problems. Numerical results show the framework's effectiveness, efficiency, accuracy, and robustness, even for problems non-symmetric operators and piecewise smooth solutions. 

The novelty of the Prox-PINNs framework opens up several possibilities for future investigation.

\begin{itemize}

	\item \textbf{Theoretical foundations}: The empirical success of Prox-PINNs motivates further investigation into their theoretical underpinnings. Rigorous analysis of convergence properties, stability, and error estimation would strengthen the mathematical justification of the framework. 
	
	\item \textbf{Algorithmic enhancements}: Integrating adaptive sampling strategies (e.g., \cite{gao2023failure,wu2023comprehensive}) and advanced optimization techniques for training (e.g., \cite{kiyani2025which,muller2023achieving}) with Prox-PINNs promises to enhance computational efficiency and solution accuracy.
	
	\item \textbf{Uncertainty quantification (UQ)}: Ensuring reliability in real-world applications requires robust methods to quantify uncertainties associated with data noise, model hyperparameters, and numerical approximations. Recent advances in UQ for PINNs \cite{psaros2023uncertainty,zhang2019quantifying,zou2025uncertainty,zou2024neuraluq} offer a foundation for adapting these techniques to Prox-PINNs.  
	
	\item \textbf{Extensions}: Expanding the Prox-PINNs framework to more challenging classes of VIs, such as stochastic EVIs \cite{kornhuber2014multilevel} and parabolic VIs, could broaden its applicability. Of particular interest are problems in computational finance, including American options pricing \cite{gao2020primal,jaillet1990variational}.
\end{itemize}

\section*{Acknowledgments}
The work of Y. Song was supported by a start-up grant from Nanyang Technological University.
The work of Z. Tan was supported in part by the National Key R\&D Program of China (Grant Number: 2024YFA1012503), the National Natural Science Foundation of China (Grant No. 12401542) and the Fundamental Research Funds for the Central Universities (Grant No. 20720240133).The work of H. Yue was supported by the National Natural Science Foundation of China (Grant No. 12301399).


%
%
%
\end{document}